\documentclass{article}

\usepackage{arxiv}

\usepackage[utf8]{inputenc} % allow utf-8 input
\usepackage[T1]{fontenc}    % use 8-bit T1 fonts
\usepackage{hyperref}       % hyperlinks
\usepackage{url}            % simple URL typesetting
\usepackage{booktabs}       % professional-quality tables
\usepackage{amsfonts}       % blackboard math symbols
\usepackage{nicefrac}       % compact symbols for 1/2, etc.
\usepackage{microtype}      % microtypography
\usepackage{lipsum}		% Can be removed after putting your text content
\usepackage{graphicx}
\usepackage{natbib}
\usepackage{doi}
\usepackage{multirow}
\usepackage{subfigure}
\usepackage{amsmath}
\newtheorem{definition}{Definition}

\title{An efficient class of increasingly high-order ENO schemes with multi-resolution}

\author{Hua Shen \thanks{This work was supported by the National Natural Science Foundation of China (Contract No. 11901602).} \\
        School of Mathematical Sciences\\
        University of Electronic Science and Technology of China\\
        Chengdu, Sichuan, 611731, China \\
	\texttt{huashen@uestc.edu.cn} \\
}

\renewcommand{\shorttitle}{\textit{arXiv} Template}

\hypersetup{
pdftitle={An efficient class of increasingly high-order ENO schemes with multi-resolution},
pdfauthor={Hua Shen},
pdfkeywords={high-order schemes, ENO,  WENO, finite difference scheme, hyperbolic conservation laws},
}

\begin{document}
\maketitle

\begin{abstract}
  We construct an efficient class of increasingly high-order (up to 17th-order) essentially non-oscillatory schemes 
  with multi-resolution (ENO-MR) for solving hyperbolic conservation laws. 
  The candidate stencils for constructing ENO-MR schemes range from first-order one-point stencil
  increasingly up to the designed very high-order stencil. 
  The proposed ENO-MR schemes adopt a very simple and efficient strategy
  that only requires the computation of the highest-order derivatives of a part of candidate stencils.
  Besides simplicity and high efficiency, ENO-MR schemes are completely parameter-free and essentially scale-invariant.
  Theoretical analysis and numerical computations show that ENO-MR schemes achieve designed high-order convergence in smooth regions
  which may contain high-order critical points (local extrema) and retain ENO property for strong shocks.
  In addition, ENO-MR schemes could capture complex flow structures very well.
\end{abstract}

% keywords can be removed

%%%%%%%%%%%%%%%%%%%%%%%%%%%%%%%%%%%%%%%%%%%%%%%%%%%%%%%%%%%%%%%%%%%%%%%%%%%%%%%%%%%%%%%%%%%%%%%%%%%%%%%%%%%%%%%%%%%

\section{Introduction}\label{sec:intro}
Hyperbolic conservation laws are ubiquitous in nature.
Although the corresponding partial differential equations (PDE) only involve first-order derivatives,
their solutions are not easy because they usually contain both discontinuities and sophisticated structures with multi-scales.
Therefore, high-order numerical schemes with excellent shock-capturing and multi-resolution properties are desired 
for solving hyperbolic conservation laws.
The essentially non-oscillatory (ENO) schemes and weighted ENO (WENO) schemes are cutting-edge high-order shock-capturing schemes
and achieve great success in practice.

Harten and his coworkers \cite{Harten1987ENO} proposed the first ENO scheme which 
adaptively selects the smoothest stencil from several local candidates of the same size
such that the constructed scheme achieves uniformly high order in smooth regions 
and suppresses non-physical oscillations near discontinuities.
Shu and Osher \cite{Shu1988EfficientENO, Shu1989EfficientENOII} significantly improve 
the efficiency of ENO schemes by the method of lines.
In their approach, the cell average reconstructions are replaced by numerical flux reconstructions
which can be implemented in a dimension-by-dimension manner, and the time is discretized 
by an efficient Runge-Kutta scheme which preserves strong stability.
The drawback of ENO schemes is that they only achieve the order of a small sub-stencil 
even when the whole global stencil is smooth.
To improve ENO schemes, Liu \emph{et al}. \cite{Liu1994WENO} proposed a WENO scheme
which combines all the candidate stencils via nonlinear weights.
The weight becomes very tiny if the corresponding stencil contains discontinuities
so that the constructed schemes achieve non-oscillatory features near discontinuities.
However, these WENO schemes still cannot achieve the optimal order when the global stencils are smooth.
Jiang and Shu \cite{Jiang_Shu1996WENO} systematically analyzed WENO schemes
and provided a general framework for calculating smoothness indicators and nonlinear weights 
so that the constructed WENO schemes can retrieve the optimal order in smooth regions 
while retaining oscillatory at discontinuities.
Henrick \emph{et al}. \cite{Henrick2005WENO_M} discovered that the WENO schemes proposed by Jiang and Shu 
lose accuracy at critical points (local extrema) if the parameter $\varepsilon$,
which is used to avoid dividing by 0, is set to a very small number.
They proposed a simple mapping function to modify the nonlinear to cure this issue.
Borges \emph{et al}. \cite{Borges2008WENO_Z} further improved the accuracy of WENO schemes 
by introducing a global smoothness and proposed the so-called WENO-Z type nonlinear weights.
Yamaleev and Carpenter \cite{YamaleevESWENO1, YamaleevESWENO2} proposed energy stable WENO (ESWENO)  
for systems of linear hyperbolic equations with both continuous and discontinuous solutions. 
The stability of $L_2$ norm is explicitly achieved (even when a purely downwind stencil is included)
by enforcing a nonlinear summation-by-parts (SBP) condition at each instant in time.
Hu \emph{et al}. \cite{Hu2010WENO_CU} constructed an adaptive central-upwind scheme
by adding an extra purely downwind stencil which is assigned with
the smoothness indicator of the global central stencil. 
When the central stencil is nonsmooth, the effect of the purely downwind stencil 
is automatically excluded and the upwind shock-capturing WENO scheme is recovered.
Fu \emph{et al.} \cite{Fu2016TENO,Fu2017TENO2,Fu2018TENO_AO,Fu2021TENO_AA,Fu2023TENOReview} 
proposed a family of targeted ENO (TENO) schemes
which find out candidate stencils containing discontinuities 
and reconstruct an optimal flux based on the remaining stencil by using a scale separation strategy.
In this way, TENO schemes significantly reduce numerical dissipation 
while maintaining ENO properties at discontinuities.
Balsara and Shu \cite{Balsara2000Monotonicity_WENO} found that very high-order extensions of WENO schemes 
may not always be monotonicity preserving.
They constructed increasingly
high-order (up to 11th-order) monotonicity preserving WENO by coupling the nonlinear weights 
of Jiang and Shu \cite{Jiang_Shu1996WENO}
and the monotonicity preserving bounds of Suresh and Huynh \cite{Suresh1997Monotonicity}.
Gerolymos \emph{et al}. \cite{Gerolymos2009VHO_WENO} extended the mapped WENO scheme of 
Henrick \emph{et al}. \cite{Henrick2005WENO_M} to very-high-order (up to 17th-order)
by coupling with a recursive-order-reduction (ROR) algorithm.
Wu \emph{et al}. \cite{Wu2021VHO_WENO} improved the efficiency of the very-high-order WENO
schemes by using efficient smoothness indicators \cite{Wu2020SI4Sine} 
and a predetermined order reduction method.
Nowadays, there are many variants of ENO/WENO schemes 
and other shock-capturing high-order schemes adopting ENO/WENO idea,
such as the weighted compact nonlinear schemes (WCNS) \cite{Deng2000WCNS},
the central WENO (CWENO) schemes \cite{Bianco1999CENO, Levy1999CWENO, Levy2000CWENO_ANM},
the hybrid compact-WENO schemes \cite{Pirozzoli2002HybridWENO,Ren2003hybridCompatWENO},
the Hermite WENO (HWENO) schemes \cite{Qiu2004HermiteWENO, Qiu2005HermiteWENO2D}, 
the  WENO-ADER schemes \cite{Toro2002ADER, Titarev2004WENO_FV, Balsara2013ADER}, 
the $P_NP_M$ schemes \cite{Dumbser2008PNPM, Dumbser2009PNPM, Dumbser2010PNPM}, and so forth.
It is impossible to review all of them here.
For more details on the development of ENO/WENO schemes,
one may refer to the review by Shu \cite{Shu2009SIAMReviewWENO, Shu2016ReviewWENO_DG}.

The classical ENO/WENO schemes are constructed based on candidate stencils of an equal size
which cannot achieve optimal accuracy in some cases.
Taking five-point schemes for example,
ENO schemes only choose the smoothest stencil among the three three-point candidate sub-stencils,
and WENO schemes combine all three sub-stencils by nonlinear weights 
which may recover the optimal accuracy on the global stencil.
However, in some cases, a non-smooth five-point global stencil may contain a smooth four-point sub-stencil,
and in some other cases, all the three three-point sub-stencils may be non-smooth.
Therefore, it is a natural choice to construct numerical schemes based on candidate stencils of unequal sizes.
Levy \emph{et al}. \cite{Levy2000CompactCWENO} first proposed a compact central WENO reconstruction
which is a nonlinear combination of two two-point stencils and one three-point stencil.
In smooth regions, the scheme approaches the three-point stencil, 
while near discontinuities, the scheme approaches the smoother two-point stencil for shock-capturing.
The advantage of this approach is that it can recover the optimal high-order accuracy
without constraint on the choice of linear weights like classical WENO schemes.
In recent years, CWENO schemes based on the reconstruction technique by using stencils of unequal sizes
have been substantially analyzed and developed by Puppo, Semplice, Dumbser, and Cravero \emph{et al}. 
\cite{Puppo2021Quinpi_CWENO, Semplice2016CWENO_Adaptive, Semplice2021CWENO_Boundary, Semplice2020EfficientWENO_AO,Dumbser2017CWENO,Cravero2016JSC,Cravero2018CWENO, Cravero2018Cool_WENO}.
Cravero, Semplice, and Visconti \cite{Cravero2019CWENO} comprehensively reviewed such CWENO schemes.
Inspired by the weighting strategy of Levy \emph{et al}. \cite{Levy2000CompactCWENO},
Zhu \emph{el al.} \cite{Zhu2016NewWENOFD, Zhu2018NewMultiResolutionWENO, Zhu2019NewMultiResolutionWENO_Tri, Zhu2020NewMultiResolutionWENO_Tet} 
and Balsara \emph{et al.} \cite{Balsara2016WENO_AO, Balsara2020WENO_AO_Unstructure}
constructed WENO schemes with adaptive order (WENO-AO) and multi-resolution (WENO-MR).
Shen \cite{Shen2021WENO_AOA} proposed a simplified weighting strategy for combining polynomials of different degree
and adopted it to construct weighted compact central schemes (WCCS) \cite{Shen2021WCCS}.
Fu \emph{et al.} \cite{Fu2018TENO_AO} constructed high-order schemes with adaptive order 
based on candidate stencils of unequal sizes under the framework of TENO schemes.
Shen \cite{Shen2023ENO_AO} constructed a class of ENO schemes with adaptive order (ENO-AO)
by using newly defined smoothness indicators which measure the minimum error
between the reconstructed polynomial and the target function at the cells adjacent to corresponding stencils.
The computation of the new smoothness indicator is simpler than
that of Jiang and Shu's smoothness indicator \cite{Jiang_Shu1996WENO},
so ENO-AO schemes cost less CPU times than WENO-AO and TENO-AO schemes 
adopting Jiang and Shu's smoothness indicator.

Multi-level ENO/WENO schemes constructed on candidate stencils of unequal sizes
naturally have some advantages over those schemes constructed on candidate stencils of an equal size,
so they have attracted a lot of attention and have undergone major development recently.
However, there are still some challenges to be resolved.
Firstly, the high-order accuracy is easily 
polluted by low-order candidate stencils, if special care is not exerted.
Arbogast \emph{et al.} \cite{Arbogast2018WENO_AO} proved that 
the highest order $r$ and the lowest order $s$ of multi-level reconstructions
must satisfy $r \le 2s-1$ when WENO-JS weights are used;
and there is no constraint on $r$ and $s$ when WENO-Z weights are used,
providing that enough intermediate levels or large enough power parameters in the weights are used.
Secondly, when extended to very high-order situations, multi-level reconstruction may contain
a lot of candidate stencils of different levels which significantly increases the complexity.
We notice that existing multi-level ENO/WENO schemes rarely include low-order 
one-point and two-point candidate stencils.
The existing high-order multi-level schemes,
e.g., 9th-order WENO-AO(9,5,3) of Balsara \emph{et al.} \cite{Balsara2016WENO_AO},
10th-order TENO10-AA scheme of Fu \cite{Fu2021TENO_AA}, 
and 9th-order multi-resolution WENO scheme of Zhu and Shu \cite{Zhu2018NewMultiResolutionWENO},
do not include candidate stencils of all levels.
Zhu \emph{et al.} \cite{Zhu2016NewWENOFD,Zhu2018NewMultiResolutionWENO,Zhu2019NewMultiResolutionWENO_Tri,Zhu2020NewMultiResolutionWENO_Tet} 
constructed WENO schemes by combining first- and second-order polynomials with high-order polynomials.
To achieve optimal accuracy in smooth regions,
the exponent of global smoothness indicators in the numerator must be
larger than that appears in the denominator.
As a result, the weights are not dimensionless, 
and the performance of the constructed schemes is problem-dependent \cite{Balsara2016WENO_AO}. 
The ENO-AO schemes of Shen \cite{Shen2023ENO_AO} consider
one-point and two-point stencils as candidates,
but they have a sensitive parameter $\varepsilon$ as well.
When the smoothness indicator of a high-order stencil is smaller than $\varepsilon$,
it will be considered to be a smooth stencil which guarantees optimal convergence order in smooth regions.
Apparently, the choice of $\varepsilon$ may affect the convergence order.
Some studies particularly discussed the effect of $\varepsilon$ \cite{Don2013WENO_D, Wang2019WENO_A, Don2022WENO_Si}.
In addition, ENO-AO schemes select the smoothest stencil from all candidates.
When extending to very high-order cases, the complexity of the selecting procedure will significantly increase.

In this study, we further simplify and improve the selecting strategy of ENO-AO schemes \cite{Shen2023ENO_AO}
and construct an efficient class of very-high-order ENO schemes with multi-resolution (ENO-MR)
which have the following advantages:
\begin{itemize}
  \item [1.] The candidates of ENO-MR schemes consist of all levels ranging from
  first-order up to the designed high-order (up to 17th-order).
  \item [2.] The selecting strategy only requires the computation of the highest-order derivatives
   of a part of candidate stencils which is beneficial to efficient implementation.
  \item [3.] ENO-MR schemes are endowed with an adaptive order-reduction mechanism, 
  so very-high-order extensions do not require additional techniques to improve the robustness
  as very-high-order WENO schemes of Gerolymos \emph{et al}. \cite{Gerolymos2009VHO_WENO} 
  and Wu \emph{et al}. \cite{Wu2021VHO_WENO} do. 
  \item [4.] ENO-MR schemes are completely parameter-free and essentially scale-invariant.
  \item [5.] ENO-MR schemes achieve designed high-order convergence at high-order critical points (local extrema)
             without any additional treatment.
\end{itemize}

%%%%%%%%%%%%%%%%%%%%%%%%%%%%%%%%%%%%%%%%%%%%%%%%%%%%%%%%%%%%%%%%%%%%
\section{Finite difference ENO/WENO schemes}\label{SEC:ENO_WENO_FD}
\subsection{A semi-discretized conservative finite difference scheme}
We consider the one-dimensional scalar conservation law,
\begin{equation}\label{Eq:1DHCL_Eq}
    \frac{\partial u}{\partial t}+\frac{\partial f(u)}{\partial x}=0, t\in[0,\infty),
\end{equation}
in the spatial domain $[x_L,x_R]$ that is discretized into uniform intervals
by $x_j=x_L+(j-1)h$ ($j=1 \text{ to } N+1$), where $h=(x_R-x_L)/N$.
Define a function $\mathcal{H}(x)$ implicitly as \cite{Shu1988EfficientENO}
\begin{equation}\label{Eq:h_x}
  f(u(x))=\frac{1}{h}\int_{x-h/2}^{x+h/2}\mathcal{H}(\xi)d\xi.
\end{equation}
Then we can get a semi-discretized conservative finite difference scheme as 
\begin{equation}\label{Eq:1DWENOFD}
  \frac{du_j(t)}{dt}=\mathcal{L}(u_j(t))=-\frac{\hat{f}_{j+1/2}-\hat{f}_{j-1/2}}{h},
\end{equation}
where the numerical flux $\hat{f}_{j\pm1/2}$ is an approximation of 
the function $\mathcal{H}(x)$ at $x_{j\pm1/2}$.
The attractive feature of this approach is that
it can be easily extended to multi-dimensional cases 
in a dimension-by-dimension manner.

To take account of the upwind mechanism
which can improve the robustness of the scheme,
we usually split the flux into two parts as
\begin{equation}\label{Eq:Flux_Split}
  f(u)=f^+(u)+f^-(u),
\end{equation}
where $\frac{df^+(u)}{du}\ge0$ and $\frac{df^-(u)}{du}\le0$.
Here, we adopt the global Lax–Friedrichs flux splitting method which is expressed as
\begin{equation}\label{Eq:LF_Flux_Split}
  f^\pm(u)=\frac{1}{2}(f(u)\pm\alpha u),
\end{equation}
where $\alpha=\mbox{max}\left|\frac{df(u)}{du}\right|$
and the maximum is taken over all mesh points on one axis-aligned line.
%%%%%%%%%%%%%%%%%%%%%%%%%%%%%%%%%%%%%%%%%%%%%%%%%%%%%%%%%%%%%%%%%%%%%%%%%%%%%%%%%
\subsection{Review of ENO/WENO flux reconstructions}
\subsubsection{Classical ENO/WENO reconstructions on candidate stencils of an equal size}
Now, the remaining task is finding a proper way to calculate $\hat{f}_{j\pm1/2}^\pm(u)$.
To this end, we define a stencil $S_{j-m}^{j+n}$ as a set of successive intervals including $I_j$, 
i.e.,$S_{j-m}^{j+n}:=\{I_{j-m},...,I_j,...,I_{j+n}\}$ ($m\ge 0,n\ge 0$).
Then, we can trivially reconstruct a polynomial $P_{j-m}^{j+n}(x)$ of degree $m+n$ to approximate $\mathcal{H}(x)$ 
by enforcing the conservation conditions on $S_{j-m}^{j+n}$, i.e.,
\begin{equation}\label{Eq:Polynomial_Construction}
  \frac{1}{h}\int_{I_k} P_{j-m}^{j+n}(\xi)d\xi=f_k, \quad k=j-m \text{ to }j+n.
\end{equation}
Finally, the fluxes are simply calculate as $\hat{f}_{j\pm1/2}=P_{j-m}^{j+n}\left(x_{j\pm1/2}\right)$
which can be written in terms of $f_k$ as
\begin{equation}\label{Eq:ReconstructedFlux}
  \hat{f}_{j+1/2}=\sum_{l=-m}^{n} a_lf_{j+l}.
\end{equation}
We note that we do not distinguish between $f^+$ and $f^-$ hereafter,
because they are calculated by the same principle.

We have many choices for the stencil.
When $f(u(x))$ is sufficiently smooth on $S_{j-m}^{j+n}$, 
$P_{j-m}^{j+n}$ is a $(m+n+1)$th-order approximation to $\mathcal{H}(x)$.
Otherwise, $P_{j-m}^{j+n}$ is not a good approximation to $\mathcal{H}(x)$.
When reconstructing fluxes, we should avoid stencils containing discontinuities.
The $(2r-1)$-point ENO scheme \cite{Harten1987ENO} selects the smoothest stencil 
from the pool of $r$-point candidate stencils 
$\{S_{j-r+1}^{j},S_{j-r+2}^{j+1},...,S_{j-1}^{j+r-2},S_{j}^{j+r-1}\}$, 
so it can only achieve at most $r$th-order accuracy.
Classical WENO schemes use a weighted combination of the reconstructions on the same candidate stencils
which can retrieve the $(2r-1)$-order reconstruction on the global stencil $S_{j-r+1}^{j+r-1}$ in smooth regions
while retaining ENO properties at discontinuities.
The formulae of standard $(2r-1)$-point WENO reconstructions are expressed as
 \begin{subequations}\label{Eq:standard_WENO}
  \begin{equation}
    P_{WENO}(x)=\sum_{k=0}^{r-1}\omega_k^{(r)}P_{j-r+k+1}^{j+k}(x),
  \end{equation}
  \text{with}
  \begin{equation}
    \omega_k^{(r)}=\frac{\alpha_k^{(r)}}{\sum_{l=0}^{r-1}\alpha_l^{(r)}}, k=0, ..., r-1.
  \end{equation}
 \end{subequations}
 The nonlinear weights of Jiang and Shu \cite{Jiang_Shu1996WENO} are calculated as
 \begin{subequations}\label{Eq:WENO_JS_Weight}
  \begin{equation}
    \alpha_k^{(r,JS)}=\frac{d_k^{(r)}}{\left(\beta_k^{(r)}+\varepsilon\right)^p}, k=0, ..., r-1,
  \end{equation}
  \text{with the smoothness indicators defined as}
  \begin{equation}\label{Eq:Smoothness_Indicator}
    \beta_k^{(r)}=\sum_{l=1}^{r-1}h^{2l-1}\int_{x_{i-1/2}}^{x_{i+1/2}}\left(\frac{d^lP_{j-r+k+1}^{j+k}(x)}{dx^l}\right)^2dx, k=0, ..., r-1.
  \end{equation}
 \end{subequations}
 The WENO-Z type nonlinear weights \cite{Borges2008WENO_Z} are calculated as 
 \begin{subequations}\label{Eq:WENO_Z_Weight}
  \begin{equation}
    \alpha_k^{(r,Z)}=d_k^{(r)}\left(1+\left(\frac{\tau_{2r-1}}{\beta_k^{(r)}+\varepsilon}\right)^p\right), k=0, ..., r-1,
  \end{equation}
  \text{with the global smoothness indicators defined as \cite{Castro2011High_Order_WENO_Z}}
  \begin{equation}
    \tau_{2r-1}=\begin{cases}
      \left|\beta_0^{(r)}-\beta_{r-1}^{(r)}\right|, \text{mod}(r,2)=1,\\
      \\
      \left|\beta_0^{(r)}-\beta_1^{(r)}-\beta_{r-2}^{(r)}+\beta_{r-1}^{(r)}\right|, \text{mod}(r,2)=0.
    \end{cases}
  \end{equation}
 \end{subequations}
 We note that in both WENO-JS and WENO-Z weights, there are common parameters,
 namely $d_k^{(r)}$, $\varepsilon$, and $p$.
  $d_k^{(r)}$ are optimal linear weights which guarantee the optimal $(2r-1)$th-order in smooth regions.
  They may become negative or even non-existent in some cases.
 $\varepsilon$ is a tiny positive number used to avoid dividing by 0.
 The convergence rates of WENO schemes are found to be sensitive to $\varepsilon$ \cite{Don2013WENO_D}.
Gerolymos \emph{et al}. \cite{Gerolymos2009VHO_WENO} found that
$p=2$ is sufficient to keep ENO property for $r\le6$ cases, 
but we have to use $p=r$ for $r>6$ cases.
When solving Euler equations by very-high-order WENO schemes, 
the interaction between characteristic fields and/or the absence of a zone of smoothness of $r$
points to choose a stencil from may cause serious non-physical oscillations \cite{Harten1987ENO}.
Order-reduction techniques are employed to retain ENO property \cite{Gerolymos2009VHO_WENO,Wu2021VHO_WENO}.

%%%%%%%%%%%%%%%%%%%%%%%%%%%%%%%%%%%%%%%%%%%%%%%%%%%%%%%%%%%%%%%%%%%%%%%%%%%%%%%%%
\subsubsection{ENO/WENO reconstructions on candidate stencils of unequal sizes}
Levy \emph{et al}. \cite{Levy2000CompactCWENO} first proposed a CWENO reconstruction
based on two two-point stencils and one three-point stencil.
In recent years, many WENO reconstructions with adaptive order (WENO-AO) 
with similar forms have been proposed 
based on this framework.
Here, we briefly summarize the WENO-AO schemes proposed by Balsara \emph{et al.} \cite{Balsara2016WENO_AO}.
The candidate stencils of WENO-AO($2r-1$,3) consist of 
$\{S_{j-r+1}^{j+r-1},S_{j-2}^{j},S_{j-1}^{j+1},S_{j}^{j+2}\}$ ($r\ge3$).
The final reconstructed polynomial is expressed as
\begin{equation}\label{Eq:WENO_AO}
  \begin{aligned}
    P_{AO(2r-1,3)}(x)&=\frac{\omega_{r-1}^{(2r-1)}}{d_{r-1}^{(2r-1)}}\left(P_{j-r+1}^{j+r-1}(x)-d_0^{(3)}P_{j-2}^{j}(x)-d_1^{(3)}P_{j-1}^{j+1}(x)-d_2^{(3)}P_{j}^{j+2}(x)\right)\\
                     &+\omega_0^{(3)}P_{j-2}^{j}(x)+\omega_1^{(3)}P_{j-1}^{j+1}(x)+\omega_2^{(3)}P_{j}^{j+2}(x),
  \end{aligned}
\end{equation}
where $\omega$ can be calculated 
via WENO-JS weights or WENO-Z weights.
The corresponding linear weights can be set as
\begin{equation}\label{WENO_AO_linear_weights}
  d_{r-1}^{(2r-1)}=\gamma_{Hi},\quad d_0^{(3)}=d_2^{(3)}=(1-\gamma_{Hi})(1-\gamma_{Lo})/2,\quad d_1^{(3)}=(1-\gamma_{Hi})\gamma_{Lo},
\end{equation}
where $\gamma_{Hi}\in [0.85,0.95]$ and $\gamma_{Lo}\in [0.85,0.95]$.
Since the global $(2r-1)$-point stencil itself is a candidate stencil,
WENO-AO reconstructions can recover the optimal $(2r-1)$th-order 
without particular constraints on the choice of linear weights.
Furthermore, unlike WENO reconstructions based on candidate stencils of an equal size,
the smallest candidate stencils of WENO-AO always contain three points 
regardless of the size of the global stencil.
We can also include two-point stencils, for example 
\cite{Levy2000CompactCWENO,Zhu2016NewWENOFD,Zhu2018NewMultiResolutionWENO,Zhu2019NewMultiResolutionWENO_Tri,Zhu2020NewMultiResolutionWENO_Tet}.

Combining several WENO-AO($2r-1$,3) reconstructions together,
Balsara \emph{et al.} \cite{Balsara2016WENO_AO} further proposed 
recursively defined WENO-AO reconstructions.
For example, the WENO-AO(7,5,3) reconstruction can be expressed as
\begin{subequations}
  \begin{equation}
    P_{AO(7,5,3)}(x)=\frac{\omega_{AO(7,3)}}{\gamma_{Hi}}\left(P_{AO(7,3)}(x)-(1-\gamma_{Hi})P_{AO(5,3)}(x)\right)+\omega_{AO(5,3)}P_{AO(5,3)},
  \end{equation}
  \text{where,}
  \begin{equation}
    \omega_{AO(7,3)}=\frac{\alpha_{AO(7,3)}}{\alpha_{AO(7,3)}+\alpha_{AO(5,3)}}, \quad
    \omega_{AO(5,3)}=\frac{\alpha_{AO(5,3)}}{\alpha_{AO(7,3)}+\alpha_{AO(5,3)}}
  \end{equation}
  \begin{equation}\label{Eq:WENO_AO_Weights}
    \alpha_{AO(7,3)}=\gamma_{Hi}\left(1+\frac{\sigma}{\beta_3^{(7)}+\varepsilon}\right),
    \alpha_{AO(5,3)}=(1-\gamma_{Hi})\left(1+\frac{\sigma}{\beta_2^{(5)}+\varepsilon}\right)
  \end{equation}
  \text{with}
  \begin{equation}
    \sigma=\left|\beta_3^{(7)}-\beta_2^{(5)}\right|.
  \end{equation}
\end{subequations}

In theory, we can construct WENO-AO schemes by combining candidate stencils of arbitrary points.
However, the nonlinear weights depend on the smoothness indicators defined by Eq. (\ref{Eq:Smoothness_Indicator})
which is computationally expensive for high-order polynomials.
Therefore, it is not practical to recursively construct very high-order WENO-AO schemes with many levels.

In order to improve the efficiency, Shen \cite{Shen2023ENO_AO} proposed new smoothness indicators of $S_{j-m}^{j+n}$
\begin{subequations}\label{Eq:New_Smoothness_Indicator}
  \begin{equation}\label{SubEq:IS}
      IS_{j-m}^{j+n}={\rm MIN}(\delta_L,\delta_R),
  \end{equation}
  \text{where,}
  \begin{equation}\label{SubEq:deltaL}
    \delta_L=\begin{cases}
      \frac{1}{2}\left(\left|f_j-f_{j-1}\right|+\left|f_{j-1}-f_{j-2}\right|\right), \quad m=n=0,\\
      \left|\frac{1}{h}\int_{I_{j-m-1}} P_{j-m}^{j+n}(\xi)d\xi-f_{j-m-1}\right|, \quad {\rm otherwise},
    \end{cases}
  \end{equation}
  \begin{equation}\label{SubEq:deltaR}
    \delta_R=\begin{cases}
      \frac{1}{2}\left(\left|f_j-f_{j+1}\right|+\left|f_{j+1}-f_{j+2}\right|\right), \quad m=n=0,\\
      \left|\frac{1}{h}\int_{I_{j+n+1}} P_{j-m}^{j+n}(\xi)d\xi-f_{j+n+1}\right|, \quad {\rm otherwise}.
    \end{cases}
  \end{equation}
  \end{subequations}
  The new smoothness indicators measure the minimum error
between the reconstructed polynomial and the target function
at the cells adjacent to $S_{j-m}^{j+n}$ which is not sensitive to the stencil size.
Obviously, the new smoothness indicators are computationally cheaper than
Jiang and Shu's smoothness indicators which are defined by Eq. (\ref{Eq:Smoothness_Indicator}).
Utilizing the new smoothness indicators, a $(2r-1)$-point ENO-AO schemes are constructed in the following way:
\begin{enumerate}
  \item [1.] We choose the linearly stable sub-stencils of $S_{j-r+1}^{j+r-1}$ as candidates;
  \item [2.] If the smoothness indicator of a high-order candidate stencil is smaller than $\varepsilon$, 
  we use it to reconstruct the flux;
  \item[3.] Otherwise, we choose the stencil with the smallest smoothness indicator among all candidates 
  to reconstruct the flux.
\end{enumerate}

%%%%%%%%%%%%%%%%%%%%%%%%%%%%%%%%%%%%%%%%%%%%%%%%%%%%%%%%%%%%%%%%%%%%%%%%%%%%%%%%
\subsection{Time discretization}
Once we have the semi-discretized scheme Eq. (\ref{Eq:1DWENOFD}),
we use a stable ordinary differential equation (ODE) solver to get a fully discretized scheme.
Theoretically, we prefer an ODE solver that owns the same order as the spatial operator.
However, very high-order ODE solvers for general nonlinear cases are either not readily developed,
or readily available, or else require a very large number of stages \cite{Butcher2006RK}.
For general nonlinear problems, the popular method is the third-order 
strong-stability-preserving Runge-Kutta (SSP-RK3) scheme \cite{Shu1988EfficientENO},
i.e., 
\begin{subequations}\label{Eq:TVD3rdRK}
  \begin{equation}\label{SubEq:RKStage1}
    u_j^{(1)}=u_j^n+\Delta t\mathcal{L}(u^n),
  \end{equation}
  \begin{equation}\label{SubEq:RKStage2}
    u_j^{(2)}=\frac{3}{4}u_j^n+\frac{1}{4}\left(u_j^{(1)}+\Delta t\mathcal{L}(u^{(1)})\right),
  \end{equation}
  \begin{equation}\label{SubEq:RKStage3}
    u_j^{n+1}=\frac{1}{3}u_j^n+\frac{2}{3}\left(u_j^{(2)}+\Delta t\mathcal{L}(u^{(2)})\right),
  \end{equation}
  \end{subequations}
where $\mathcal{L}$ represents the spatial operator.
For linear problems, Gottlieb \emph{et al.} \cite{Gottlieb2001SSPRK, Gottlieb2003lSSPRK, Gottlieb2005lSSPRK}
developed a family of linear strong-stability-preserving Runge-Kutta schemes ($\ell$SSP-RK($m,m-1$))
which can achieve $(m-1)$th-order with $m$ stages.
The formulae of $\ell$SSP RK($m,m-1$) are given as
\begin{subequations}\label{Eq:lSSPRK}
  \begin{equation}
    u_j^{(0)}=u_j^n,
  \end{equation}
  \begin{equation}
    u_j^{(s)}=u_j^{(s-1)}+\frac{1}{2}\Delta t\mathcal{L}(u^{(s-1)}), s=1, ..., m-1,
  \end{equation}
  \begin{equation}
    u_j^{n+1}=\sum_{k=0}^{m-2}\alpha_{m,k}u_j^{(k)}+\alpha_{m,m-1}\left(u_j^{(m-1)}+\frac{1}{2}\Delta t\mathcal{L}(u^{(m-1)})\right),
  \end{equation}
\end{subequations}
where the coefficients are recursively computed as
\begin{subequations}\label{Eq:lSSPRK_Coe}
  \begin{equation}
    \alpha_{2,0}=0, \alpha_{2,1}=1,
  \end{equation}
  \begin{equation}
    \alpha_{m,k}=\frac{2}{k}\alpha_{m-1,k-1}, k=1, ..., m-2,
  \end{equation}
  \begin{equation}
    \alpha_{m,m-1}, \frac{2}{m}\alpha_{m-1,m-2}, \alpha_{m,0}=1-\sum_{k=1}^{m-1} \alpha_{m,k}.
  \end{equation}
\end{subequations}

In this study, we use SSP-RK3 for nonlinear problems
and $\ell$SSP-RK($m,m-1$) coupled with $(m-1)$th-order spatial operators for linear problems.

%%%%%%%%%%%%%%%%%%%%%%%%%%%%%%%%%%%%%%%%%%%%%%%%%%%%%%%%%%%%%%%%%%%%%
\section{Increasingly high-order ENO schemes with multi-resolution}
We construct increasingly high-order ENO schemes with multi-resolution
in the framework of conservative finite difference scheme as described in Section \ref{SEC:ENO_WENO_FD}.
This study mainly focuses on the flux reconstructions,
the other procedures remain unchanged.
\subsection{Candidate stencils of multi-levels}

For a $(2r-1)$-point scheme, there are $r^2$ stencils in total that can be used to reconstruct $f_{j+1/2}$.
However, as ENO-AO schemes of Shen \cite{Shen2023ENO_AO},
we only choose linearly stable stencils as candidates to guarantee stability.
Particularly, the candidates of a 17-point scheme include the following 29 stencils:
\begin{equation}\label{Eq:Candidates}
  \begin{aligned}
    &S_{j-8}^{j+8},S_{j-7}^{j+8},S_{j-8}^{j+7},S_{j-7}^{j+7},S_{j-8}^{j+6},S_{j-6}^{j+7},S_{j-7}^{j+6},S_{j-6}^{j+6},S_{j-7}^{j+5},\\
    &S_{j-5}^{j+6},S_{j-6}^{j+5},S_{j-5}^{j+5},S_{j-6}^{j+4},S_{j-4}^{j+5},S_{j-5}^{j+4},S_{j-4}^{j+4},S_{j-5}^{j+3},S_{j-3}^{j+4},S_{j-4}^{j+3},\\
    &S_{j-3}^{j+3},S_{j-2}^{j+3},S_{j-3}^{j+2},S_{j-2}^{j+2},S_{j-1}^{j+2},S_{j-2}^{j+1},S_{j-1}^{j+1},S_{j}^{j+1},S_{j-1}^{j},S_{j}^{j}.
  \end{aligned}
\end{equation}
In the pool of candidate stencils, 
there is one candidate of 1-point, 3-point, 5-point, 7-point, and 17-point stencils,
and there are two candidates of other levels.
Tables \ref{Table:Flux_reconstruction1} and \ref{Table:Flux_reconstruction2} 
give the coefficients of $a_l$ in Eq. (\ref{Eq:ReconstructedFlux})
for computing $\hat{f}^+_{j+1/2}$ on all candidate stencils.
We note that the formulae of $\hat{f}_{j+1/2}^-$ and $\hat{f}_{j+1/2}^+$ 
are symmetric with respect to $x_{j+1/2}$,
so we only provide the formulae of $\hat{f}_{j+1/2}^+$.

\begin{table}[htbp]
  
  \caption{The coefficients $a_l$ in Eq. (\ref{Eq:ReconstructedFlux}) for computing $\hat{f}_{j+1/2}^+$ on different candidate stencils.}
  \label{Table:Flux_reconstruction1}
  \centering
  \begin{tabular}{|c|c|c|c|c|c|c|c|c|c|}
    % after \\: \hline or \cline{col1-col2} \cline{col3-col4} ...
     \hline
      $l$&$S_{j-8}^{j+8}$&$S_{j-7}^{j+8}$&$S_{j-8}^{j+7}$&$S_{j-7}^{j+7}$&$S_{j-8}^{j+6}$&$S_{j-6}^{j+7}$&$S_{j-7}^{j+6}$&$S_{j-6}^{j+6}$&$S_{j-7}^{j+5}$  \\
      \hline
      -8&$\frac{56}{12252240}$&                      &$\frac{7}{720720}$&                     &$\frac{8}{360360}$&&&&\\
      \hline
      -7&$\frac{-1015}{12252240}$&$\frac{-7}{720720}$&$\frac{-119}{720720}$&$\frac{-7}{360360}$&$\frac{-127}{360360}$&                 &$\frac{-15}{360360}$&                 &$\frac{-35}{360360}$\\
     \hline
     -6&$\frac{8777}{12252240}$&$\frac{121}{720720}$&$\frac{961}{720720}$&$\frac{113}{360360}$&$\frac{953}{360360}$&$\frac{15}{360360}$&$\frac{225}{360360}$&$\frac{30}{360360}$&$\frac{485}{360360}$\\
     \hline
     -5&$\frac{-48343}{12252240}$&$\frac{-999}{720720}$&$\frac{-4919}{720720}$&$\frac{-867}{360360}$&$\frac{-4507}{360360}$&$\frac{-230}{360360}$&$\frac{-1595}{360360}$&$\frac{-425}{360360}$&$\frac{-3155}{360360}$\\
     \hline
     -4&$\frac{191561}{12252240}$&$\frac{5273}{720720}$&$\frac{18013}{720720}$&$\frac{4229}{360360}$&$\frac{15149}{360360}$&$\frac{1681}{360360}$&$\frac{7141}{360360}$&$\frac{2851}{360360}$&$\frac{12861}{360360}$\\
     \hline
     -3&$\frac{-588127}{12252240}$&$\frac{-20207}{720720}$&$\frac{-50783}{720720}$&$\frac{-14881}{360360}$&$\frac{-38905}{360360}$&$\frac{-7874}{360360}$&$\frac{-22889}{360360}$&$\frac{-12164}{360360}$&$\frac{-37189}{360360}$\\
     \hline
     -2&$\frac{1491041}{12252240}$&$\frac{61329}{720720}$&$\frac{117385}{720720}$&$\frac{41175}{360360}$&$\frac{81215}{360360}$&$\frac{27161}{360360}$&$\frac{57191}{360360}$&$\frac{37886}{360360}$&$\frac{82931}{360360}$\\
     \hline
     -1&$\frac{-3409855}{12252240}$&$\frac{-162895}{720720}$&$\frac{-242975}{720720}$&$\frac{-98965}{360360}$&$\frac{-150445}{360360}$&$\frac{-77944}{360360}$&$\frac{-122989}{360360}$&$\frac{-97249}{360360}$&$\frac{-157309}{360360}$\\
     \hline
     0&$\frac{8842385}{12252240}$&$\frac{477745}{720720}$&$\frac{567835}{720720}$&$\frac{261395}{360360}$&$\frac{312875}{360360}$&$\frac{237371}{360360}$&$\frac{288851}{360360}$&$\frac{263111}{360360}$&$\frac{323171}{360360}$\\
     \hline
      1&$\frac{7481025}{12252240}$&$\frac{477745}{720720}$&$\frac{397665}{720720}$&$\frac{216350}{360360}$&$\frac{176310}{360360}$&$\frac{237371}{360360}$&$\frac{192326}{360360}$&$\frac{211631}{360360}$&$\frac{166586}{360360}$\\
      \hline
      2&$\frac{-2320767}{12252240}$&$\frac{-162895}{720720}$&$\frac{-106839}{720720}$&$\frac{-63930}{360360}$&$\frac{-39906}{360360}$&$\frac{-77944}{360360}$&$\frac{-47914}{360360}$&$\frac{-58639}{360360}$&$\frac{-33614}{360360}$\\
     \hline
     3&$\frac{797985}{12252240}$&$\frac{61329}{720720}$&$\frac{30753}{720720}$&$\frac{20154}{360360}$&$\frac{9234}{360360}$&$\frac{27161}{360360}$&$\frac{12146}{360360}$&$\frac{16436}{360360}$&$\frac{6426}{360360}$\\
     \hline
     4&$\frac{-241599}{12252240}$&$\frac{-20207}{720720}$&$\frac{-7467}{720720}$&$\frac{-5326}{360360}$&$\frac{-1686}{360360}$&$\frac{-7874}{360360}$&$\frac{-2414}{360360}$&$\frac{-3584}{360360}$&$\frac{-854}{360360}$\\
     \hline
     5&$\frac{58281}{12252240}$&$\frac{5273}{720720}$&$\frac{1353}{720720}$&$\frac{1044}{360360}$&$\frac{204}{360360}$&$\frac{1681}{360360}$&$\frac{316}{360360}$&$\frac{511}{360360}$&$\frac{56}{360360}$\\
     \hline
     6&$\frac{-10263}{12252240}$&$\frac{-999}{720720}$&$\frac{-159}{720720}$&$\frac{-132}{360360}$&$\frac{-12}{360360}$&$\frac{-230}{360360}$&$\frac{-20}{360360}$&$\frac{-35}{360360}$&\\
     \hline
     7&$\frac{1161}{12252240}$&$\frac{121}{720720}$&$\frac{9}{720720}$&$\frac{8}{360360}$&                             &$\frac{15}{360360}$&&&\\
     \hline
     8&$\frac{-63}{12252240}$&$\frac{-7}{720720}$&&&&&&&\\
     \hline
    \end{tabular}
\end{table}

\begin{table}[htbp]
  
  \caption{The coefficients $a_l$ in Eq. (\ref{Eq:ReconstructedFlux}) for computing $\hat{f}_{j+1/2}^+$ on different candidate stencils.}
  \label{Table:Flux_reconstruction2}
  \centering
  \begin{tabular}{|c|c|c|c|c|c|c|c|c|c|c|}
    % after \\: \hline or \cline{col1-col2} \cline{col3-col4} ...
     \hline
     $l$&$S_{j-5}^{j+6}$&$S_{j-6}^{j+5}$&$S_{j-5}^{j+5}$&$S_{j-6}^{j+4}$&$S_{j-4}^{j+5}$&$S_{j-5}^{j+4}$&$S_{j-4}^{j+4}$&$S_{j-5}^{j+3}$&$S_{j-3}^{j+4}$&$S_{j-4}^{j+3}$\\
     \hline
     -6&                  &$\frac{5}{27720}$&                    &$\frac{12}{27720}$&&&&&&\\
     \hline
     -5&$\frac{-5}{27720}$&$\frac{-65}{27720}$&$\frac{-10}{27720}$&$\frac{-142}{27720}$&               &$\frac{-2}{2520}$&&$\frac{-5}{2520}$&&\\
     \hline
     -4&$\frac{67}{27720}$&$\frac{397}{27720}$&$\frac{122}{27720}$&$\frac{782}{27720}$&$\frac{2}{2520}$&$\frac{22}{2520}$&$\frac{4}{2520}$&$\frac{49}{2520}$&                              &$\frac{3}{840}$\\
     \hline
     -3&$\frac{-428}{27720}$&$\frac{-1528}{27720}$&$\frac{-703}{27720}$&$\frac{-2683}{27720}$&$\frac{-23}{2520}$&$\frac{-113}{2520}$&$\frac{-41}{2520}$&$\frac{-221}{2520}$&$\frac{-3}{840}$&$\frac{-27}{840}$\\
     \hline
     -2&$\frac{1772}{27720}$&$\frac{4247}{27720}$&$\frac{2597}{27720}$&$\frac{6557}{27720}$&$\frac{127}{2520}$&$\frac{367}{2520}$&$\frac{199}{2520}$&$\frac{619}{2520}$&$\frac{29}{840}$&$\frac{113}{840}$\\
     \hline
     -1&$\frac{-5653}{27720}$&$\frac{-9613}{27720}$&$\frac{-7303}{27720}$&$\frac{-12847}{27720}$&$\frac{-473}{2520}$&$\frac{-893}{2520}$&$\frac{-641}{2520}$&$\frac{-1271}{2520}$&$\frac{-139}{840}$&$\frac{-307}{840}$\\
     \hline
     0&$\frac{18107}{27720}$&$\frac{22727}{27720}$&$\frac{20417}{27720}$&$\frac{25961}{27720}$&$\frac{1627}{2520}$&$\frac{2131}{2520}$&$\frac{1879}{2520}$&$\frac{2509}{2520}$&$\frac{533}{840}$&$\frac{743}{840}$\\
     \hline
     1&$\frac{18107}{27720}$&$\frac{14147}{27720}$&$\frac{15797}{27720}$&$\frac{11837}{27720}$&$\frac{1627}{2520}$&$\frac{1207}{2520}$&$\frac{1375}{2520}$&$\frac{955}{2520}$&$\frac{533}{840}$&$\frac{365}{840}$\\
     \hline
     2&$\frac{-5653}{27720}$&$\frac{-3178}{27720}$&$\frac{-4003}{27720}$&$\frac{-2023}{27720}$&$\frac{-473}{2520}$&$\frac{-233}{2520}$&$\frac{-305}{2520}$&$\frac{-125}{2520}$&$\frac{-139}{840}$&$\frac{-55}{840}$\\
     \hline
     3&$\frac{1772}{27720}$&$\frac{672}{27720}$&$\frac{947}{27720}$&$\frac{287}{27720}$&$\frac{127}{2520}$&$\frac{37}{2520}$&$\frac{55}{2520}$&$\frac{10}{2520}$&$\frac{29}{840}$&$\frac{5}{840}$\\
     \hline
     4&$\frac{-428}{27720}$&$\frac{-98}{27720}$&$\frac{-153}{27720}$&$\frac{-21}{27720}$&$\frac{-23}{2520}$&$\frac{-3}{2520}$&$\frac{-5}{2520}$&&$\frac{-3}{840}$&\\
     \hline
     5&$\frac{67}{27720}$&$\frac{7}{27720}$&$\frac{12}{27720}$&&$\frac{2}{2520}$&&&&&\\
     \hline
     6&$\frac{-5}{27720}$&&&&&&&&&\\
     \hline
     $l$&$S_{j-3}^{j+3}$&$S_{j-2}^{j+3}$&$S_{j-3}^{j+2}$&$S_{j-2}^{j+2}$&$S_{j-1}^{j+2}$&$S_{j-2}^{j+1}$&$S_{j-1}^{j+1}$&$S_{j}^{j+1}$&$S_{j-1}^{j}$&$S_{j}^{j}$\\
     \hline
     -3&$\frac{-3}{420}$&              &$\frac{-1}{60}$&&&&&&&\\
     \hline
     -2&$\frac{25}{420}$&$\frac{1}{60}$&$\frac{7}{60}$&$\frac{2}{60}$&&$\frac{1}{12}$&&&&\\
     \hline
     -1&$\frac{-101}{420}$&$\frac{-8}{60}$&$\frac{-23}{60}$&$\frac{-13}{60}$&$\frac{-1}{12}$&$\frac{-5}{12}$&$\frac{-1}{6}$&&$\frac{-1}{2}$&\\
     \hline
     0&$\frac{319}{420}$&$\frac{37}{60}$&$\frac{57}{60}$&$\frac{47}{60}$&$\frac{7}{12}$&$\frac{13}{12}$&$\frac{5}{6}$&$\frac{1}{2}$&$\frac{3}{2}$&1\\
     \hline
     1&$\frac{214}{420}$&$\frac{37}{60}$&$\frac{22}{60}$&$\frac{27}{60}$&$\frac{7}{12}$&$\frac{3}{12}$&$\frac{2}{6}$&$\frac{1}{2}$&&\\
     \hline
     2&$\frac{-38}{420}$&$\frac{-8}{60}$&$\frac{-2}{60}$&$\frac{-3}{60}$&$\frac{-1}{12}$&&&&&\\
     \hline
     3&$\frac{4}{420}$&$\frac{1}{60}$&&&&&&&&\\
     \hline
    \end{tabular}
\end{table}

%%%%%%%%%%%%%%%%%%%%%%%%%%%%%%%%%%%%%%%%%%%%%%%%%%%%%%%%%%%%%%
\subsection{Efficient ENO reconstructions with multi-resolution}\label{SEC:ENO_MR_Rec}

WENO-AO reconstructions \cite{Levy2000CompactCWENO,Zhu2016NewWENOFD, Zhu2018NewMultiResolutionWENO, Zhu2019NewMultiResolutionWENO_Tri, 
Zhu2020NewMultiResolutionWENO_Tet, Balsara2016WENO_AO, Balsara2020WENO_AO_Unstructure}
depend on the smoothness indicators of Jiang and Shu \cite{Jiang_Shu1996WENO} defined by
Eq. (\ref{Eq:Smoothness_Indicator}) which becomes very complicated and computationally expensive,
so very high-order WENO-AO reconstructions with many levels are not practical.
ENO-AO reconstructions of Shen \cite{Shen2023ENO_AO} have much simpler
smoothness indicators as defined by Eq. (\ref{Eq:New_Smoothness_Indicator}),
but they need to select the smoothest stencil from the pool of candidates.
When the number of candidates becomes very large,
the selection procedure also becomes cumbersome.
In this study, we propose a class of ENO reconstructions with multi-resolution (ENO-MR)
that are efficient in selecting the stencil with optimal order from many candidates.
The selection procedure is as follows:
\begin{itemize}
  \item [Step 1.] We define a baseline smoothness indicator as
  \begin{subequations}\label{Eq:Baseline_IS}
    \begin{equation}\label{Eq:Baseline_IS0}
       IS_0=\text{MIN}\left(IS_L,IS_R\right),
    \end{equation}
    \text{with}
    \begin{equation}\label{Eq:Baseline_ISL}
      IS_L=\text{MAX}\left(\left|f_j-f_{j-1}\right|,\left|f_j-2f_{j-1}+f_{j-2}\right|\right),
   \end{equation}
   \begin{equation}\label{Eq:Baseline_ISR}
      IS_R=\text{MAX}\left(\left|f_j-f_{j+1}\right|,\left|f_j-2f_{j+1}+f_{j+2}\right|\right).
 \end{equation}
\end{subequations}
\item [Step 2.] We define smoothness indicators for $S_{j-m}^{j+n}$ as
\begin{equation}\label{Eq:CurrentIS}
  IS_{j-m}^{j+n}=\left|\frac{d^{m+n}P_{j-m}^{j+n}(x)}{dx^{m+n}}\right|h^{m+n},
\end{equation}
 where $P_{j-m}^{j+n}(x)$ is the polynomial reconstructed on $S_{j-m}^{j+n}$ as Eq. (\ref{Eq:Polynomial_Construction}).
 Once $P_{j-m}^{j+n}(x)$ is determined, $IS_{j-m}^{j+n}$ can be calculated in terms of $f_k$ as
 \begin{equation}\label{Eq:ISCoe}
  IS_{j-m}^{j+n}=\left|\sum_{l=-m}^{n} b_lf_{j+l}\right|.
\end{equation}
The coefficients $b_l$ of all candidate stencils are given by Tables \ref{Table:ISCoe1} and \ref{Table:ISCoe2}.
 \item [Step 3.] We compare the smoothness indicators of candidate stencils 
 in sequence as listed by Eq. (\ref{Eq:Candidates}) (from high-order to low-order) with the baseline smoothness indicator.
 \begin{itemize}
  \item [Step 3.1.] If any $IS_{j-m}^{j+n}$ ($m\ge1$ and $n\ge1$) is smaller than the baseline $IS_0$,
  we directly use the reconstructed flux on $S_{j-m}^{j+n}$. 

  \item[Step 3.2.] If all $IS_{j-m}^{j+n}$ ($m\ge1$ and $n\ge1$) are larger than the baseline $IS_0$,
 we use the minmod function to select a low-order stencil from $\{S_j^{j+1}$, $S_{j-1}^{j},S_j^j\}$, i.e.,
 \begin{subequations}
  \begin{equation}
    f_{j+1/2}=f_j+\frac{1}{2}\text{minmod}\left(f_{j+1}-f_j,f_j-f_{j-1}\right),
   \end{equation}
   \text{where the minmod function is defined as}
   \begin{equation}
    \text{minmod}(a,b)=\begin{cases}
      a,\quad \text{if } ab>0 \text{ and } |a|\le|b|,\\
      b,\quad \text{if } ab>0 \text{ and } |a|>|b|,\\
      0,\quad \text{if } ab\le0.
    \end{cases}
   \end{equation}
 \end{subequations}
 \end{itemize}
\end{itemize}

\begin{table}[htbp]
  
  \caption{The coefficients $b_l$ in Eq. (\ref{Eq:ISCoe}) for computing smoothness indicators on different candidate stencils.}
  \label{Table:ISCoe1}
  \centering
  \begin{tabular}{|c|c|c|c|c|c|c|c|c|c|}
    % after \\: \hline or \cline{col1-col2} \cline{col3-col4} ...
     \hline
      $l$&$S_{j-8}^{j+8}$&$S_{j-7}^{j+8}$&$S_{j-8}^{j+7}$&$S_{j-7}^{j+7}$&$S_{j-8}^{j+6}$&$S_{j-6}^{j+7}$&$S_{j-7}^{j+6}$&$S_{j-6}^{j+6}$&$S_{j-7}^{j+5}$  \\
      \hline
      -8&1&    &-1& &1&&&&\\
      \hline
      -7&-16&-1&15&1&-14&  &-1&    &1\\
     \hline
     -6&120&15&-105&-14&91&-1&13&1&-12\\
     \hline
     -5&-560&-105&455&91&-364&13&-78&-12&66\\
     \hline
     -4&1820&455&-1365&-364&1001&-78&286&66&-220\\
     \hline
     -3&-4368&-1365&3003&1001&-2002&286&-715&-220&495\\
     \hline
     -2&8008&3003&-5005&-2002&3003&-715&1287&495&-792\\
     \hline
     -1&-11440&-5005&6435&3003&-3423&1287&-1716&-792&924\\
     \hline
     0&12870&6435&-6435&-3423&3003&-1716&1716&924&-792\\
     \hline
      1&-11440&-6435&5005&3003&-2002&1716&-1287&-792&495\\
      \hline
      2&8008&5005&-3003&-2002&1001&-1287&715&495&-220\\
     \hline
     3&-4368&-3003&1365&1001&-364&715&-286&-220&66\\
     \hline
     4&1820&1365&-455&-364&91&-286&78&66&-12\\
     \hline
     5&-560&-455&105&91&-14&78&-13&-12&1\\
     \hline
     6&120&105&-15&-14&1&-13&1&1&\\
     \hline
     7&-16&-15&1&1&    &1&&&\\
     \hline
     8&1&1&&&&&&&\\
     \hline
    \end{tabular}
\end{table}

\begin{table}[htbp]
  
  \caption{The coefficients $b_l$ in Eq. (\ref{Eq:ISCoe}) for computing smoothness indicators on different candidate stencils.}
  \label{Table:ISCoe2}
  \centering
  \begin{tabular}{|c|c|c|c|c|c|c|c|c|c|c|}
    % after \\: \hline or \cline{col1-col2} \cline{col3-col4} ...
     \hline
     $l$&$S_{j-5}^{j+6}$&$S_{j-6}^{j+5}$&$S_{j-5}^{j+5}$&$S_{j-6}^{j+4}$&$S_{j-4}^{j+5}$&$S_{j-5}^{j+4}$&$S_{j-4}^{j+4}$&$S_{j-5}^{j+3}$&$S_{j-3}^{j+4}$&$S_{j-4}^{j+3}$\\
     \hline
     -6&   &-1& &1&&&&&&\\
     \hline
     -5&-1&11&1&-10&    &-1&&1&&\\
     \hline
     -4&11&-55&-10&45&-1&9&1&-8&      &-1\\
     \hline
     -3&-55&165&45&-120&9&-36&-8&28&-1&7\\
     \hline
     -2&165&-330&-120&210&-36&84&28&-56&7&-21\\
     \hline
     -1&-330&462&210&-252&84&-126&-56&70&-21&35\\
     \hline
     0&462&-462&-252&210&-126&126&70&-56&35&-35\\
     \hline
     1&-462&330&210&-120&126&-84&-56&28&-35&21\\
     \hline
     2&330&-165&-120&45&-84&36&28&-8&21&-7\\
     \hline
     3&-165&55&45&-10&36&-9&-8&1&-7&1\\
     \hline
     4&55&-11&-10&1&-9&1&1&&1&\\
     \hline
     5&-11&1&1&&1&&&&&\\
     \hline
     6&1&&&&&&&&&\\
     \hline
     $l$&$S_{j-3}^{j+3}$&$S_{j-2}^{j+3}$&$S_{j-3}^{j+2}$&$S_{j-2}^{j+2}$&$S_{j-1}^{j+2}$&$S_{j-2}^{j+1}$&$S_{j-1}^{j+1}$&$S_{j}^{j+1}$&$S_{j-1}^{j}$&$S_{j}^{j}$\\
     \hline
     -3&1&   &-1&&&&&&&\\
     \hline
     -2&-6&-1&5&1&&-1&&&&\\
     \hline
     -1&15&5&-10&-4&-1&3&1&&-1&\\
     \hline
     0&-20&-10&10&6&3&-3&-2&-1&1&\\
     \hline
     1&15&10&-5&-4&-3&1&1&1&&\\
     \hline
     2&-6&-5&1&1&1&&&&&\\
     \hline
     3&1&1&&&&&&&&\\
     \hline
    \end{tabular}
\end{table}

\subsection{Properties of ENO-MR reconstructions}
\begin{definition}
  For a smooth function $f(x)$, 
  if $f^{'}(x_c)=f^{''}(x_c)...=f^{(n_c)}(x_c)=0$ but $f^{(n_c+1)}(x_c)\neq 0$ $(n_c\ge 1)$,
  we call $x_c$ is a critical point of order $n_c$.
\end{definition}

\begin{definition}
  Assume $\mathcal{R}$ is a flux reconstruction operator, i.e., $f_{j+1/2}=\mathcal{R}(f_{j-m},...,f_{j+n})$.
  If $\mathcal{R}(\lambda f_{j-m},...,\lambda f_{j+n})=\lambda\mathcal{R}(f_{j-m},...,f_{j+n})$
  for any given function $f$ and $\lambda\in\mathbb{R}$, we call $\mathcal{R}$ is a scale-invariant operator.
\end{definition}

We demonstrate the properties of ENO-MR reconstructions in the following three aspects:

First, we show that ENO-MR reconstructions on the global stencil $S_{j-r+1}^{j+r-1}$ 
achieve optimal convergence order of $2r-1$
if $f$ is smooth on $S_{j-r+1}^{j+r-1}$ and contains critical points at most of order $2r-2$.
Perform Taylor expansion at $x_j$, we can write 
Eqs. (\ref{Eq:Baseline_ISL}), (\ref{Eq:Baseline_ISR}), and (\ref{Eq:CurrentIS})
in the following forms:
\begin{subequations}\label{Eq:IS_Taylor}
  \begin{equation}\label{Eq:ISL_Taylor}
    IS_L=\text{MAX}\left(\left|\sum_{k=1}^{\infty}\frac{(-1)^{k+1}}{k!}f^{(k)}_jh^k\right|,
    \left|\sum_{k=2}^{\infty}\frac{(-1)^k(2^k-2)}{k!}f^{(k)}_jh^k\right|\right),
 \end{equation}
 \begin{equation}\label{Eq:ISR_Taylor}
    IS_R=\text{MAX}\left(\left|\sum_{k=1}^{\infty}\frac{1}{k!}f^{(k)}_jh^k\right|,
    \left|\sum_{k=2}^{\infty}\frac{(2^k-2)}{k!}f^{(k)}_jh^k\right|\right),
\end{equation}
\begin{equation}\label{Eq:ISr_Taylor}
  IS_{j-r+1}^{j+r-1}=\left|f^{(2r-2)}_jh^{2r-2}+\mathcal{O}(h^{2r-1})\right|.
\end{equation}
\end{subequations}

If $x_j$ is a critical point of order $n_c<(2r-2)$, Eqs. (\ref{Eq:Baseline_IS0}) and (\ref{Eq:IS_Taylor}) indicate that
\begin{equation}
  IS_0=\mathcal{O}(h^{n_c+1}),\quad IS_{j-r+1}^{j+r-1}=\mathcal{O}(h^{2r-2}).
\end{equation}

When $n_c\le(2r-4)$, we have $IS_{j-r+1}^{j+r-1}<IS_0$ 
providing the mesh size $h$ is sufficiently small.
As a result, ENO-MR reconstructions described in Section \ref{SEC:ENO_MR_Rec}
can correctly select the optimal stencil $S_{j-r+1}^{j+r-1}$ in these cases.
In addition, when $IS_{j-r+1}^{j+r-1}<IS_0$, 
we directly adopt $S_{j-r+1}^{j+r-1}$ to reconstruct the flux,
and do not need to calculate the smoothness of other candidate stencils.
Therefore, ENO-MR reconstructions are efficient even if there are many candidate stencils
because the proportion of smooth regions is usually much greater than that of discontinuities in practice.

Next, we show that ENO-MR reconstructions can exclude all the stencils 
containing discontinuities if $f$ has one discontinuity on $S_{j-r+1}^{j+r-1}$ 
which is a common case in practice.
Eqs. (\ref{Eq:Baseline_ISL}) and (\ref{Eq:Baseline_ISR}) indicate that
$IS_L$ and $IS_R$ respectively measure the smoothness of $S_{j-2}^j$ and $S_{j}^{j+2}$.
Since $f$ has only one discontinuity on $S_{j-r+1}^{j+r-1}$,
at most one of $S_{j-2}^j$ and $S_{j}^{j+2}$ contains the discontinuity.
As a result, $IS_0=\mathcal{O}(h)$ regardless of the location of the discontinuity.
Meanwhile, if a candidate stencil $S_{j-m}^{j+n}$ contains a discontinuity,
then $IS_{j-m}^{j+n}=\mathcal{O}(\Delta f)$ where $\Delta f$ is the jump of $f$ at the discontinuity.
For example, when the discontinuity is located in $(x_{j-1},x_{j})$ and $h\ll|\Delta f|$, we have
\begin{equation}
  IS_{j-1}^{j+1}=|f_{j-1}-2f_j+f_{j+1}|=|\underbrace{f_{j-1}-f_j}_{\Delta f}
                                    \underbrace{-f_j+f_{j+1}}_{\mathcal{O}(h)}|\approx |\Delta f|>IS_0.
\end{equation}
Therefore, all candidate stencils containing the discontinuity are abandoned by
the ENO-MR selection strategy, thereby achieving ENO properties at discontinuities.

Finally, we show that ENO-MR reconstructions are essentially scale-invariant operators.
Most WENO schemes are scale-dependent due to the sensitive parameter $\varepsilon$
appearing in the denominator of nonlinear weights 
(see Eqs. (\ref{Eq:WENO_JS_Weight}), (\ref{Eq:WENO_Z_Weight}), and (\ref{Eq:WENO_AO_Weights})).
The reason behind this is that the smoothness indicators $\beta$ depend on the scale of $f$.
More specifically, $\beta(\lambda f)=\lambda^2\beta(f)$ according to Eq. (\ref{Eq:Smoothness_Indicator}).
If we use a fixed sensitive parameter $\varepsilon$, nonlinear weights 
defined by Eqs. (\ref{Eq:WENO_JS_Weight}), (\ref{Eq:WENO_Z_Weight}), and (\ref{Eq:WENO_AO_Weights})
apparently depend on $\lambda$.
To fix this issue, we need to modify the nonlinear weights properly \cite{Don2022WENO_Si}.
However, ENO-MR reconstructions are completely parameter-free and essentially scale-invariant.
The new smoothness indicators are absolute values of linear combinations of the point values of $f$, 
so we have 
\begin{equation}
  IS_0(\lambda f)=|\lambda| IS_0(f),\quad IS_{j-m}^{j+n}(\lambda f)=|\lambda| IS_{j-m}^{j+n}(f).
\end{equation}
Obviously, when the scale of $f$ changes, 
the inequality relationships between $IS_0$ and $IS_{j-m}^{j+n}$ are invariant.
Therefore, ENO-MR reconstructions will select the same stencil for flux reconstruction
when the scale of $f$ changes.
That is to say, ENO-MR reconstructions are scale-invariant.

 %%%%%%%%%%%%%%%%%%%%%%%%%%%%%%%%%%%%%%%%%%
\section{Numerical examples}\label{SEC:_NumExam}
In this section, we use some benchmarks to test the performance of the proposed 
5th-order ($r=3$), 9th-order ($r=5$), 13th-order ($r=7$), and 17th-order ($r=9$) ENO-MR schemes.
The results of WENO-AO(5,3) and WENO-AO(9,5,3) schemes \cite{Balsara2016WENO_AO}
are used as references.
Following Balsara \emph{et al.} \cite{Balsara2016WENO_AO}, 
we use WENO-Z weights with $\varepsilon=10^{-12}$, $p=2$, and $\gamma_{Lo}=\gamma_{Hi}=0.85$ for WENO-AO schemes.
For the linear advection equations,
time is discretized by the $\ell$SSP-RK($m,m-1$) scheme with $\Delta t=h$,
where the order of time $m-1$ equals the order of ENO-MR space reconstruction.
For the other equations,
time is discretized by the SSP-RK3 scheme.
If there is no special instruction, we set $CFL=0.3$.
Quadruple precision is used for convergence tests.
For computations of the Euler equations,
the Roe average at the cell face is adopted for characteristic decomposition.
%%%%%%%%%%%%%%%%%%%%%%%%%%%%%%%%%%%%%%%%%%%%%%%%%%%%%%%%%%%%%%%%%%
\subsection{Numerical examples for the 1D linear advection equation}
We consider the 1D linear advection equation
\begin{equation}\label{Eq:1DAdvEq}
 \frac{\partial u}{\partial t}+\frac{\partial u}{\partial x}=0.
\end{equation}

First, we test the convergence rate by setting 
the initial condition as $u_0(x)=\lambda sin^\alpha(\pi x)$, $x\in[-1,1]$.
It is easy to check that $x=-1, 0, 1$ are the $(\alpha-1)$th-order critical points of $\lambda sin^\alpha(\pi x)$.
Periodic boundary conditions are applied to left and right boundaries.
Tables. \ref{Table:1DAdvectionConvergence5th}-\ref{Table:1DAdvectionConvergence13th17th}
shows the numerical errors of WENO-AO and ENO-MR schemes at $t=2$,
where $L_1$ is the average error and $L_\infty$ is the maximum error.
When $\lambda=1$, WENO-AO schemes can not achieve optimal convergence rates on relatively coarse meshes
and achieve `super convergence' when the mesh size $h$ reduces to some thresholds ($1/800$ for $\alpha=3$).
This is because the order of smoothness indicators $\beta_k^{(r)}$ in Eq. (\ref{Eq:WENO_Z_Weight}) will increase at critical points
but are still much larger than the sensitive parameter $\varepsilon$ when $h$ is large.
As a result, the order of the difference between WENO-Z nonlinear weights and optimal will reduce.
When $h$ reduces to some thresholds so that $\beta_k^{(r)}$ and $\tau_{2r-1}$ become much smaller than $\varepsilon$,
then WENO-Z nonlinear weights are almost equal to optimal weights.
Therefore, WENO-AO schemes achieve `super convergence'.
When $\lambda$ increases to $10^6$,  $\beta_k^{(r)}$ and $\tau_{2r-1}$ significantly increase 
but $\varepsilon$ remains fixed, so WENO-AO schemes require smaller $h$ to achieve `super convergence'.
However, when $\lambda=10^6$, the numerical errors of ENO-MR schemes are exactly scaled by a factor of $10^6$
comparing with the results of $\lambda=1$.
This demonstrates the scale-invariant property of ENO-MR schemes.
When $\alpha=3$, the solution contains 2nd-order critical points,
all ENO-MR schemes can achieve optimal order.
When $\alpha=4$, the solution contains 3rd-order critical points,
all ENO-MR schemes except ENO-MR5 can achieve optimal order.
This is consistent with the previous analysis.

\begin{table}[htbp]
 \footnotesize
 \caption{Numerical errors of the advection of $\lambda sin^\alpha(\pi x)$ at $t=2$ computed by 5th-order WENO-AO and ENO-MR schemes.}
 \label{Table:1DAdvectionConvergence5th}
 \centering
 \begin{tabular}{|c|c|c|c|c|c|c|c|c|}
   \hline
   % after \\: \hline or \cline{col1-col2} \cline{col3-col4} ...
     & \multicolumn{4}{|c|}{WENO-AO(5,3),$\lambda=1,\alpha=3$} & \multicolumn{4}{|c|}{WENO-AO(5,3),$\lambda=1,\alpha=4$} \\
    \hline
     $h$       & $L_1$  &Order & $L_\infty$  &Order    & $L_1$  &Order & $L_\infty$  &Order    \\
    \hline
    1/100              &2.20E-05 &-     &3.05E-04      &-        &2.52E-06 &-     &2.67E-05      &- \\

    1/200             &2.38E-06 &3.21  &5.86E-05     &2.38     &9.46E-08 &4.74  &1.90E-06      &3.81 \\

    1/400             &2.14E-07 &3.48  &9.60E-06      &2.61     &9.81E-10 &6.59  &2.91E-09      &9.35 \\

    1/800             &6.97E-11 &11.59  &5.41E-09      &10.79     &3.06E-11 &5.00  &5.10E-11    &5.83 \\

    1/1600             &3.49E-13 &7.64  &1.08E-12      &12.29     &9.56E-13 &5.00  &1.59E-12      &5.00 \\
    \hline

    & \multicolumn{4}{|c|}{WENO-AO(5,3),$\lambda=10^6,\alpha=3$} & \multicolumn{4}{|c|}{WENO-AO(5,3),$\lambda=10^6,\alpha=4$}  \\
    \hline
     $h$       & $L_1$  &Order & $L_\infty$  &Order    & $L_1$  &Order & $L_\infty$  &Order    \\
    \hline
    1/100              &2.20E+01 &-     &3.05E+02     &-        &2.53E+00 &-     &2.68E+01      &- \\

    1/200             &2.38E+00 &3.21  &5.87E+01     &2.38     &1.31E-01 &4.27  &2.58E+00      &3.38 \\

    1/400             &2.63E-01 &3.18  &1.14E+01      &2.36     &7.24E-03 &4.18  &2.69E-01      &3.26 \\

    1/800             &2.66E-02 &3.31  &2.05E+00      &2.48     &4.74E-04 &3.93  &3.68E-02   &2.87 \\

    1/1600             &3.17E-03 &3.07  &4.58E-01      &2.16     &3.61E-05 &3.72  &5.68E-03      &2.70 \\
    \hline

    & \multicolumn{4}{|c|}{ENO-MR5,$\lambda=1,\alpha=3$} & \multicolumn{4}{|c|}{ENO-MR5,$\lambda=1,\alpha=4$} \\
    \hline
     $h$       & $L_1$  &Order & $L_\infty$  &Order    & $L_1$  &Order & $L_\infty$  &Order    \\
    \hline
    1/100              &3.54E-07 &-     &5.61E-07      &-        &2.71E-06 &-     &2.13E-05      &- \\

    1/200             &1.11E-08 &5.00  &1.76E-08     &4.99     &1.00E-07 &4.76  &1.76E-06      &3.60 \\

    1/400             &3.48E-10 &5.00  &5.49E-10      &5.00     &3.81E-09 &4.71  &1.29E-07      &3.77 \\

    1/800             &1.09E-11 &5.00  &1.72E-11      &5.00     &1.48E-10 &4.69  &9.95E-09    &3.70\\

    1/1600             &3.40E-13 &5.00  &5.36E-13     &5.00     &5.65E-12 &4.71  &7.36E-10      &3.76 \\
    \hline

    & \multicolumn{4}{|c|}{ENO-MR5,$\lambda=10^6,\alpha=3$} & \multicolumn{4}{|c|}{ENO-MR5,$\lambda=10^6,\alpha=4$}  \\
    \hline
     $h$       & $L_1$  &Order & $L_\infty$  &Order    & $L_1$  &Order & $L_\infty$  &Order    \\
    \hline
    1/100              &3.54E-01 &-     &5.61E-01      &-        &2.71E+00 &-     &2.13E+01      &- \\

    1/200             &1.11E-02 &5.00  &1.76E-02     &4.99     &1.00E-01 &4.76  &1.76E+00      &3.60 \\

    1/400             &3.48E-04 &5.00  &5.49E-04      &5.00     &3.81E-03 &4.71  &1.29E-01      &3.77 \\

    1/800             &1.09E-05 &5.00  &1.72E-05      &5.00     &1.48E-04 &4.69  &9.95E-03    &3.70\\

    1/1600             &3.40E-07 &5.00  &5.36E-07     &5.00     &5.65E-06 &4.71  &7.36E-04      &3.76 \\
    \hline
 \end{tabular}
\end{table}

\begin{table}[htbp]
  \footnotesize
  \caption{Numerical errors of the advection of $\lambda sin^\alpha(\pi x)$ at $t=2$ computed by 9th-order WENO-AO and ENO-MR schemes.}
  \label{Table:1DAdvectionConvergence9th}
  \centering
  \begin{tabular}{|c|c|c|c|c|c|c|c|c|}
    \hline
    % after \\: \hline or \cline{col1-col2} \cline{col3-col4} ...
      & \multicolumn{4}{|c|}{WENO-AO(9,5,3),$\lambda=1,\alpha=3$} & \multicolumn{4}{|c|}{WENO-AO(9,5,3),$\lambda=1,\alpha=4$} \\
     \hline
      $h$       & $L_1$  &Order & $L_\infty$  &Order    & $L_1$  &Order & $L_\infty$  &Order    \\
     \hline
     1/50              &1.70E-04 &-     &1.37E-03      &-        &4.79E-05 &-     &2.67E-04      &- \\
 
     1/100             &1.43E-05 &3.57  &2.13E-04     &2.69     &2.30E-06 &4.38  &2.40E-05      &3.48 \\
 
     1/200             &1.17E-06 &3.61  &3.33E-05      &2.68     &1.00E-07 &4.52  &1.95E-06      &3.62 \\
 
     1/400             &8.39E-08 &3.80  &4.70E-06      &2.83     &6.13E-10 &7.35  &2.80E-08    &6.12 \\
 
     1/800             &1.05E-09 &6.32  &1.25E-7      &5.23     &1.35E-13 &12.15  &1.28E-11      &11.10 \\
     \hline
 
     & \multicolumn{4}{|c|}{WENO-AO(9,5,3),$\lambda=10^6,\alpha=3$} & \multicolumn{4}{|c|}{WENO-AO(9,5,3),$\lambda=10^6,\alpha=4$}  \\
     \hline
      $h$       & $L_1$  &Order & $L_\infty$  &Order    & $L_1$  &Order & $L_\infty$  &Order    \\
     \hline
     1/50              &1.70E+02 &-     &1.37E+03     &-        &4.79E+01 &-     &2.67E+02      &- \\
 
     1/100             &1.43E+01 &3.57  &2.13E+02     &2.69     &2.30E+00 &4.38  &2.40E+01      &3.48 \\
 
     1/200             &1.17E+00 &3.61  &3.34E+01      &2.67     &1.06E-01 &4.44  &2.08E+00      &3.53 \\
 
     1/400             &9.33E-02 &3.65  &5.17E+00      &2.69     &4.86E-03 &4.45  &1.84E-01   &3.50 \\
 
     1/800             &7.58E-03 &3.62  &8.29E-01      &2.64     &2.30E-04 &4.40  &1.76E-02      &3.39 \\
     \hline

     & \multicolumn{4}{|c|}{ENO-MR9,$\lambda=1,\alpha=3$} & \multicolumn{4}{|c|}{ENO-MR9,$\lambda=1,\alpha=4$} \\
     \hline
      $h$       & $L_1$  &Order & $L_\infty$  &Order    & $L_1$  &Order & $L_\infty$  &Order    \\
     \hline
     1/50              &7.01E-10 &-     &1.11E-09      &-        &6.26E-09 &-     &9.78E-09      &- \\
 
     1/100             &1.39E-12 &8.98  &2.19E-12     &8.99     &1.24E-11 &8.98  &1.94E-11      &8.98 \\
 
     1/200             &2.72E-15 &9.00  &4.28E-15      &9.00     &2.42E-14 &9.00  &3.81E-14      &8.99 \\
 
     1/400             &5.32E-18 &9.00  &8.37E-18      &9.00     &4.73E-17 &9.00  &7.46E-17    &9.00\\
 
     1/800             &1.04E-20 &9.00  &1.63E-20     &9.01     &9.24E-20 &9.00  &1.46E-19      &9.00 \\
     \hline
 
     & \multicolumn{4}{|c|}{ENO-MR9,$\lambda=10^6,\alpha=3$} & \multicolumn{4}{|c|}{ENO-MR9,$\lambda=10^6,\alpha=4$}  \\
     \hline
      $h$       & $L_1$  &Order & $L_\infty$  &Order    & $L_1$  &Order & $L_\infty$  &Order    \\
     \hline
     1/50              &7.01E-04 &-     &1.11E-03      &-        &6.26E-03 &-     &9.78E-03      &- \\
 
     1/100             &1.39E-06 &8.98  &2.19E-06     &8.99     &1.24E-05 &8.98  &1.94E-05      &8.98 \\
 
     1/200             &2.72E-09 &9.00  &4.28E-09      &9.00     &2.42E-08 &9.00  &3.81E-08      &8.99 \\
 
     1/400             &5.32E-12 &9.00  &8.37E-12      &9.00     &4.73E-11 &9.00  &7.46E-11    &9.00\\
 
     1/800             &1.04E-14 &9.00  &1.63E-14     &9.01     &9.24E-14 &9.00  &1.46E-13      &9.00 \\
     \hline
  \end{tabular}
 \end{table}

 \begin{table}[htbp]
  \footnotesize
  \caption{Numerical errors of the advection of $\lambda sin^\alpha(\pi x)$ at $t=2$ computed by 13th- and 17th-order ENO-MR schemes.}
  \label{Table:1DAdvectionConvergence13th17th}
  \centering
  \begin{tabular}{|c|c|c|c|c|c|c|c|c|}
    \hline
    % after \\: \hline or \cline{col1-col2} \cline{col3-col4} ...
      & \multicolumn{4}{|c|}{ENO-MR13,$\lambda=1,\alpha=3$} & \multicolumn{4}{|c|}{ENO-MR13,$\lambda=1,\alpha=4$} \\
     \hline
      $h$       & $L_1$  &Order & $L_\infty$  &Order    & $L_1$  &Order & $L_\infty$  &Order    \\
     \hline
     1/50              &4.62E-14 &-     &7.31E-14      &-        &1.30E-12 &-     &2.03E-12      &- \\
 
     1/100             &5.73E-18 &12.98  &9.04E-18     &12.98     &1.61E-16 &12.98  &2.53E-16      &12.97 \\
 
     1/200             &7.03E-22 &12.99  &1.11E-21      &12.99     &1.98E-20 &12.99  &3.10E-20      &13.00 \\
 
     1/400             &8.60E-26 &13.00  &1.35E-25      &13.01    &2.42E-24 &13.00  &3.80E-24    &13.00 \\
 
     1/800             &1.05E-29 &13.00  &1.66E-29      &12.99     &2.95E-28 &13.00  &4.64E-28      &13.00 \\
     \hline
 
     & \multicolumn{4}{|c|}{ENO-MR13,$\lambda=10^6,\alpha=3$} & \multicolumn{4}{|c|}{ENO-MR13,$\lambda=10^6,\alpha=4$}  \\
     \hline
      $h$       & $L_1$  &Order & $L_\infty$  &Order    & $L_1$  &Order & $L_\infty$  &Order    \\
     \hline
     \hline
     1/50              &4.62E-08 &-     &7.31E-08      &-        &1.30E-06 &-     &2.03E-06      &- \\
 
     1/100             &5.73E-12 &12.98  &9.04E-12     &12.98     &1.61E-10 &12.98  &2.53E-10      &12.97 \\
 
     1/200             &7.03E-16 &12.99  &1.11E-15      &12.99     &1.98E-14 &12.99  &3.10E-14      &13.00 \\
 
     1/400             &8.60E-20 &13.00  &1.35E-19      &13.01    &2.42E-18 &13.00  &3.80E-18    &13.00 \\
 
     1/800             &1.05E-23 &13.00  &1.66E-23      &12.99     &2.95E-22 &13.00  &4.64E-22      &13.00 \\
     \hline

     & \multicolumn{4}{|c|}{ENO-MR17,$\lambda=1,\alpha=3$} & \multicolumn{4}{|c|}{ENO-MR17,$\lambda=1,\alpha=4$} \\
     \hline
      $h$       & $L_1$  &Order & $L_\infty$  &Order    & $L_1$  &Order & $L_\infty$  &Order    \\
     \hline
     1/15              &1.97E-09 &-     &3.10E-09      &-        &1.60E-07 &-     &3.24E-07      &- \\
 
     1/30             &1.81E-14 &16.73  &2.83E-14     &16.74     &1.55E-12 & 16.66 &2.42E-12      &17.03 \\
 
     1/60             &1.45E-19 &16.93  &2.28E-19      &16.92     &1.28E-17 &16.89  &1.99E-17      &16.89 \\
 
     1/120             &1.12E-24 &16.98  &1.76E-24      &16.98     &9.97E-23 &16.97  &1.56E-22    &16.96\\
 
     1/240             &8.58E-30 &17.00  &1.35E-29     &16.99     &7.63E-28 &17.00  &1.20E-27      &16.99 \\
     \hline
 
     & \multicolumn{4}{|c|}{ENO-MR17,$\lambda=10^6,\alpha=3$} & \multicolumn{4}{|c|}{ENO-MR17,$\lambda=10^6,\alpha=4$}  \\
     \hline
      $h$       & $L_1$  &Order & $L_\infty$  &Order    & $L_1$  &Order & $L_\infty$  &Order    \\
     \hline
     1/15              &1.97E-03 &-     &3.10E-03      &-        &1.60E-01 &-     &3.24E-01      &- \\
 
     1/30             &1.81E-08 &16.73  &2.83E-08     &16.74     &1.55E-06 & 16.66 &2.42E-06      &17.03 \\
 
     1/60             &1.45E-13 &16.93  &2.28E-13      &16.92     &1.28E-11 &16.89  &1.99E-11      &16.89 \\
 
     1/120             &1.12E-18 &16.98  &1.76E-18      &16.98     &9.97E-17 &16.97  &1.56E-16    &16.96\\
 
     1/240             &8.58E-24 &17.00  &1.35E-23     &16.99     &7.63E-22 &17.00  &1.20E-21      &16.99 \\
     \hline
  \end{tabular}
 \end{table}

\begin{figure}[htbp]
 \centering
 \subfigure[WENO-AO(5,3),$\lambda=1$]{
 \label{FIG:Linear_Advection_WENO-AO-5-3-A1}
 \includegraphics[width=5.5 cm]{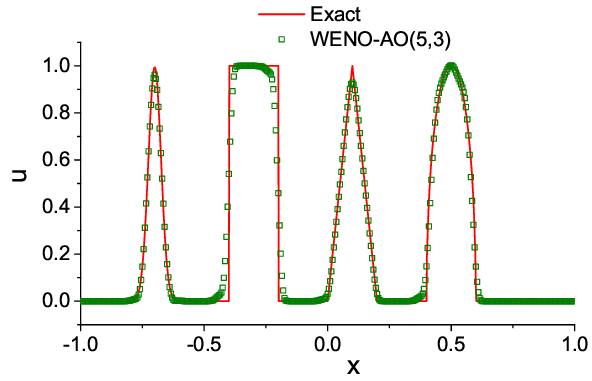}}
 \subfigure[WENO-AO(5,3),$\lambda=10^{-4}$]{
 \label{FIG:Linear_Advection_WENO-AO-5-3-A1E-4}
 \includegraphics[width=5.5 cm]{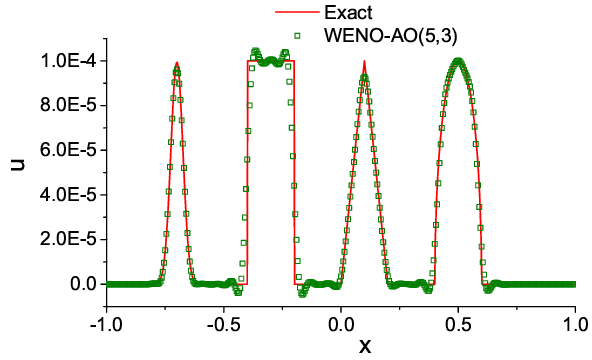}}
 \subfigure[ENO-MR5,$\lambda=1$]{
 \label{FIG:Linear_Advection_ENO-MR5-A1}
 \includegraphics[width=5.5 cm]{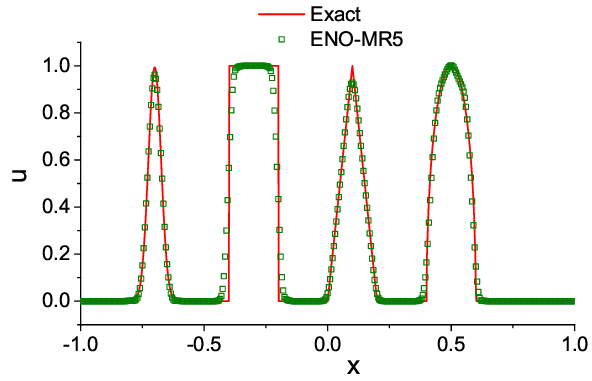}}
 \subfigure[ENO-MR5,$\lambda=10^{-4}$]{
 \label{FIG:Linear_Advection_ENO-MR5-A1E-4}
 \includegraphics[width=5.5 cm]{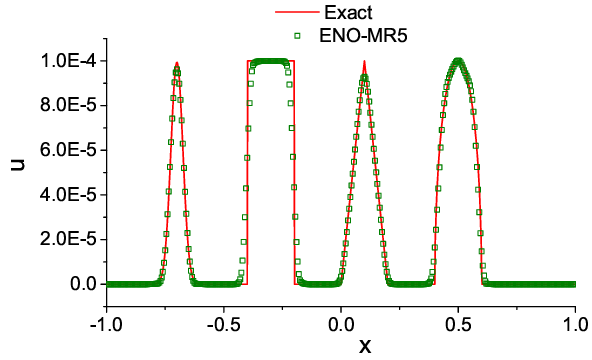}}
 \subfigure[ENO-MR13,$\lambda=1$]{
 \label{FIG:Linear_Advection_ENO-MR13-A1}
 \includegraphics[width=5.5 cm]{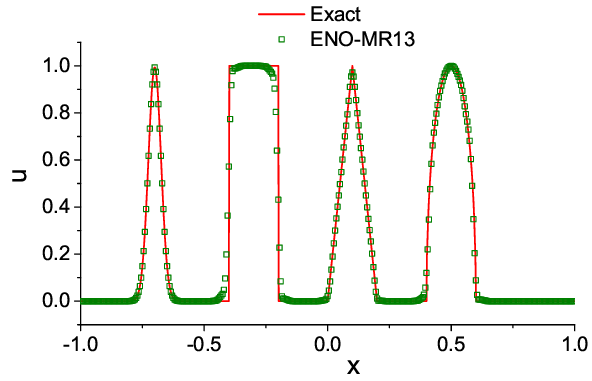}}
 \subfigure[ENO-MR17,$\lambda=10^{-4}$]{
 \label{FIG:Linear_Advection_ENO-MR17-A1E-4}
 \includegraphics[width=5.5 cm]{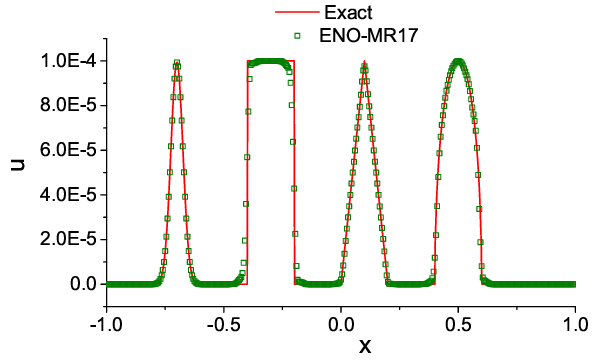}}
 \caption{The advection of a combination of Gaussians, a square wave, a sharp triangle wave,
 and a half ellipse (scaled by a factor of $\lambda$) arranged from left to right at $t=20$ 
 calculated by WENO-AO and ENO-MR schemes with $h=1/200$.}
\label{FIG:Linear_Advection}
\end{figure}

Next, we test the high-fidelity property of the proposed ENO-MR schemes for 
different shapes of solutions by advecting a combination of Gaussians, a square wave, a sharp
triangle wave, and a half ellipse arranged from left to right 
which was originally proposed by Jiang and Shu \cite{Jiang_Shu1996WENO}.
To test the scale-invariant property, the solution is scaled by a factor of $\lambda$.
Fig. \ref{FIG:Linear_Advection} shows the profile of $u$ 
at $t=20$ calculated by WENO-AO and ENO-MR schemes with $h=1/200$.
We observe that WENO-AO(5,3) scheme is highly oscillatory
when the solution is scaled by a small factor $\lambda=10^{-4}$.
WENO-AO(9,5,3) scheme obtains similar results
which are omitted to save space.
This is because the nonlinear weights approach to linear weights 
when the scale of the solution is small such that $\beta_k^{(r)}$ and $\tau_{2r-1}$ 
become much smaller than $\varepsilon$ in Eq. (\ref{Eq:WENO_Z_Weight}).
However, ENO-MR schemes are accurate for all shapes of solutions 
regardless of the scale of the solution.
This demonstrates the scale-invariant property of ENO-MR schemes again.
Furthermore, ENO-MR13 and ENO-MR17 schemes are more accurate
than WENO-AO(5,3) and ENO-MR5 schemes, 
particularly for the triangle wave.

%%%%%%%%%%%%%%%%%%%%%%%%%%%%%%%%%%%%%%%%%%%%%%%%%%%%%%%%%%%%%%%%%%%%%%%%%%%%%%%%%%%%%%%%%%%
\subsection{Numerical examples for the 1D Burgers equation}
\begin{figure}[htbp]
  \centering
  \subfigure[5th-order,$\lambda=1$]{
  \label{FIG:Burgers-5th-order-A1}
  \includegraphics[width=5.5 cm]{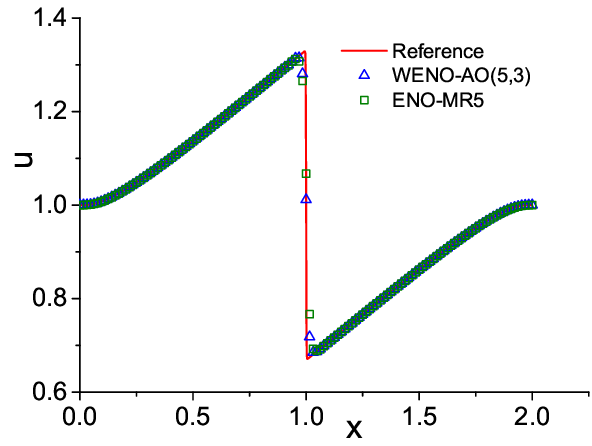}}
  \subfigure[5th-order,$\lambda=10^{-3}$]{
  \label{FIG:Burgers-5th-order-A1E-3}
  \includegraphics[width=5.5 cm]{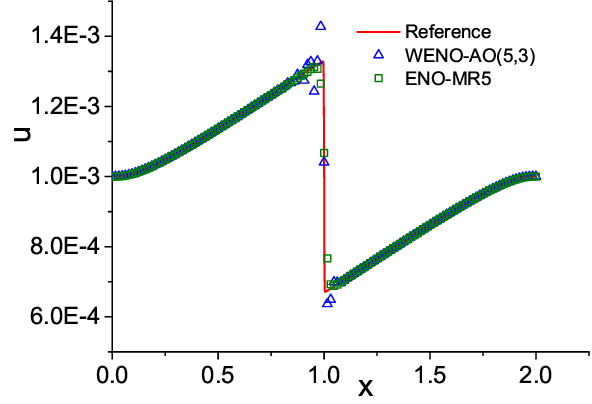}}
  \subfigure[9th-order,$\lambda=1$]{
  \label{FIG:Burgers-9th-order-A1}
  \includegraphics[width=5.5 cm]{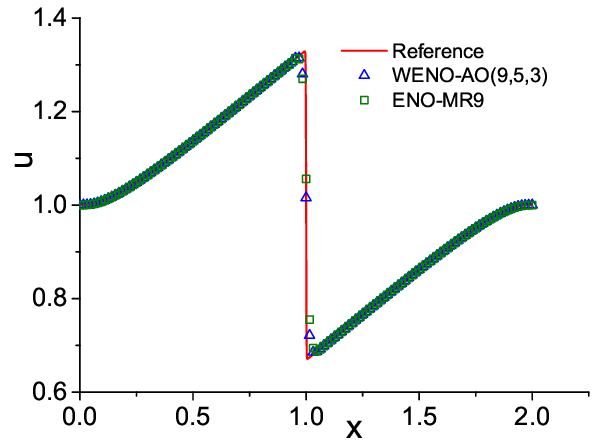}}
  \subfigure[9th-order,$\lambda=10^{-3}$]{
  \label{FIG:Burgers-9th-order-A1E-3}
  \includegraphics[width=5.5 cm]{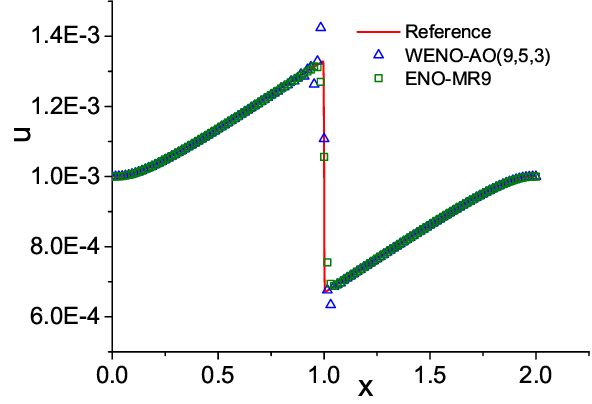}}
  \subfigure[13th- and 17th- order,$\lambda=1$]{
  \label{FIG:Burgers-13th-17th-order-A1}
  \includegraphics[width=5.5 cm]{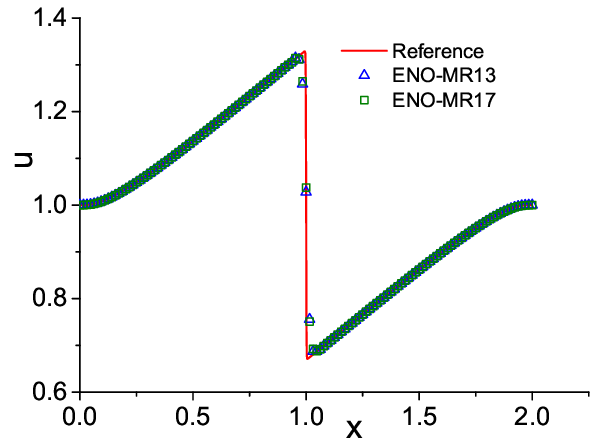}}
  \subfigure[13th- and 17th- order,$\lambda=10^{-3}$]{
  \label{FIG:Burgers-13th-17th-order-A1E-3}
  \includegraphics[width=5.5 cm]{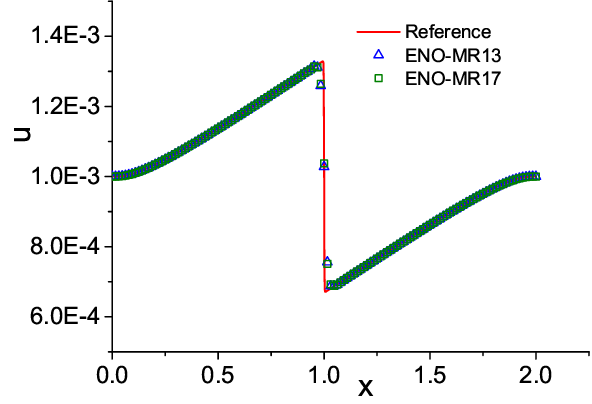}}
  \caption{Solutions of the Burgers equation at $t=2/\lambda$ 
  calculated by WENO-AO and ENO-MR schemes with $h=1/64$.}
 \label{FIG:Burgers}
 \end{figure}

We consider the 1D inviscid Burgers equation
\begin{equation}\label{Eq:1DAdvEq}
 \frac{\partial u}{\partial t}+\frac{\partial u^2/2}{\partial x}=0.
\end{equation}
The computational domain is $[0,2]$,
the initial condition is $u_0(x)=\lambda(1+0.5sin^3(\pi x))$,
and periodic boundary conditions are applied to both sides.
It is easy to check that $\frac{\partial f}{\partial u}=u$ is positive in the whole domain,
so flux-splitting is not required.
First, we run the problem to $t=0.1/\lambda$ when the solution is still smooth.
The time step is calculated by $\Delta t=h^{5/3}/\lambda$, $\Delta t=100\times h^{3}/\lambda$,
$\Delta t=10^4\times h^{13/3}/\lambda$, and $\Delta t=10^6\times h^{17/3}/\lambda$ 
for 5th-, 9th-, 13th-, and 17th-order schemes, which guarantees the same order in space and time 
and the same CFL number for different values of the scale factor $\lambda$.
Table \ref{Table:1DBurgersConvergence} shows the numerical errors of $u$ computed by
WENO-AO and ENO-MR schemes with different mesh sizes.
We observe that WENO-AO(5,3) and WENO-AO(9,5,3) lose accuracy on coarse meshes
and the corresponding errors depend on the scale factor $\lambda$, 
which is similar to the results for linear equations.
However, ENO-MR schemes almost achieve the optimal convergence order 
and are scale-invariant.
Next, we run this problem up to $t=2/\lambda$ when a strong discontinuity appears at the center.
Fig. \ref{FIG:Burgers} shows the profiles of $u$ calculated by WENO-AO and ENO-MR schemes with $h=1/64$ and $CFL=0.3$.
All schemes can sharply capture the discontinuity without spurious oscillations when $\lambda=1$,
but WENO-AO schemes generate numerical oscillations near the discontinuity when $\lambda=0.001$.
However, ENO-MR schemes retain non-oscillation regardless of the values of the scale factor.

\begin{table}[htbp]
  \footnotesize
  \caption{Numerical errors of the Burgers equation at $t=0.1/\lambda$ computed by WENO-AO and ENO-MR schemes.}
  \label{Table:1DBurgersConvergence}
  \centering
  \begin{tabular}{|c|c|c|c|c|c|c|c|c|}

    \hline
    % after \\: \hline or \cline{col1-col2} \cline{col3-col4} ...
    & \multicolumn{4}{|c|}{WENO-AO(5,3),$\lambda=1$} & \multicolumn{4}{|c|}{WENO-AO(5,3),$\lambda=10^3$} \\
     \hline
      $h$       & $L_1$  &Order & $L_\infty$  &Order    & $L_1$  &Order & $L_\infty$  &Order    \\
     \hline
     1/32              &5.46E-05 &-     &3.94E-04     &-        &5.46E-02 &-     &3.94E-01      &- \\
 
     1/64             &5.06E-06 &3.43  &5.89E-05     &2.74     &5.06E-03 &3.43  &5.89E-02      &2.74 \\
 
     1/128             &4.25E-07 &3.57  &1.17E-05      &2.33     &4.26E-04 &3.57  &1.17E-02      &2.33 \\
 
     1/256             &3.85E-08 &3.46  &1.92E-06      &2.61     &4.17E-05 &3.35  &2.02E-03    &2.53 \\

     1/512             &1.09E-10 &8.47  &1.06E-08      &7.50     &4.78E-06 &3.13  &4.96E-04     &2.03 \\
     \hline

      & \multicolumn{4}{|c|}{ENO-MR5,$\lambda=1$} & \multicolumn{4}{|c|}{ENO-MR5,$\lambda=10^3$} \\
     \hline
      $h$       & $L_1$  &Order & $L_\infty$  &Order    & $L_1$  &Order & $L_\infty$  &Order    \\
     \hline
     1/32              &1.73E-05 &-     &1.43E-04     &-        &1.73E-02 &-     &1.43E-01      &- \\
 
     1/64             &3.09E-07 &5.81  &1.73E-06     &6.37     &3.09E-04 &5.81  &1.73E-03      &6.37 \\
 
     1/128             &9.79E-09 &4.98  &5.51E-08      &4.97     &9.79E-06 &4.98  &5.51E-05      &4.97 \\
 
     1/256             &3.06E-10 &5.00  &1.73E-09      &4.99     &3.06E-07 &5.00  &1.73E-06    &4.99 \\

     1/512             &9.55E-12 &5.00  &5.39E-11      &5.00     &9.55E-09 &5.00  &5.39E-08     &5.00 \\
     \hline

     & \multicolumn{4}{|c|}{WENO-AO(9,5,3),$\lambda=1$} & \multicolumn{4}{|c|}{WENO-AO(9,5,3),$\lambda=10^3$} \\
     \hline
      $h$       & $L_1$  &Order & $L_\infty$  &Order    & $L_1$  &Order & $L_\infty$  &Order    \\
     \hline
     1/32              &8.79E-05 &-     &5.62E-04     &-        &8.79E-02 &-     &5.62E-01      &- \\
 
     1/64             &8.91E-06 &3.30  &8.30E-05     &2.76     &8.91E-03 &3.30  &8.30E-02      &2.76 \\
 
     1/128             &6.82E-07 &3.71  &1.58E-05      &2.39     &6.83E-04 &3.71  &1.58E-02      &2.39 \\
 
     1/256             &5.20E-08 &3.71  &2.55E-06      &2.63     &5.49E-05 &3.64  &2.64E-03    &2.58 \\

     1/512             &8.21E-10 &5.99  &9.00E-08      &4.82     &4.97E-06 &3.47  &4.46E-04     &2.57\\
     \hline

     & \multicolumn{4}{|c|}{ENO-MR9,$\lambda=1$} & \multicolumn{4}{|c|}{ENO-MR9,$\lambda=10^3$} \\
     \hline
      $h$       & $L_1$  &Order & $L_\infty$  &Order    & $L_1$  &Order & $L_\infty$  &Order    \\
     \hline
     1/32              &3.28E-07 &-     &2.68E-06     &-        &3.28E-04 &-     &2.68E-03      &- \\
 
     1/64             &7.86E-10 &8.71  &7.35E-09     &8.51     &7.86E-07 &8.71  &7.35E-06      &8.51 \\
 
     1/128             &1.65E-12 &8.90  &1.67E-11      &8.78     &1.65E-09 &8.90  &1.67E-08      &8.78 \\
 
     1/256             &3.28E-15 &8.98  &3.40E-14      &8.94     &3.28E-12 &8.98  &3.40E-11    &8.94 \\

     1/512             &6.45E-18 &8.99  &6.71E-17      &8.99     &6.45E-15 &8.99  &6.71E-14     &8.99 \\
     \hline

     & \multicolumn{4}{|c|}{ENO-MR13,$\lambda=1$} & \multicolumn{4}{|c|}{ENO-MR13,$\lambda=10^3$} \\
     \hline
      $h$       & $L_1$  &Order & $L_\infty$  &Order    & $L_1$  &Order & $L_\infty$  &Order    \\
     \hline
     1/32              &1.85E-07 &-     &1.27E-06     &-        &1.85E-04 &-     &1.27E-03      &- \\
 
     1/64             &2.59E-11 &12.80  &2.42E-10     &12.36     &2.59E-08 &12.80  &2.42E-07      &12.36 \\
 
     1/128             &3.64E-15 &12.80  &3.88E-14      &12.61     &3.64E-12 &12.80  &3.88E-11      &12.61 \\
 
     1/192             &1.94E-17 &12.91  &2.10E-16      &12.87     &1.94E-14 &12.91  &2.10E-13    &12.87 \\

     1/256             &4.65E-19 &12.97  &5.09E-18      &12.93     &4.65E-16 &12.97  &5.09E-15     &12.93 \\
     \hline

     & \multicolumn{4}{|c|}{ENO-MR17,$\lambda=1$} & \multicolumn{4}{|c|}{ENO-MR17,$\lambda=10^3$} \\
     \hline
      $h$       & $L_1$  &Order & $L_\infty$  &Order    & $L_1$  &Order & $L_\infty$  &Order    \\
     \hline
     1/32              &1.72E-07 &-     &1.17E-06     &-        &1.72E-04 &-     &1.17E-03      &- \\
 
     1/64             &1.61E-12 &16.71  &1.56E-11     &16.20     &1.61E-09 &16.71  &1.56E-08      &16.20 \\
 
     1/96             &2.16E-15 &16.31  &2.33E-14      &16.05     &2.16E-12 &16.31  &2.33E-11      &16.05 \\
 
     1/128             &1.83E-17 &16.58  &2.20E-16      &16.21     &1.83E-14 &16.58  &2.20E-13    &16.21 \\

     1/160             &4.39E-19 &16.72  &5.57E-18      &16.47     &4.39E-16 &16.72  &5.57E-15     &16.47 \\
     \hline

  \end{tabular}
 \end{table}

%%%%%%%%%%%%%%%%%%%%%%%%%%%%%%%%%%%%%%%%%%%%%%%%%%%%%%%%%%%%%%%%%%%%%%%%%%%%%%%%%%%%%%%%%%
\subsection{Numerical examples for the 1D Euler equations}
We consider the 1D compressible Euler equations
\begin{equation}
 \frac{\partial \mathbf{U}}{\partial t}+\frac{\partial \mathbf{F}}{\partial x}=0,
\end{equation}
where $\mathbf{U}=[\rho,\rho u,\rho e]^T$ and $\mathbf{F}=[{\rho u,\rho u^2+p,(\rho e+p)u}]^T$
with $\rho$, $u$, $p$, and $e$ representing the density, velocity, pressure, and specific total energy respectively.
The specific total energy is calculated as $e=\frac{p}{\rho(\gamma-1)}+\frac{1}{2} u^2$,
where $\gamma$ is the specific heat ratio.
Without particular instructions, $\gamma$ is set to 1.4 hereafter.

\begin{figure}[htbp]
  \centering
  \subfigure[5th-order]{
  \label{FIG:Lax-5th-order}
  \includegraphics[width=3.8 cm]{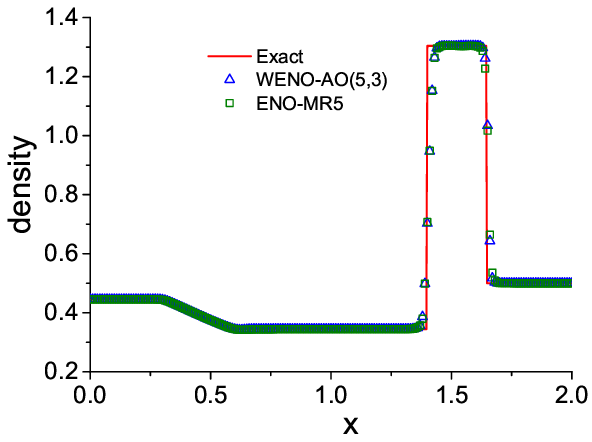}}
  \subfigure[9th-order]{
  \label{FIG:Lax-9th-order}
  \includegraphics[width=3.8 cm]{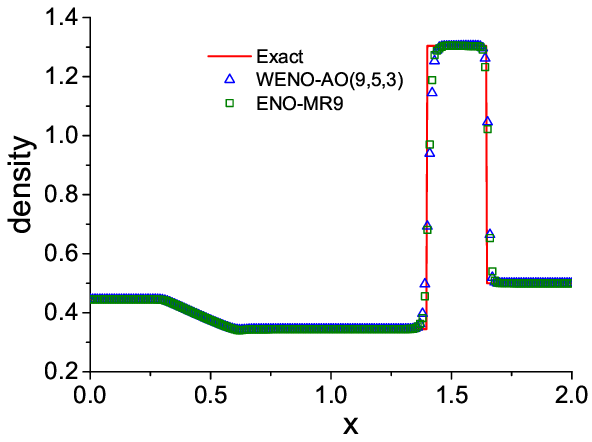}}
  \subfigure[13th- and 17th-order]{
  \label{FIG:Lax-13th-17th-order}
  \includegraphics[width=3.8 cm]{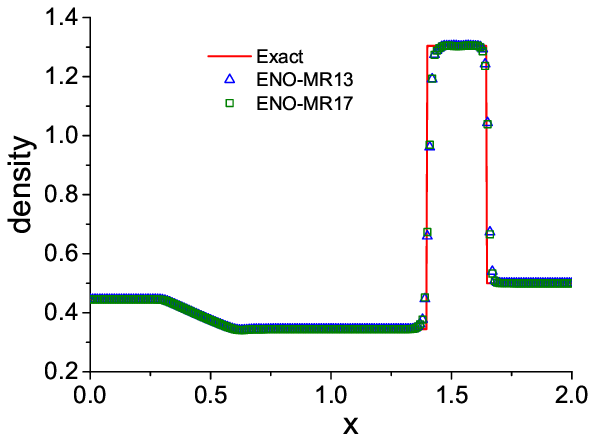}}
  \caption{Density profiles of Lax shock tube problem at $t=0.26$ 
  calculated by WENO-AO and ENO-MR schemes with $h=1/100$.}
 \label{FIG:Lax}
 \end{figure}

The first case is the Lax shock tube problem that is used to test 
the shock-capturing capability of numerical methods.
The computational domain is [0,2], and the initial condition is given by
\begin{equation*}
 (\rho,u,p)=\begin{cases}
              (0.445,0.698,3.528), & \mbox{if } x<1 \\
              (0.5,0,0.571), & \mbox{otherwise}.
            \end{cases}
\end{equation*}
Non-reflection boundary conditions are applied to the left and right boundaries.
Fig. \ref{FIG:Lax} shows the density and velocity profiles at $t=0.26$ 
calculated by WENO-AO and ENO-MR schemes with $h=1/100$.
We observe that both WENO-AO and ENO-MR schemes can capture discontinuity very well.
Even when the designed order increases to a very high level,
ENO-MR can still keep ENO properties.

\begin{figure}[htbp]
  \centering
  \subfigure[5th-order]{
  \label{FIG:Titarev-5th-order}
  \includegraphics[width=10 cm]{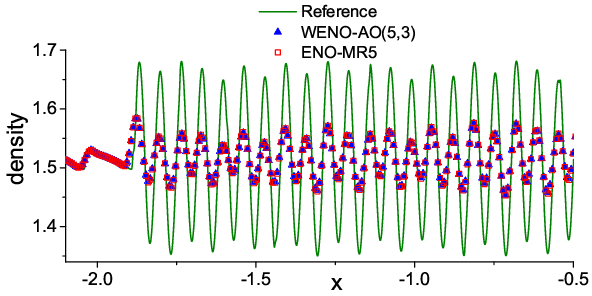}}
  \subfigure[9th-order]{
  \label{FIG:Titarev-9th-order}
  \includegraphics[width=10 cm]{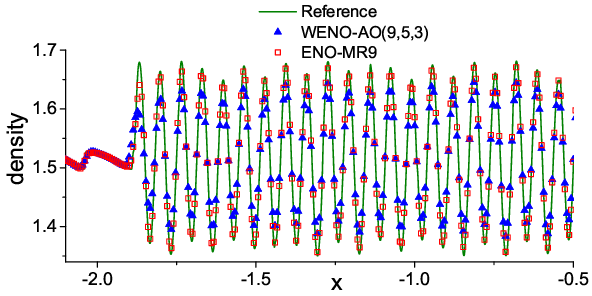}}
  \subfigure[13th- and 17th-order]{
  \label{FIG:Titarev-13th-17th-order}
  \includegraphics[width=10 cm]{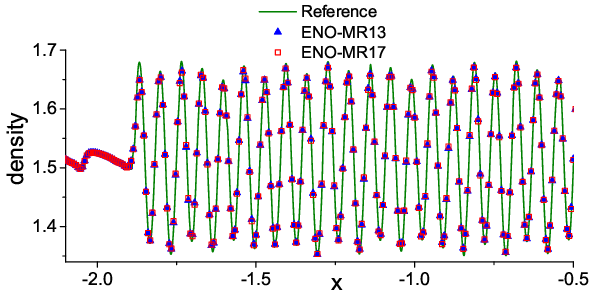}}
  \caption{Density profiles of Titarev-Toro problem at $t=5$ 
  calculated by WENO-AO and ENO-MR schemes with $h=1/150$.}
 \label{FIG:Titarev_Toro}
 \end{figure}

The second case is the Titarev and Toro \cite{Titarev2004WENO_FV} problem 
which is an upgraded version of the Shu-Osher problem \cite{Shu1988EfficientENO}.
It depicts the interaction of a Mach 1.1 moving shock with a very high-frequency entropy sine wave.
The computational domain is [-5,5], and the initial condition is given by
\begin{equation*}
 (\rho,u,p)=\begin{cases}
              (1.515695,0.523346,1.805), & \mbox{if } x<-4.5 \\
              (1+0.1\mbox{sin}(20\pi x),0,1), & \mbox{otherwise}.
            \end{cases}
\end{equation*}
Non-reflection boundary conditions are applied to the left and right sides.
Fig. \ref{FIG:Titarev_Toro} shows the density profiles in the high-frequency region at $t=5$ 
calculated by WENO-AO and ENO-MR schemes with $h=1/150$.
The reference solution is calculated by the matured WENO-Z5 scheme with $h=1/1000$.
We observe that WENO-AO(5,3) and ENO-MR5 perform similarly.
ENO-MR9 performs much better than WENO-AO(9,5,3).
ENO-MR13 and ENO-MR17 have no obvious advantage over ENO-MR9
at the current resolution.

%%%%%%%%%%%%%%%%%%%%%%%%%%%%%%%%%%%%%%%%%%%%%%%%%%%%%%%%%%%%%%%%%%%%%%%%%%%%%%%%%%%%%%%%%%%%%%%%%
\subsection{Numerical examples for the 2D Euler equations}
The 2D compressible Euler equations can be written as 
\begin{equation}\label{Eq:2D_Euler_equations}
 \frac{\partial \mathbf{U}}{\partial t}+\frac{\partial \mathbf{F}}{\partial x}+\frac{\partial \mathbf{G}}{\partial y}=0,
\end{equation}
with $\mathbf{U}=[\rho,\rho u,\rho v,\rho e]^T$, $\mathbf{F}=[{\rho u,\rho u^2+p,\rho uv,(\rho e+p)u}]^T$,
and $\mathbf{G}=[{\rho v,\rho uv,\rho v^2+p,(\rho e+p)v}]^T$,
where $\rho$, $u$, $v$, $p$, and $e$ denote the density, $x$-velocity, $y$-velocity, pressure, and specific total energy respectively.
The specific total energy is calculated as $e=\frac{p}{\rho(\gamma-1)}+\frac{1}{2} (u^2+v^2)$.

%%%%%%%%%%%%%%%%%%%%%%%%%%%%%%%%%%%%%%%%%%%%%%%%%%%%%%%%%%%%%%%%%%%%%%%%%%%%%%%%%%%%%%%%%%%%%%%%%%%%
\subsubsection{Two-dimensional Riemann problems}
\begin{figure}[htbp]
  \centering
  \subfigure[WENO-AO(5,3)]{
  \label{FIG:RP1-WENO-AO-5-3}
  \includegraphics[width=5.5 cm]{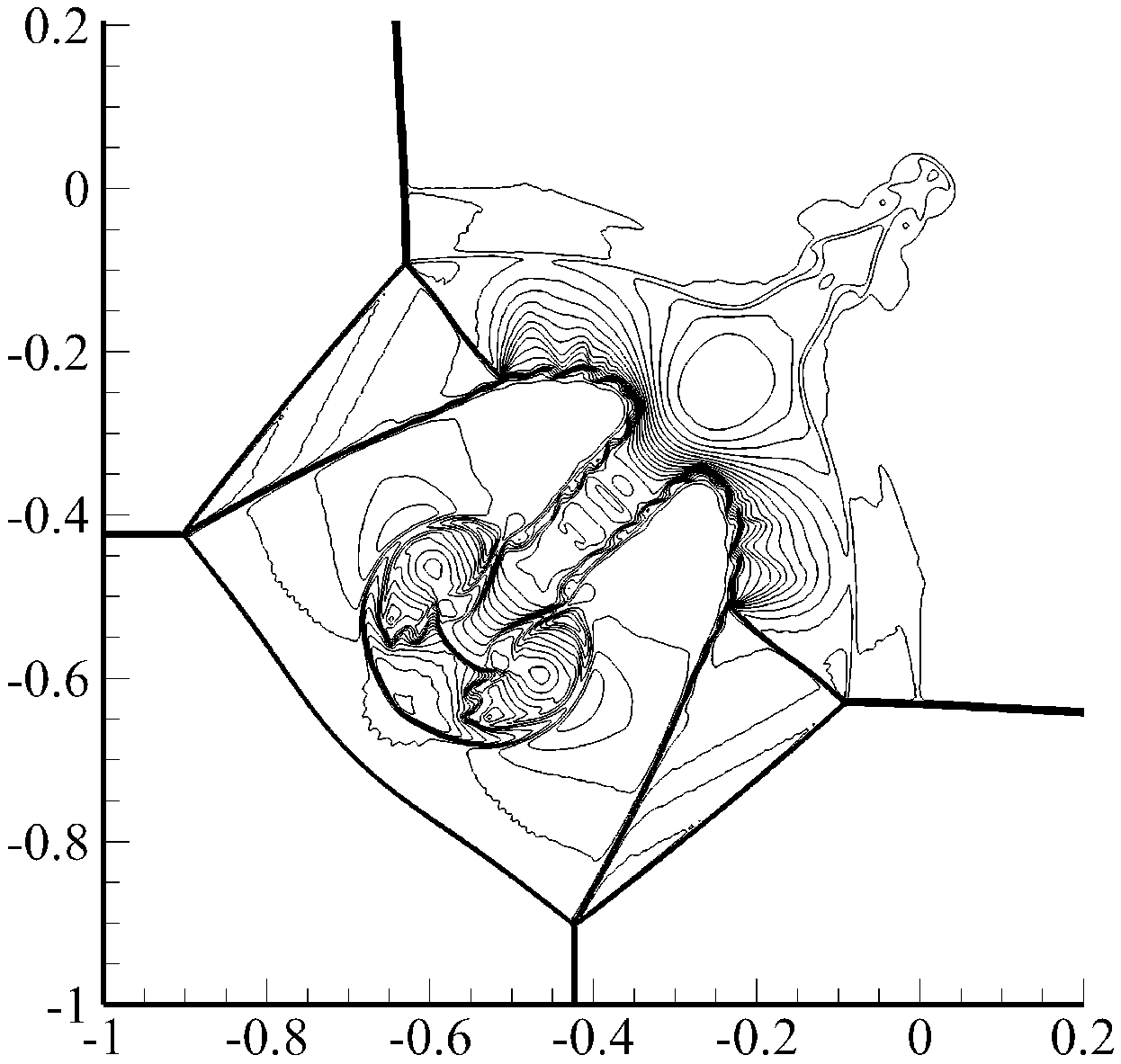}}
  \subfigure[WENO-AO(9,5,3)]{
  \label{FIG:RP1-WENO-AO-9-5-3}
  \includegraphics[width=5.5 cm]{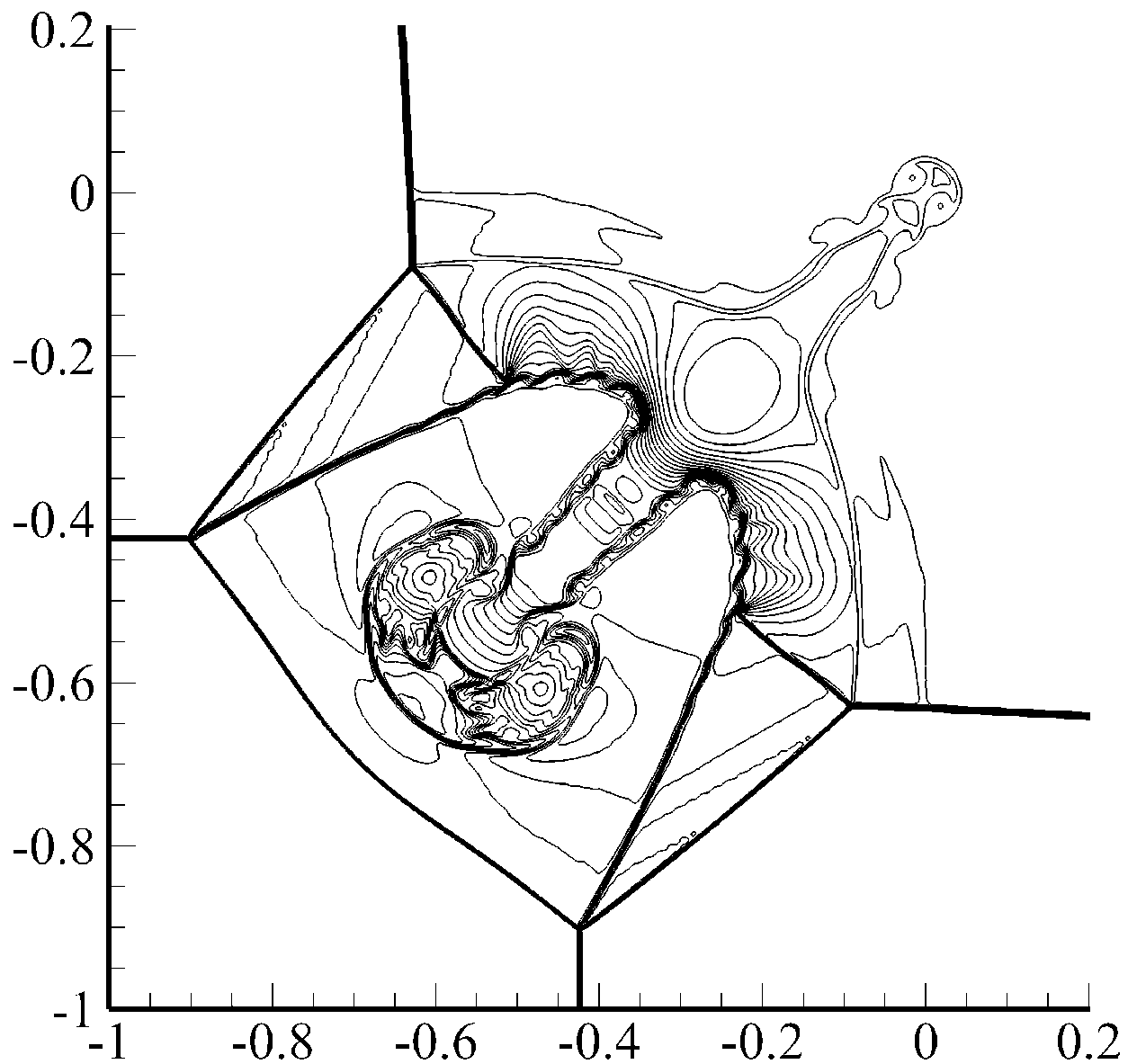}}
  \subfigure[ENO-MR5]{
  \label{FIG:RP1-ENO-MR5}
  \includegraphics[width=5.5 cm]{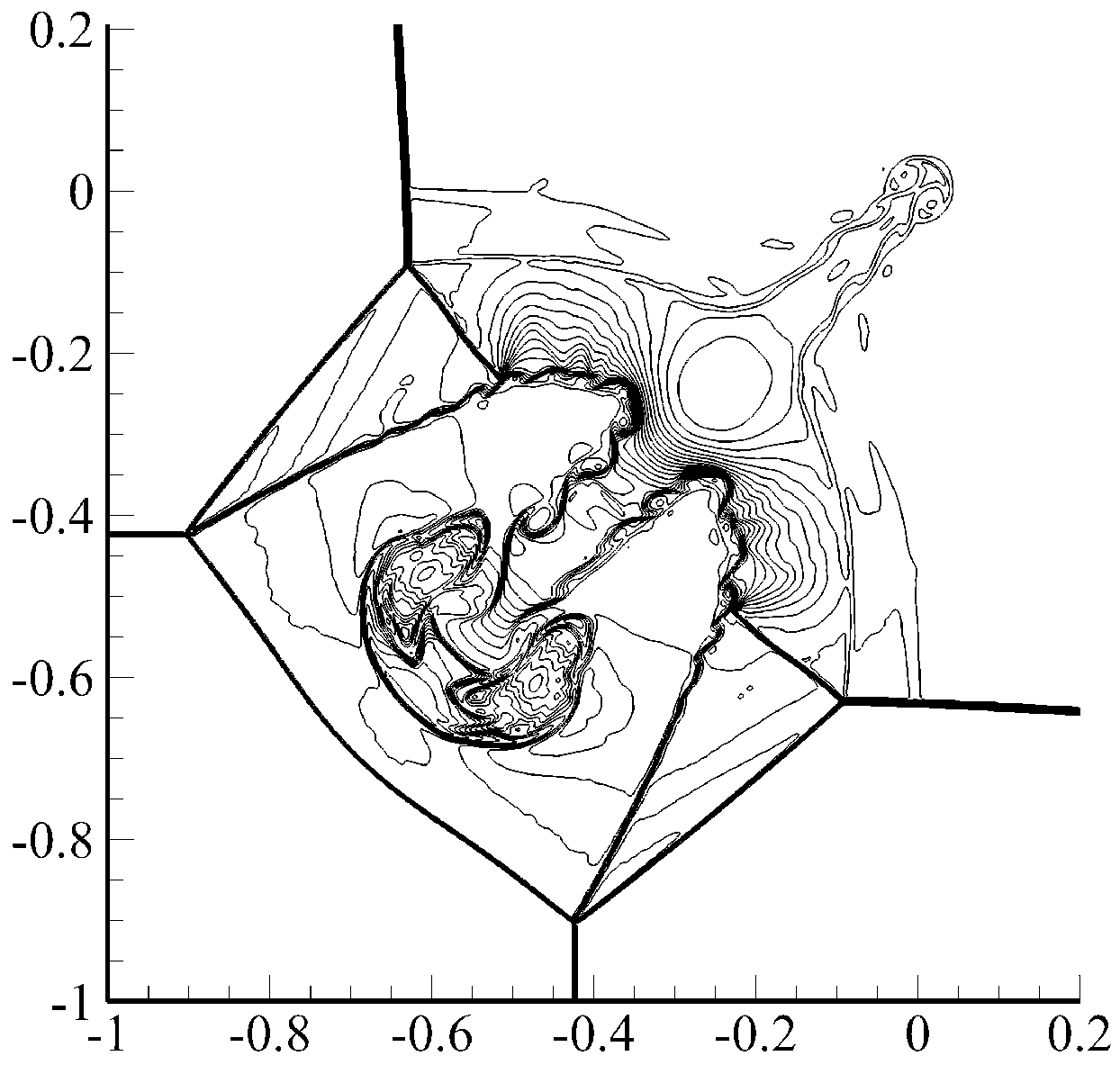}}
  \subfigure[ENO-MR9]{
  \label{FIG:RP1-ENO-MR9}
  \includegraphics[width=5.5 cm]{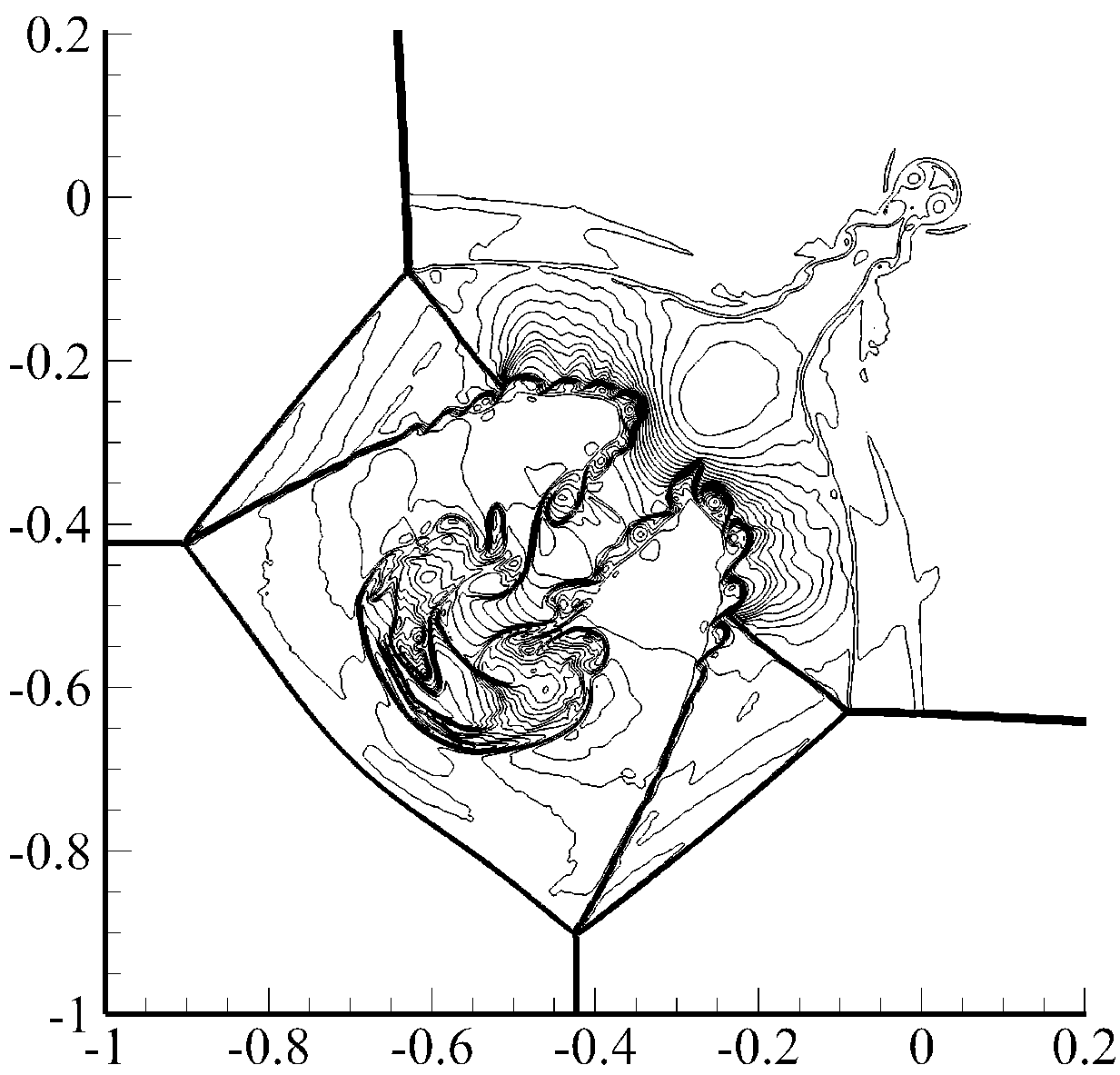}}
  \subfigure[ENO-MR13]{
  \label{FIG:RP1-ENO-MR13}
  \includegraphics[width=5.5 cm]{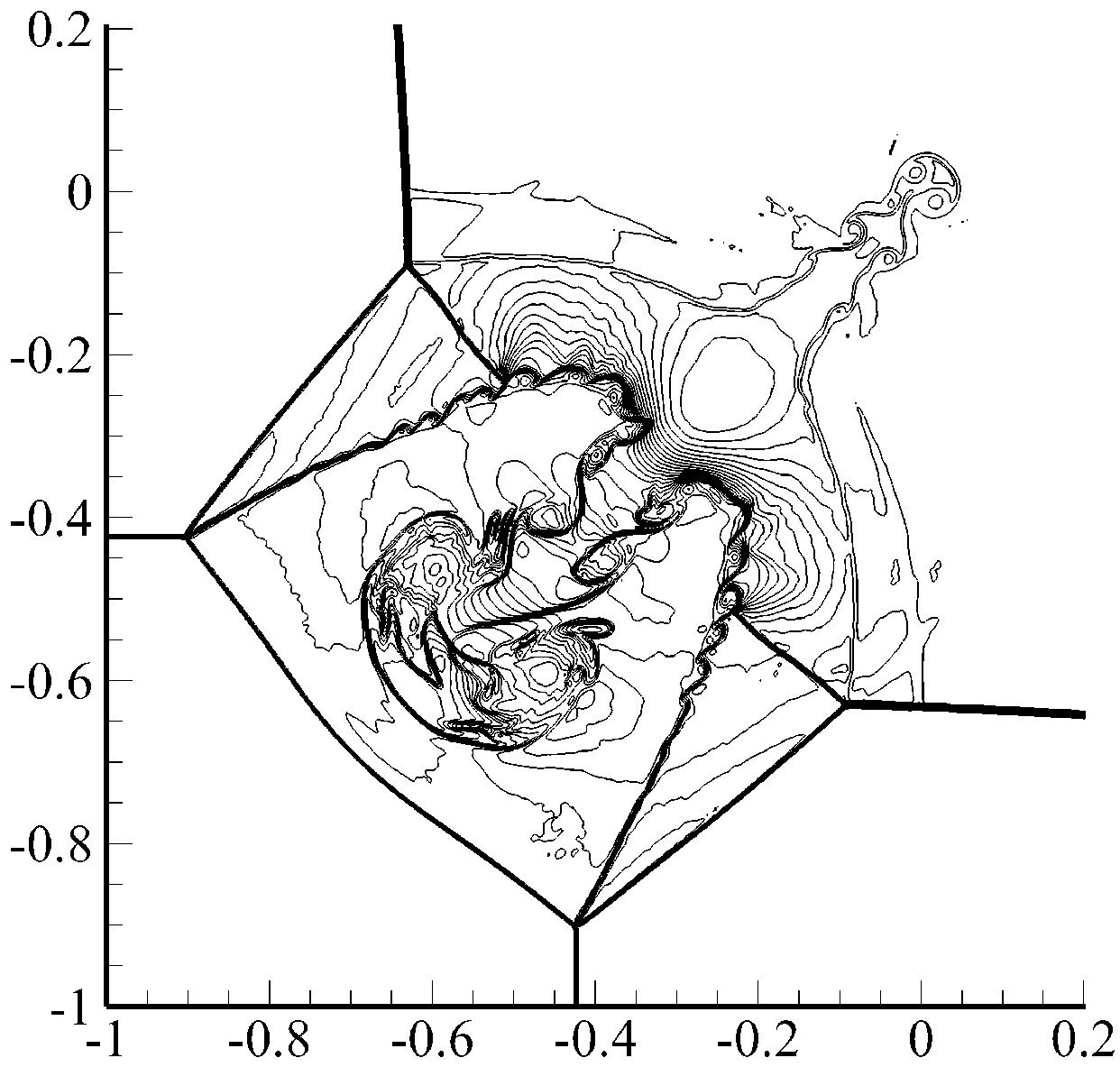}}
  \subfigure[ENO-MR17]{
  \label{FIG:RP1-ENO-MR17}
  \includegraphics[width=5.5 cm]{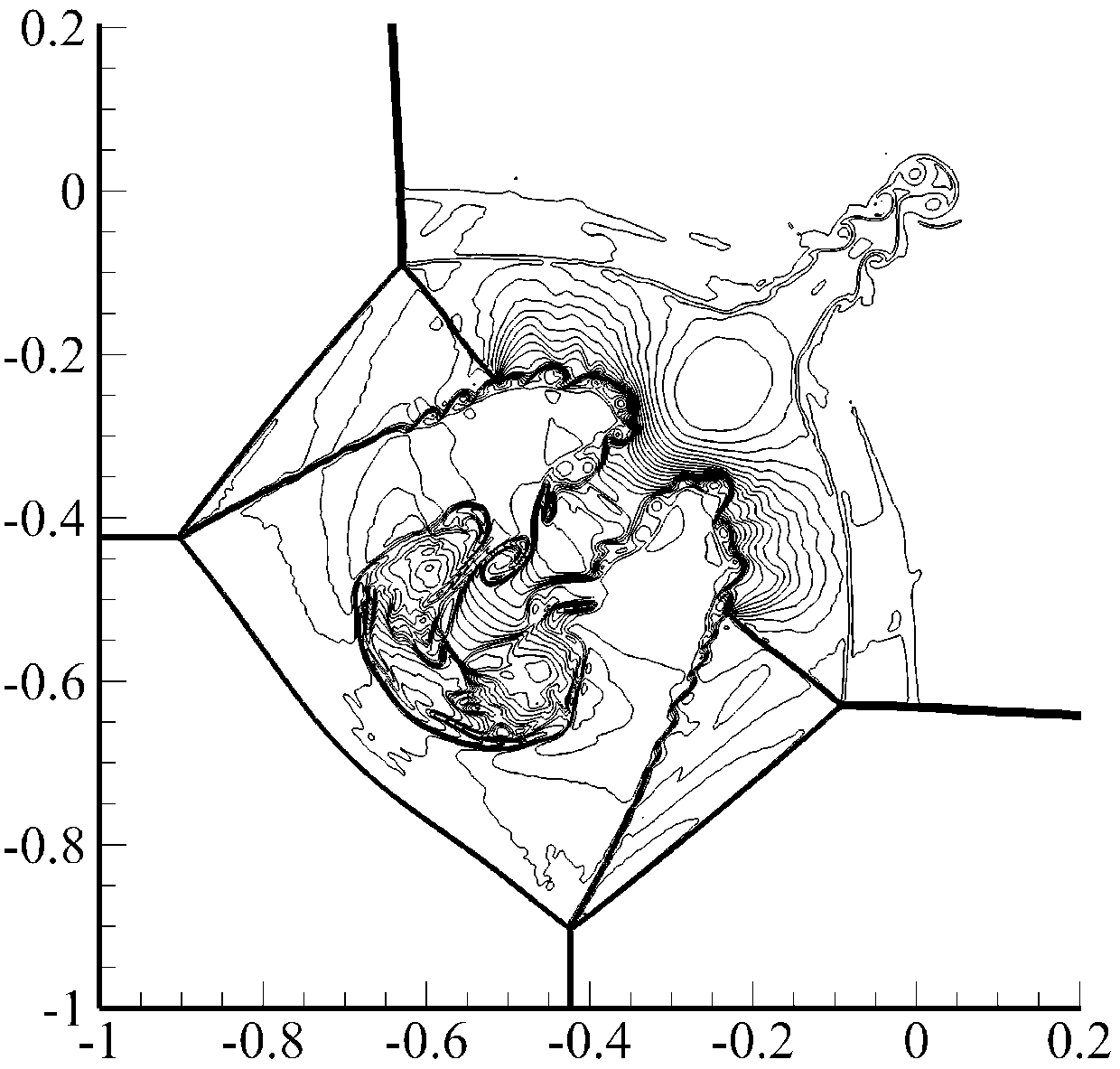}}
  
  \caption{Density contours of the first 2D Riemann problem at $t=1$ calculated 
  by WENO-AO and ENO-MR schemes with $801\times801$ mesh points.
  The contours contain 30 equidistant contours from 0.2 to 1.7.}
 \label{FIG:RP1}
 \end{figure}

The computational is initially separated into four parts
which are filled with gases of different states.
We can simulate abundant interactions of different waves in two-dimensional space
by providing different initial conditions.
In this study, we consider two typical configurations.
The first configuration is given by
\begin{equation*}
 \left(\rho, u, v, p\right)=
 \begin{cases}
   \left(0.138,1.206,1.206,0.029\right),    & \text{in }[-1, 0]\times[-1, 0],\\
   \left(0.5323,1.206,0,0.3\right),         & \text{in }[-1, 0]\times[0, 1],\\
   \left(1.5,0,0,1.5\right),                & \text{in }[0, 1] \times[0, 1],\\
   \left(0.5323,0,1.206,0.3\right),         & \text{in }[0, 1] \times[-1, 0],
 \end{cases}
\end{equation*}
which depicts interactions of two horizontally moving shocks and two vertically moving shocks.
The collision of the four normal shocks will 
result in two double-Mach reflections and an oblique shock moving along the diagonal of the computational domain.
Fig. \ref{FIG:RP1} shows density contours at $t=1$ calculated by WENO-AO and ENO-MR schemes with $801\times801$ mesh points.
At the current resolution, WENO-AO schemes start to trace out Kelvin-Helmholtz (KH) instabilities along the slip lines,
but ENO-MR schemes can clearly capture more details.

\begin{figure}[htbp]
  \centering
  \subfigure[WENO-AO(5,3)]{
  \label{FIG:RP2-WENO-AO-5-3}
  \includegraphics[width=5.5 cm]{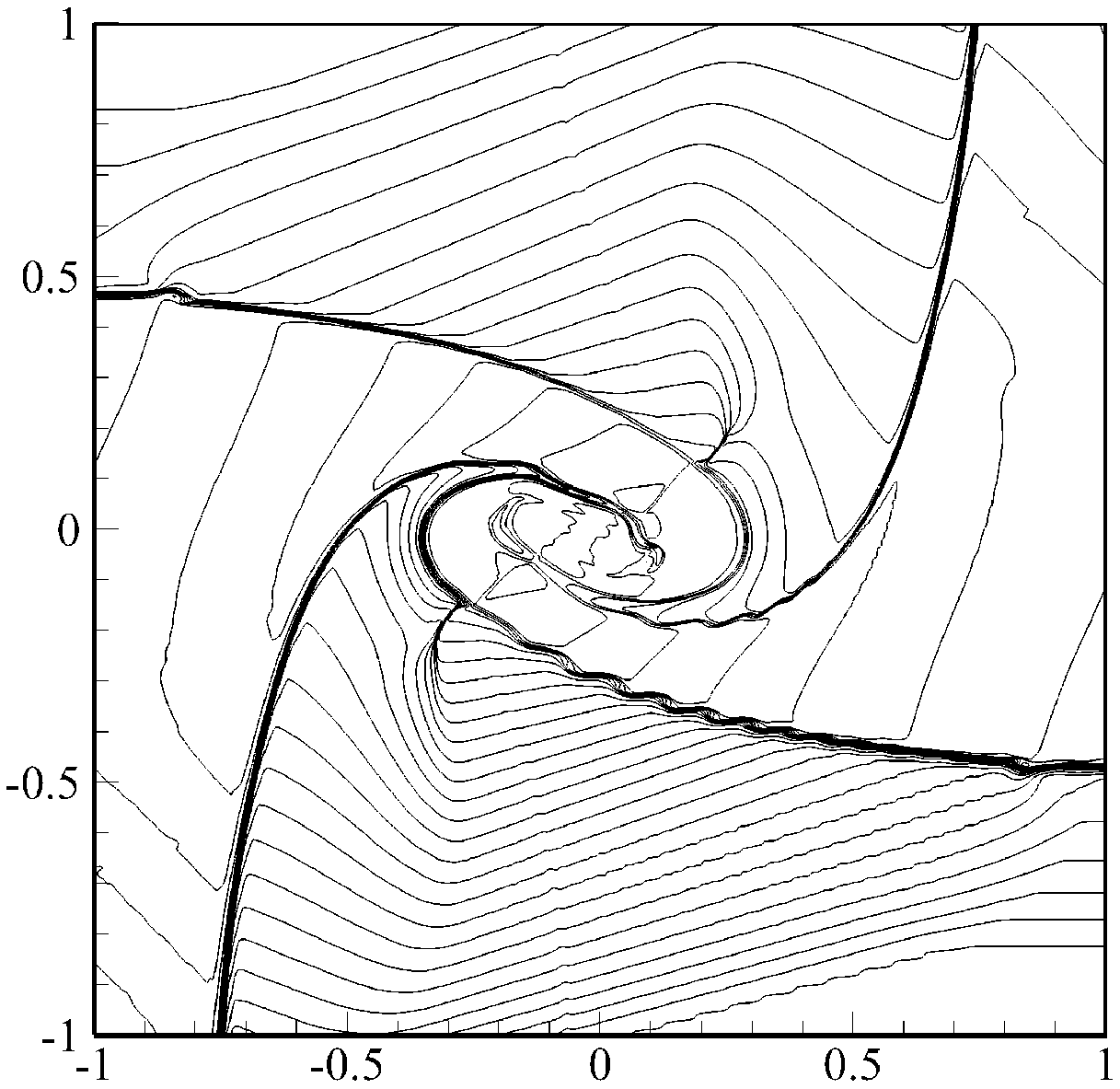}}
  \subfigure[WENO-AO(9,5,3)]{
  \label{FIG:RP2-WENO-AO-9-5-3}
  \includegraphics[width=5.5 cm]{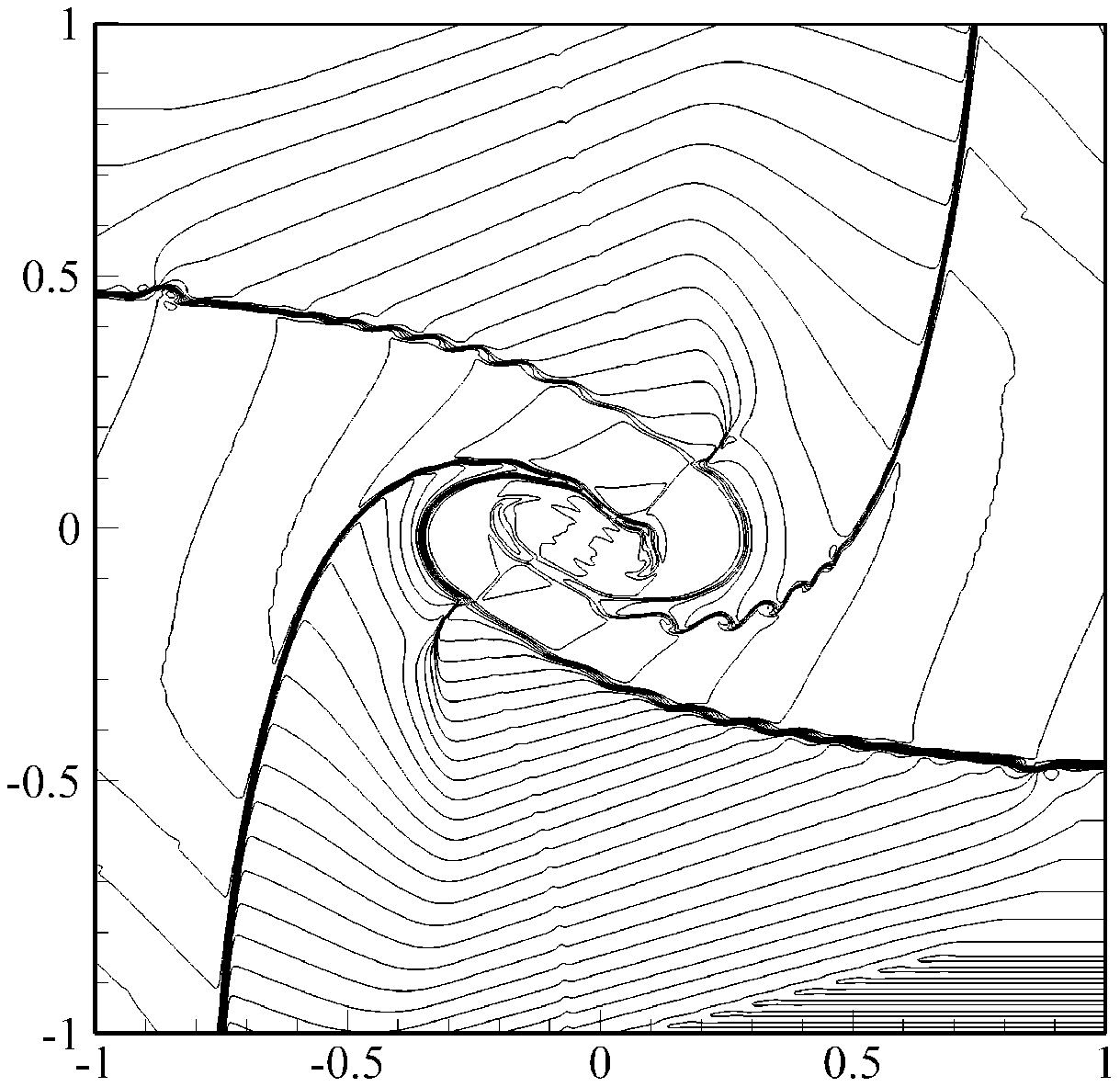}}
  \subfigure[ENO-MR5]{
  \label{FIG:RP2-ENO-MR5}
  \includegraphics[width=5.5 cm]{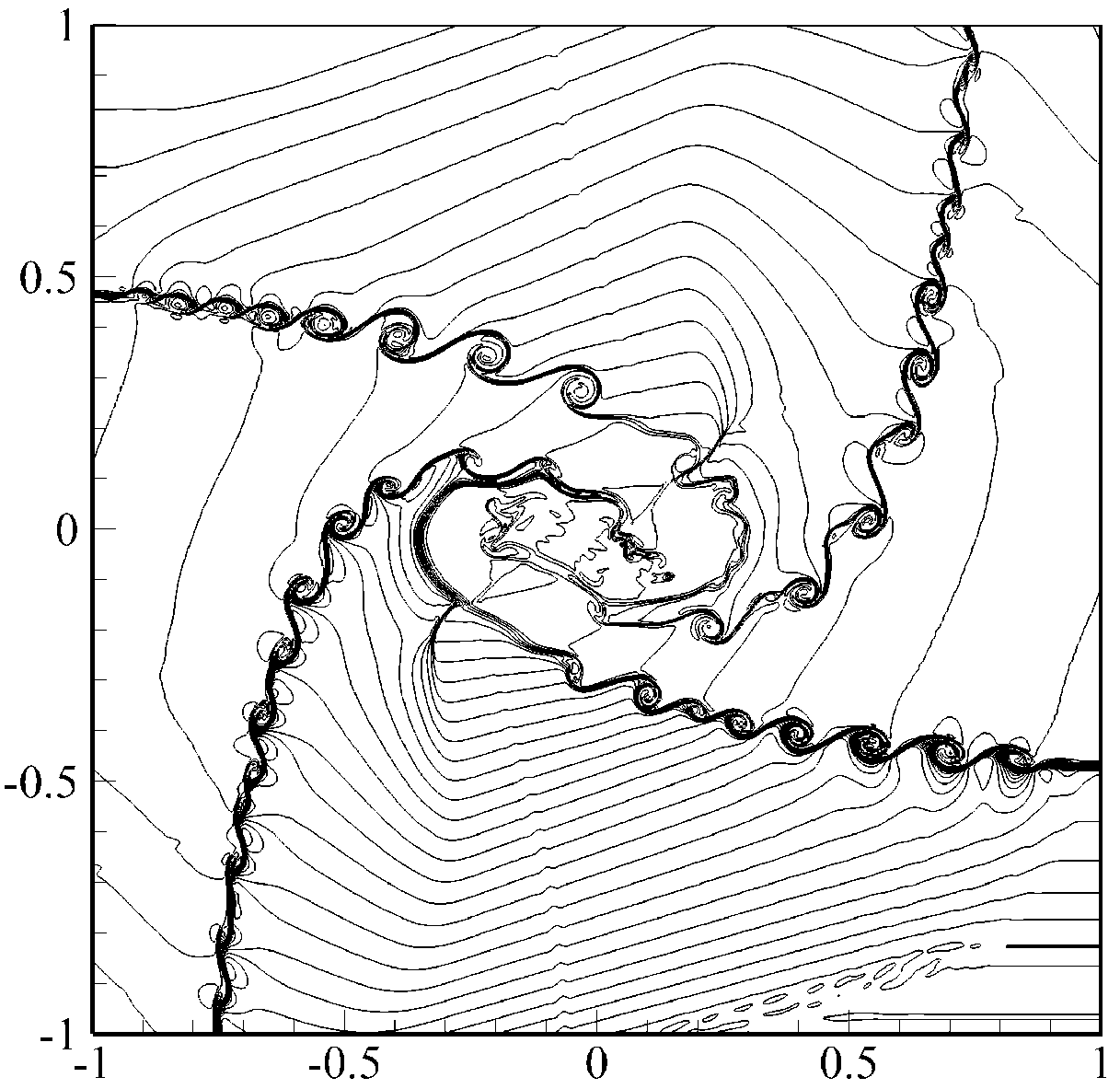}}
  \subfigure[ENO-MR9]{
  \label{FIG:RP2-ENO-MR9}
  \includegraphics[width=5.5 cm]{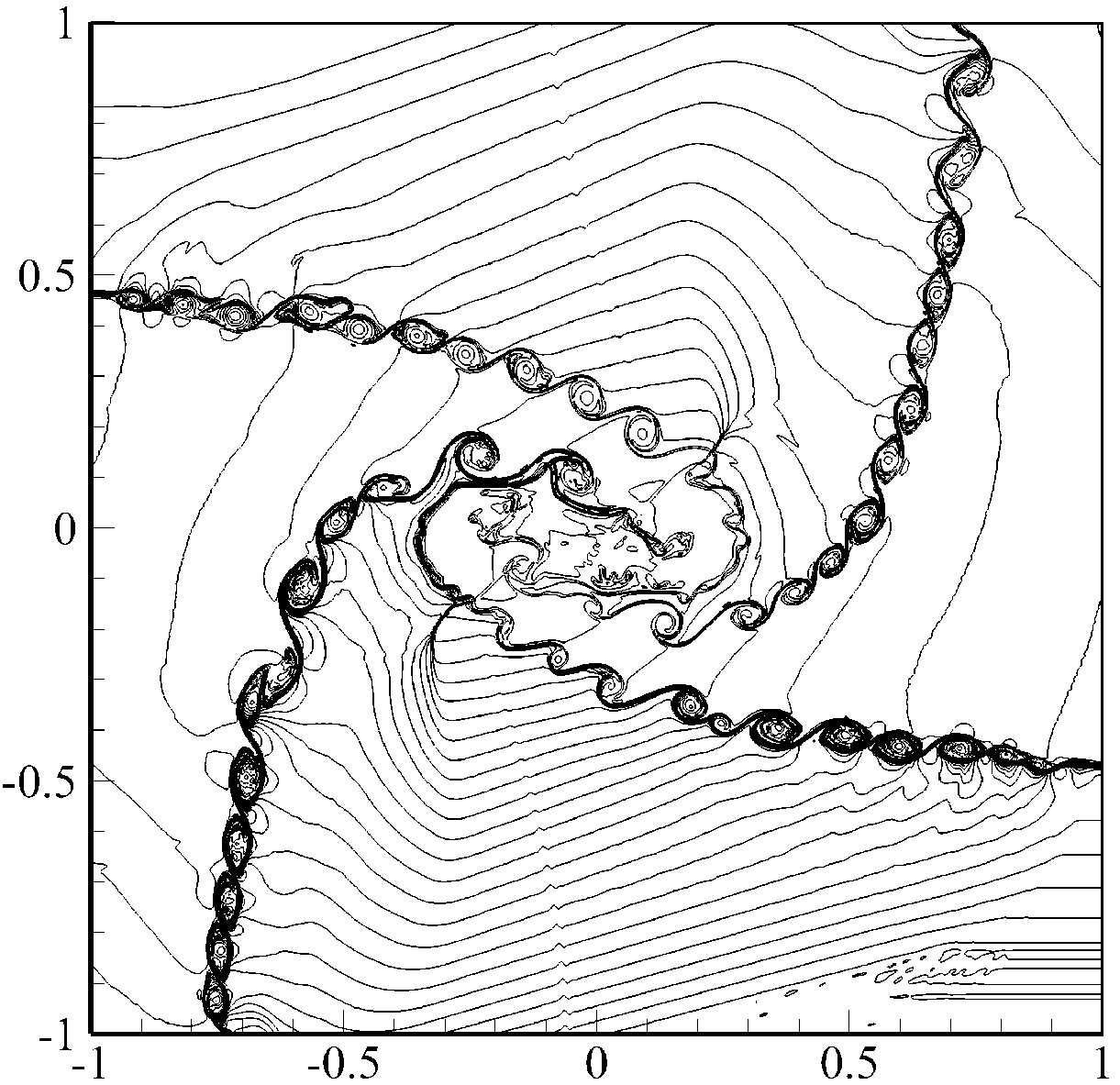}}
  \subfigure[ENO-MR13]{
  \label{FIG:RP2-ENO-MR13}
  \includegraphics[width=5.5 cm]{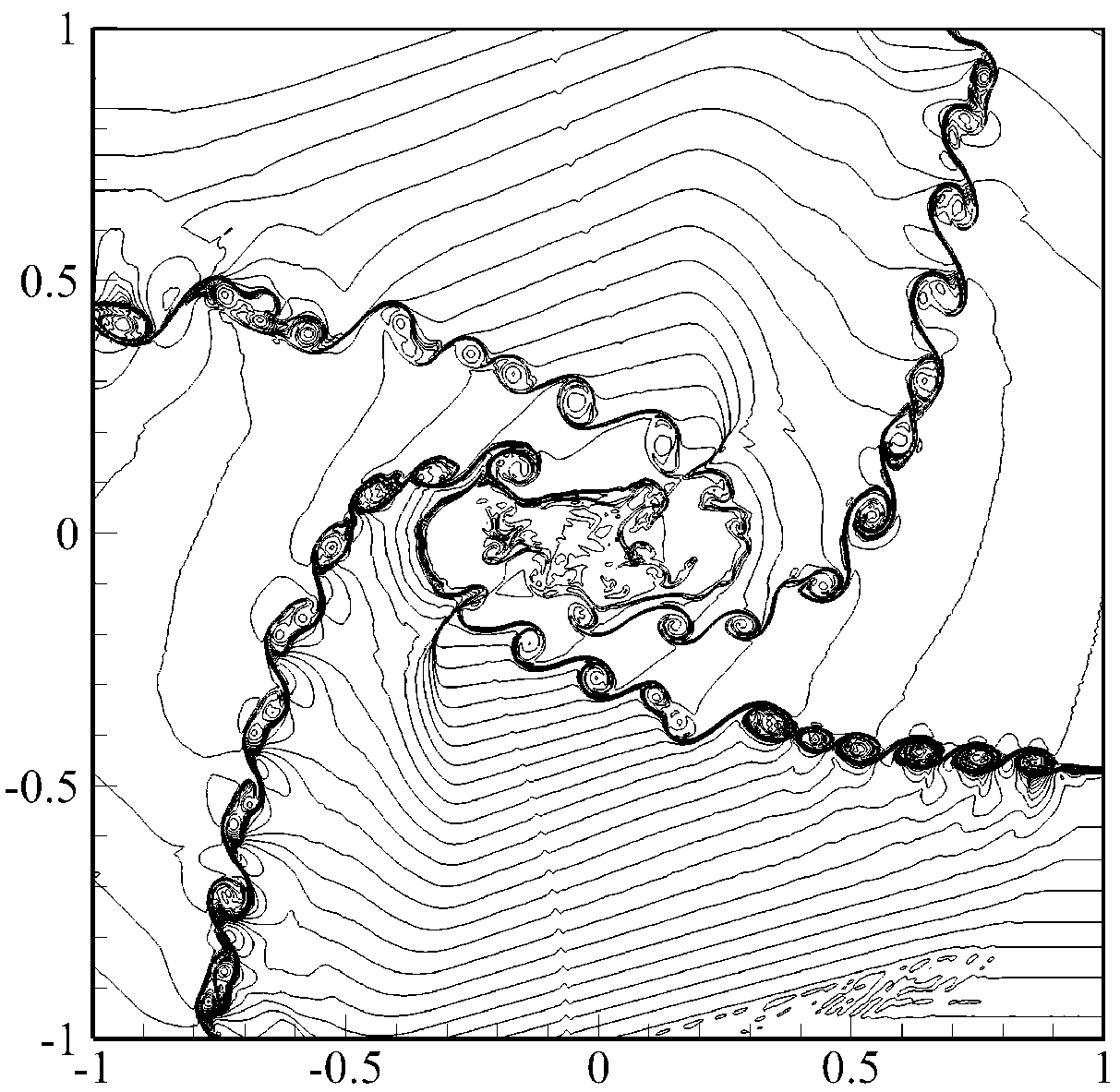}}
  \subfigure[ENO-MR17]{
  \label{FIG:RP2-ENO-MR17}
  \includegraphics[width=5.5 cm]{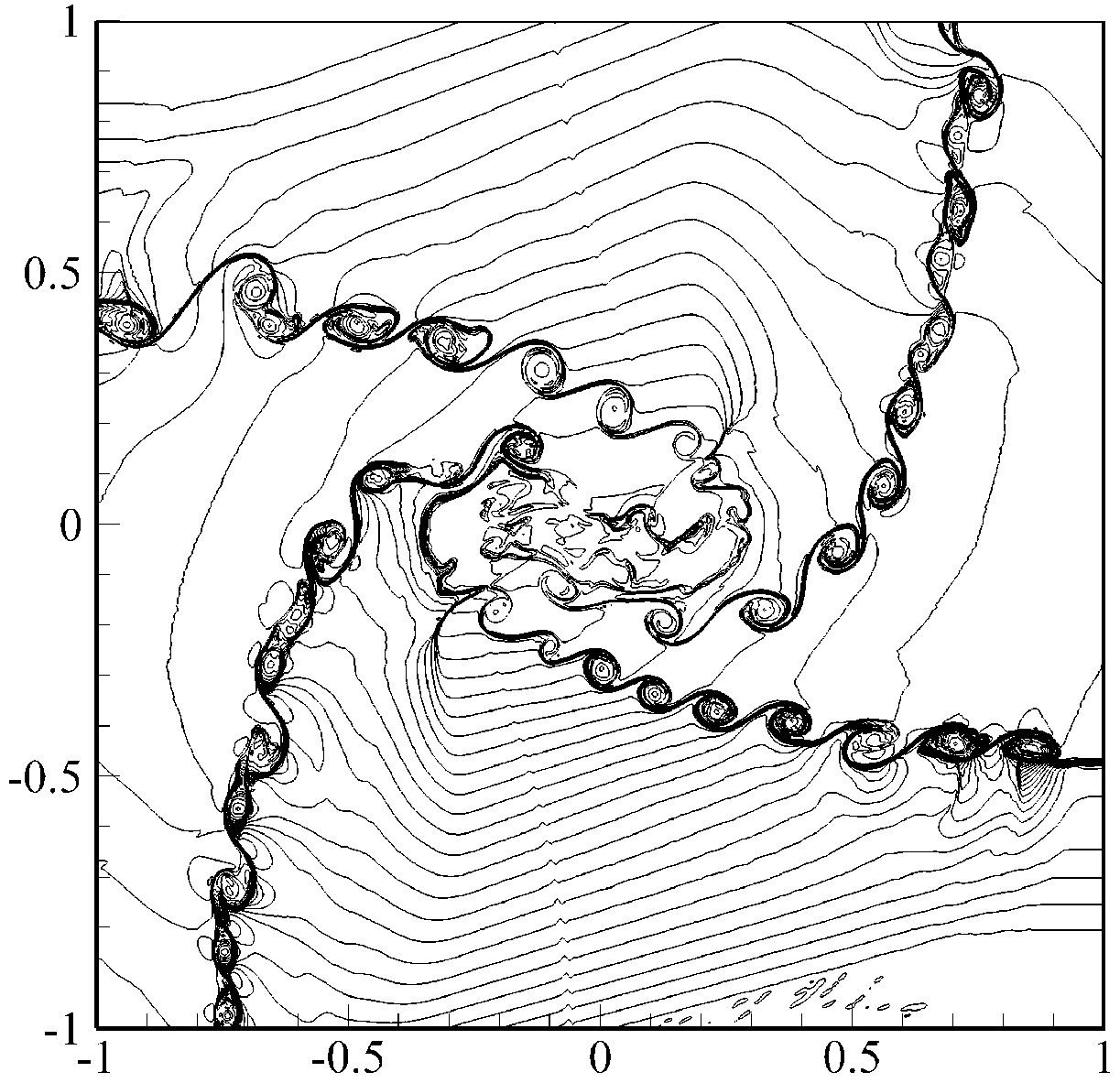}}
  
  \caption{Density contours of the second 2D Riemann problem at $t=1$ calculated 
  by WENO-AO and ENO-MR schemes with $801\times801$ mesh points.
  The contours contain 30 equidistant contours from 0.2 to 3.}
 \label{FIG:RP2}
 \end{figure}

The second configuration starts with four contact discontinuities formulated by
\begin{equation*}
 \left(\rho, u, v, p\right)=
 \begin{cases}
   \left(1,-0.75,0.5,1\right),        & \text{in }[-1, 0]\times[-1, 0],\\
   \left(2,0.75,0.5,1\right),         & \text{in }[-1, 0]\times[0, 1],\\
   \left(1,0.75,-0.5,1\right),        & \text{in }[0, 1] \times[0, 1],\\
   \left(3,-0.75,-0.5,1\right),       & \text{in }[0, 1] \times[-1, 0].
 \end{cases}
\end{equation*}
Fig. \ref{FIG:RP2} shows density contours at $t=1$ 
calculated by WENO-AO and ENO-MR schemes with $801\times801$ mesh points.
WENO-AO and ENO-MR schemes can capture the large-scale vortex at the center of the domain,
but ENO-MR schemes capture many small-scale vortices along contact discontinuities
which are almost absent in the results of WENO-AO schemes at the current resolution.

%%%%%%%%%%%%%%%%%%%%%%%%%%%%%%%%%%%%%%%%%%%%%%%%%%%%%%%%%%%%%%%%%%%%%%%%%%%%%%%%%%%%%%%%%%%%%%%%%%%%
\subsubsection{Double Mach reflection problem}
\begin{figure}[htbp]
  \centering
  \subfigure[WENO-AO(5,3)]{
  \label{FIG:DMR_Density_WENO_AO_5_3}
  \includegraphics[width=5.5 cm]{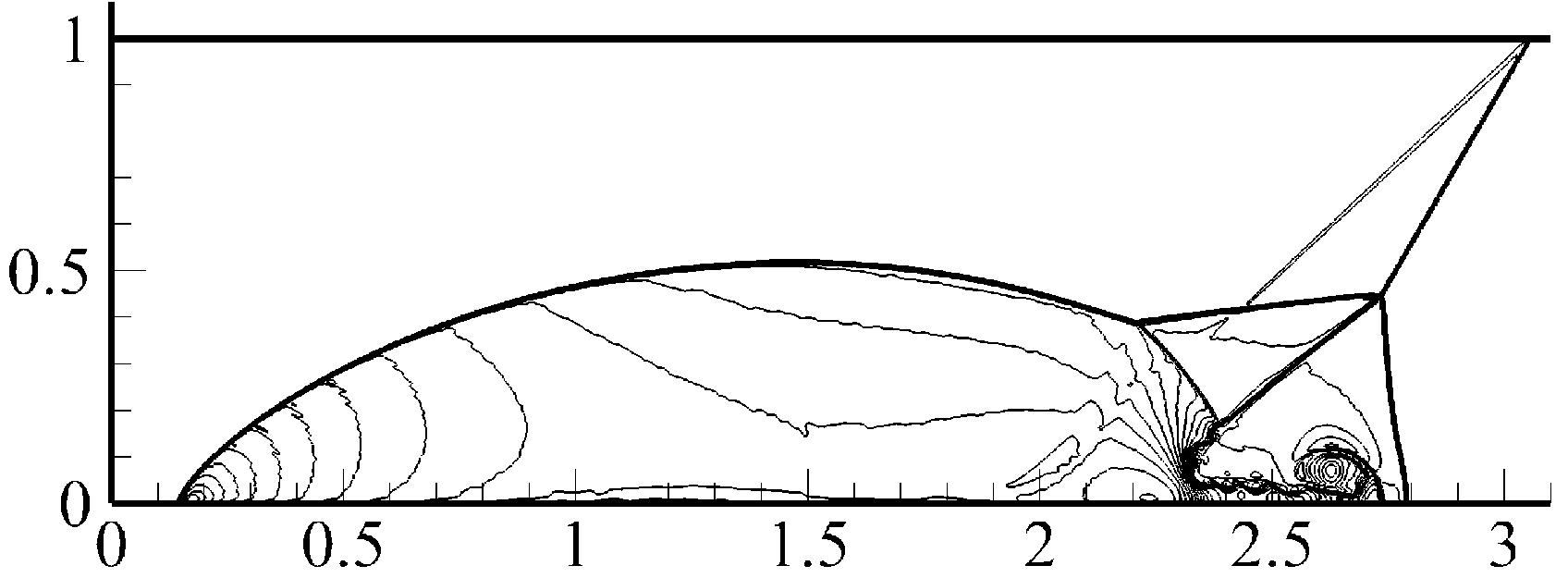}}
  \subfigure[WENO-AO(9,5,3)]{
  \label{FIG:DMR_Density_WENO_AO_9_5_3}
  \includegraphics[width=5.5 cm]{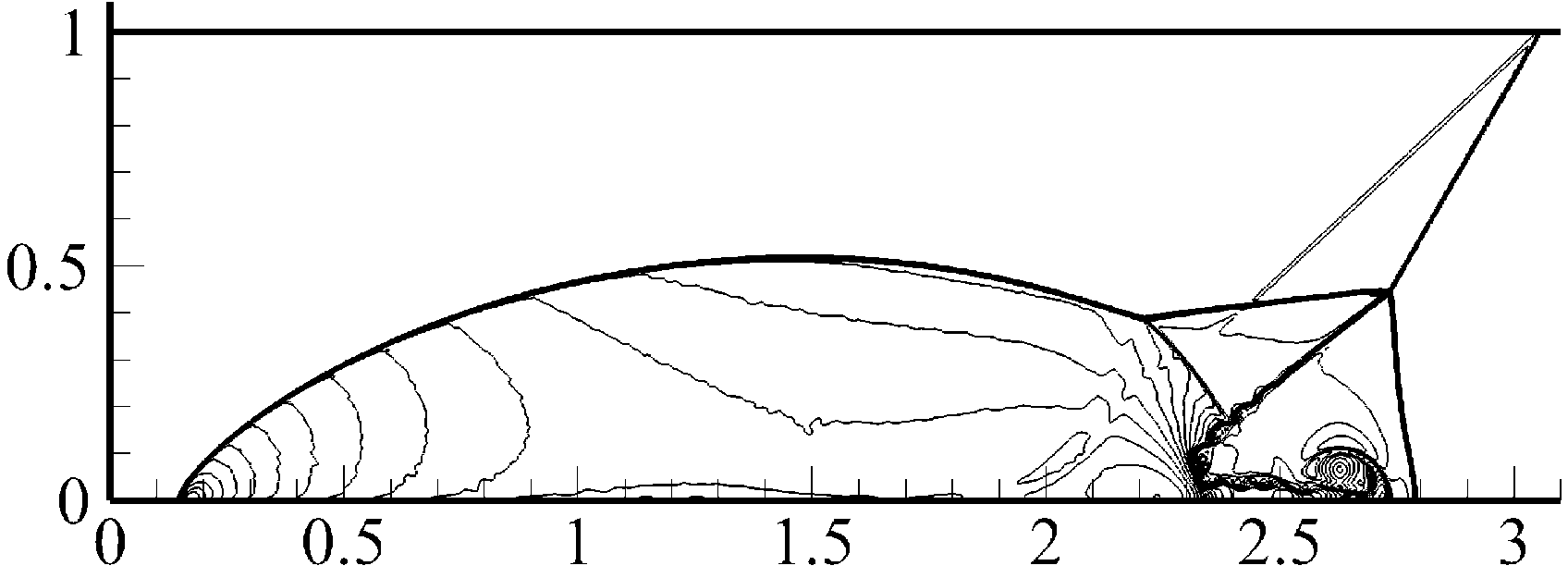}}
  \subfigure[WENO-MR5]{
  \label{FIG:DMR_Density_MR5}
  \includegraphics[width=5.5 cm]{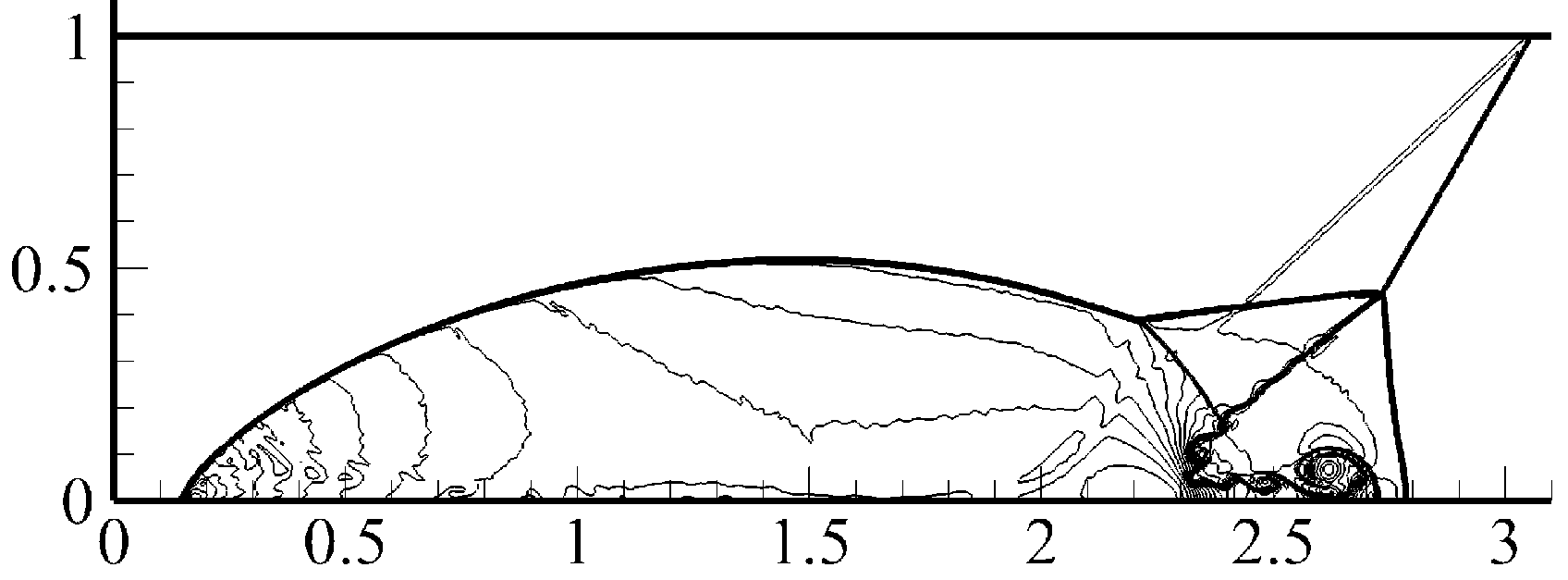}}
  \subfigure[WENO-MR9]{
  \label{FIG:DMR_Density_MR9}
  \includegraphics[width=5.5 cm]{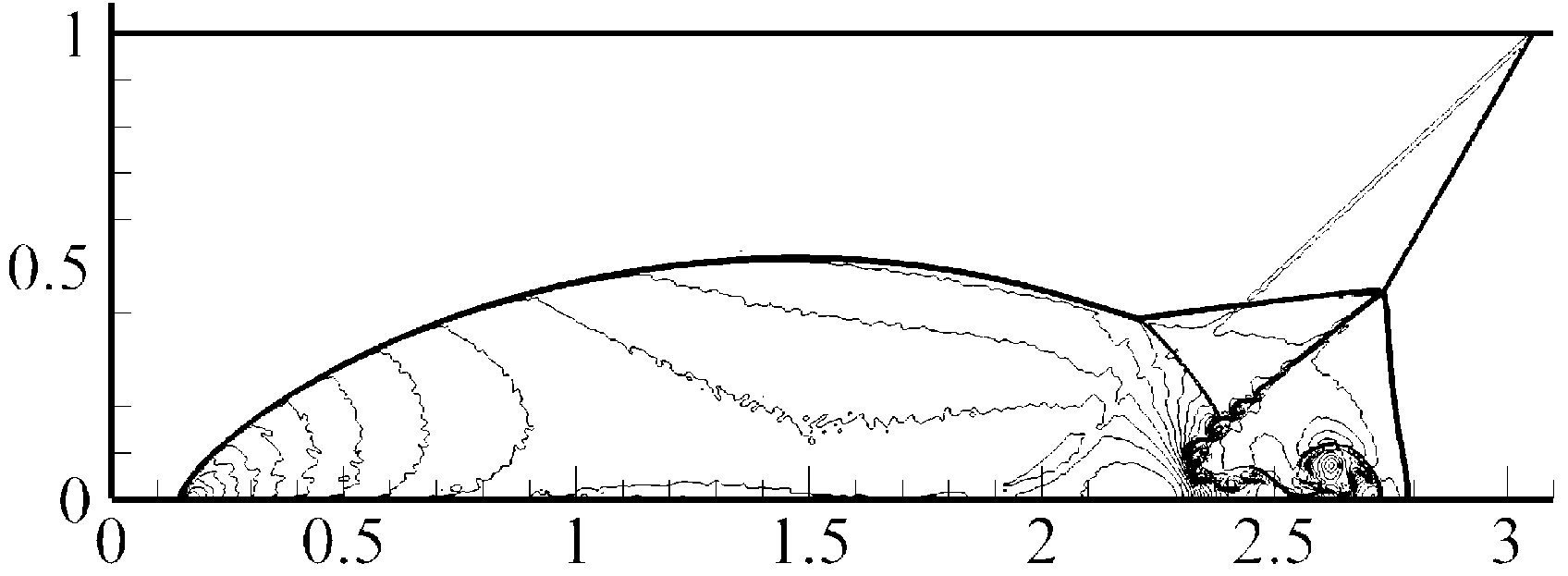}}
  \subfigure[WENO-MR13]{
  \label{FIG:DMR_Density_MR13}
  \includegraphics[width=5.5 cm]{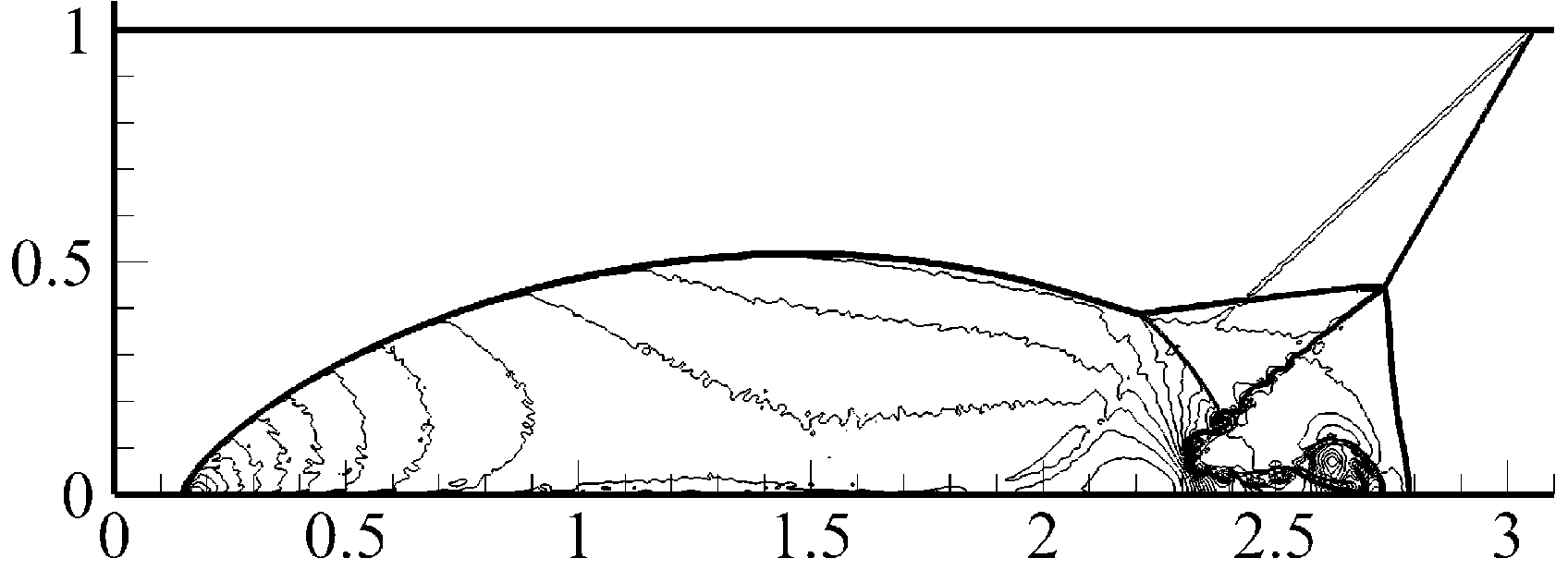}}
  \subfigure[WENO-MR17]{
  \label{FIG:DMR_Density_MR17}
  \includegraphics[width=5.5 cm]{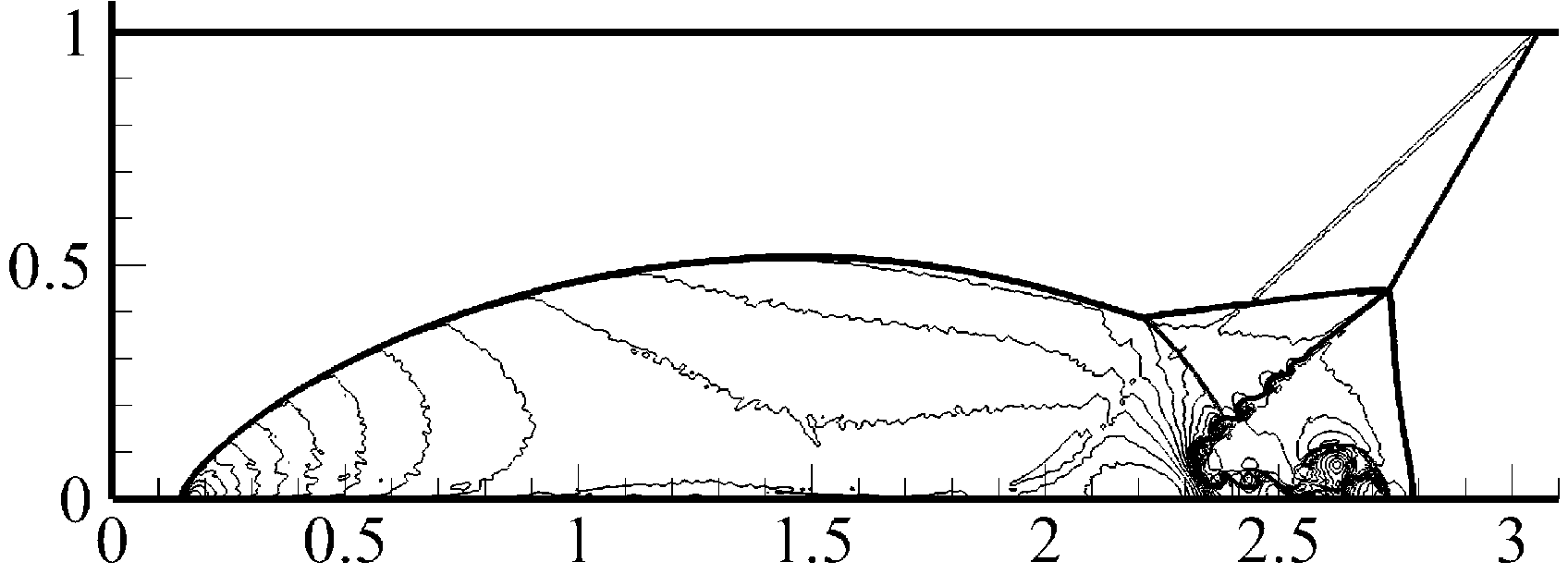}}
  \caption{Density contours of the double Mach reflection problem at $t=0.2$ calculated 
  by WENO-AO and ENO-MR schemes with $1201\times301$ mesh points.
  The density contours contain 40 equidistant contours from 2 to 21.}
 \label{FIG:DMR}
 \end{figure}

\begin{figure}
  \centering
  \subfigure[WENO-AO(5,3)]{
  \label{FIG:DMR_Density_WENO_AO_5_3_enlarge}
  \includegraphics[width=5.5 cm]{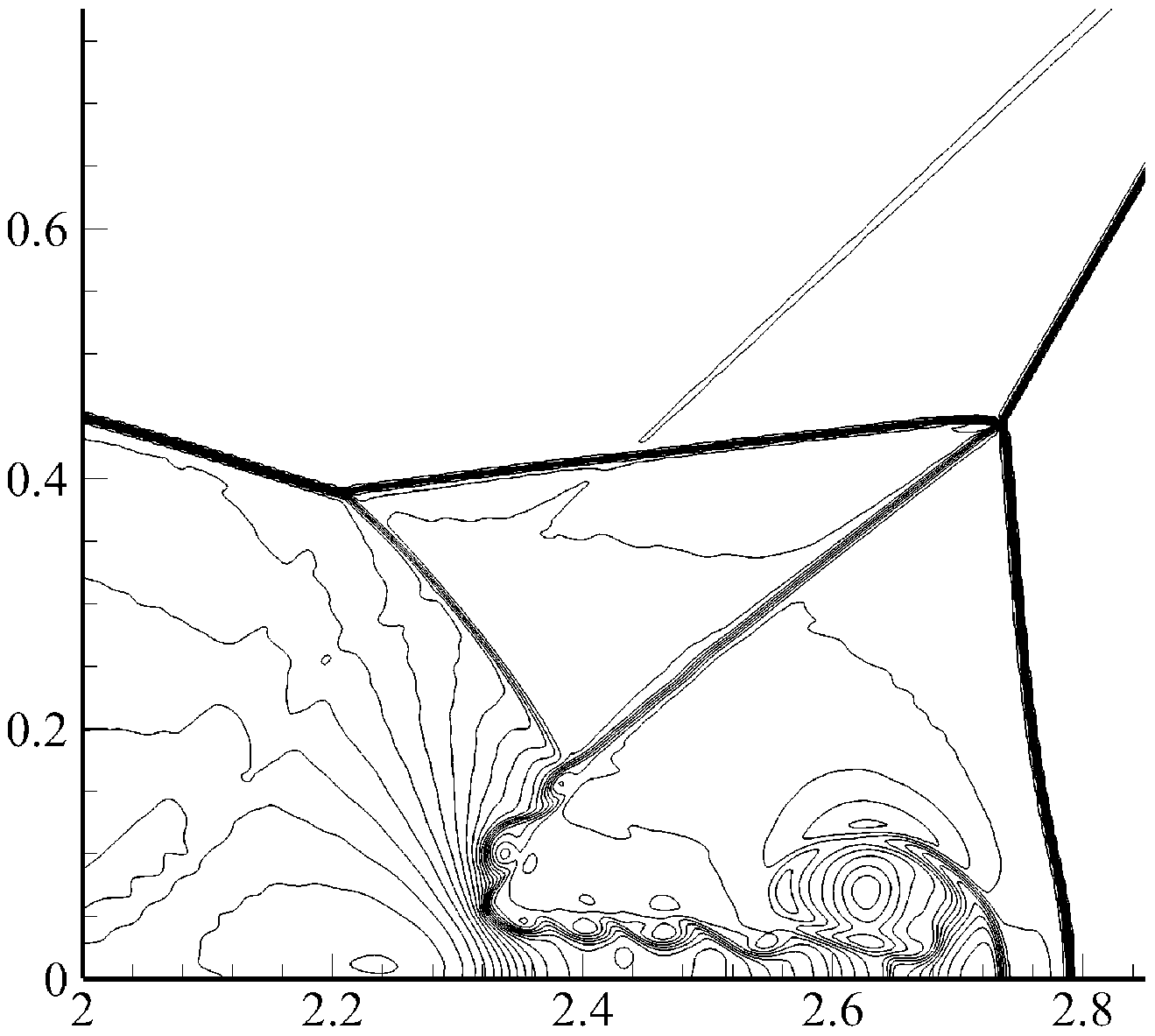}}
  \subfigure[WENO-AO(9,5,3)]{
  \label{FIG:DMR_Density_WENO_AO_9_5_3_enlarge}
  \includegraphics[width=5.5 cm]{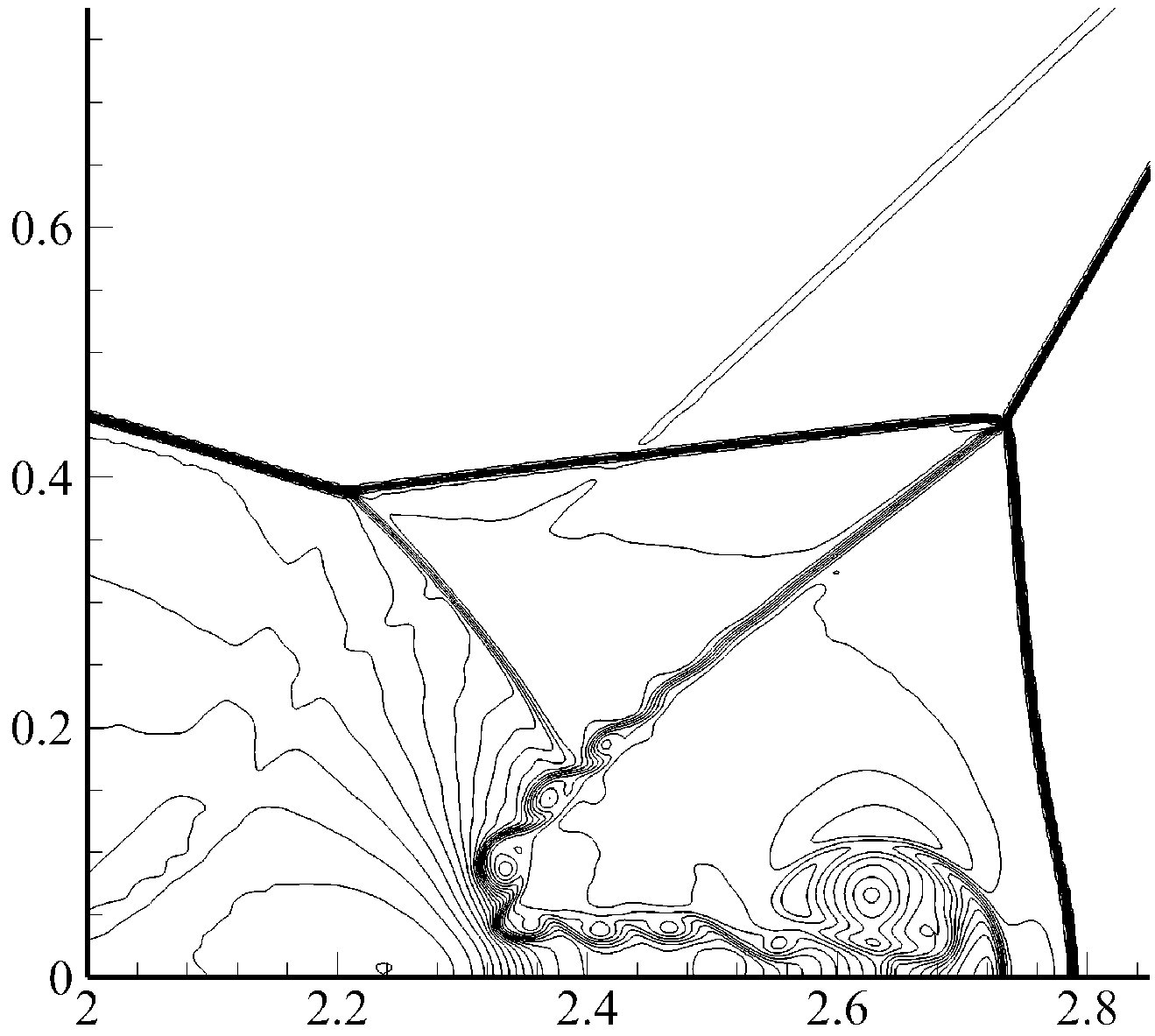}}
  \subfigure[ENO-MR5]{
  \label{FIG:DMR_Density_ENO_MR5_enlarge}
  \includegraphics[width=5.5 cm]{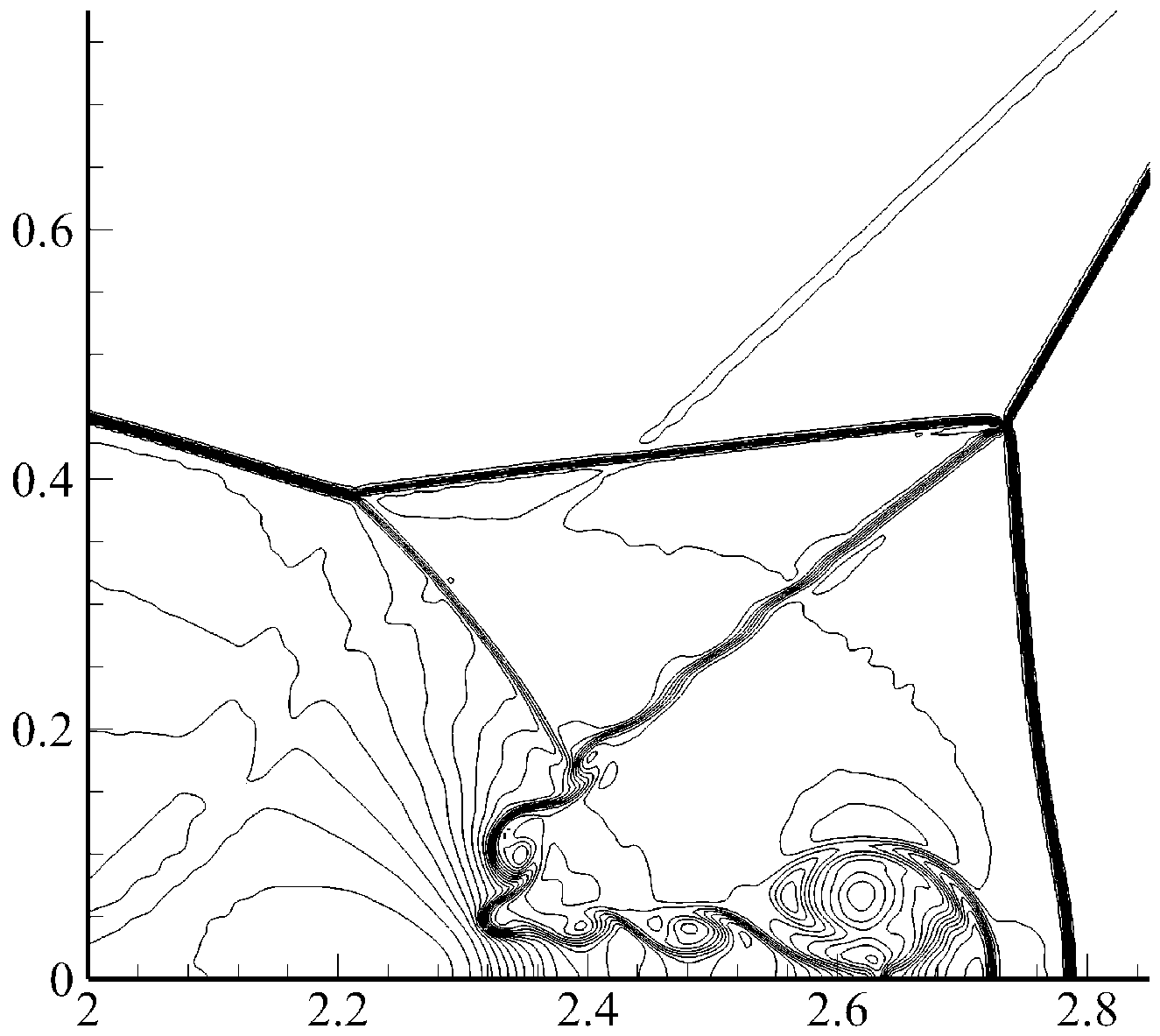}}
  \subfigure[ENO-MR9]{
  \label{FIG:DMR_Density_ENO_MR9_enlarge}
  \includegraphics[width=5.5 cm]{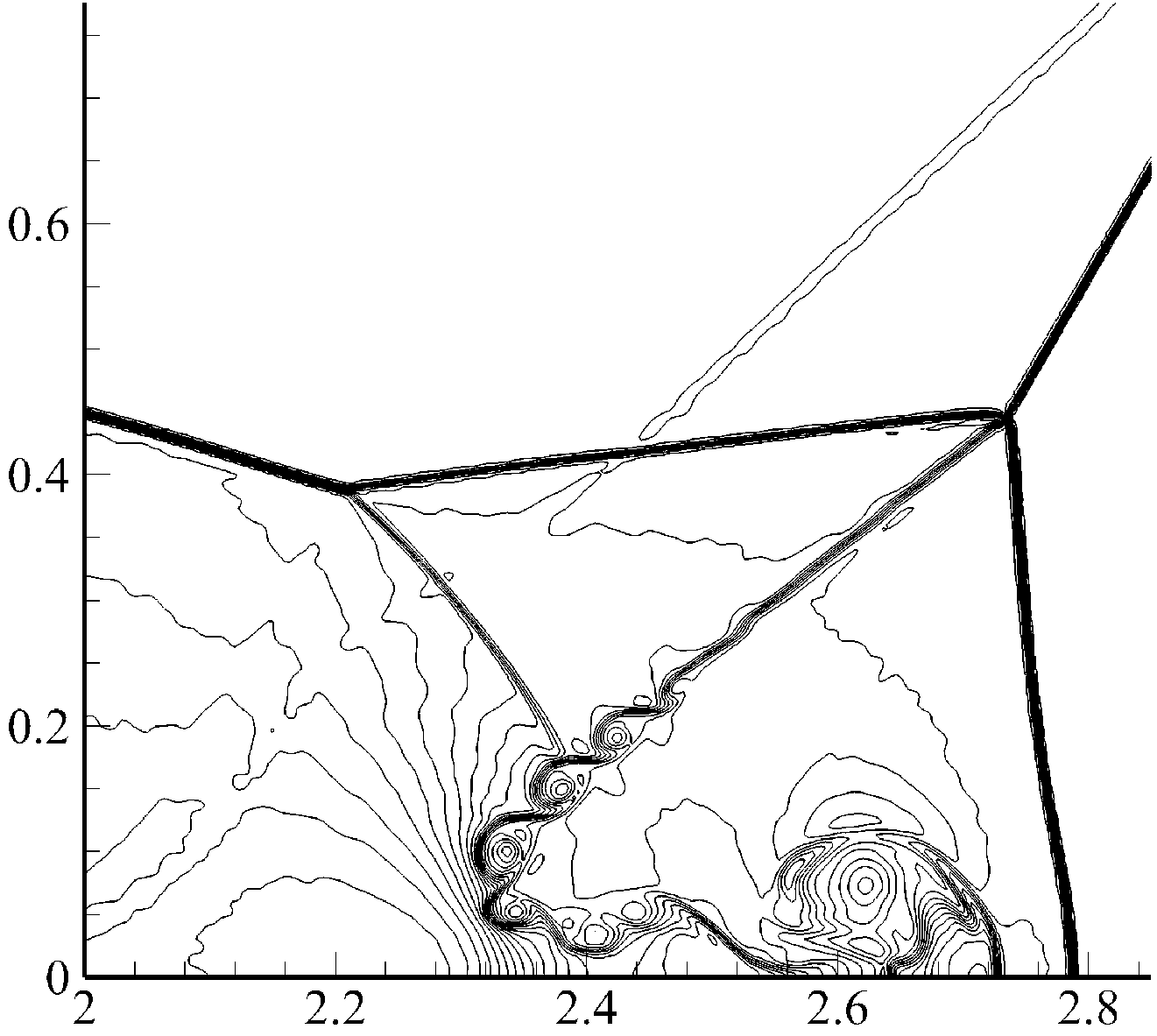}}
  \subfigure[ENO-MR13]{
  \label{FIG:DMR_Density_ENO_MR13_enlarge}
  \includegraphics[width=5.5 cm]{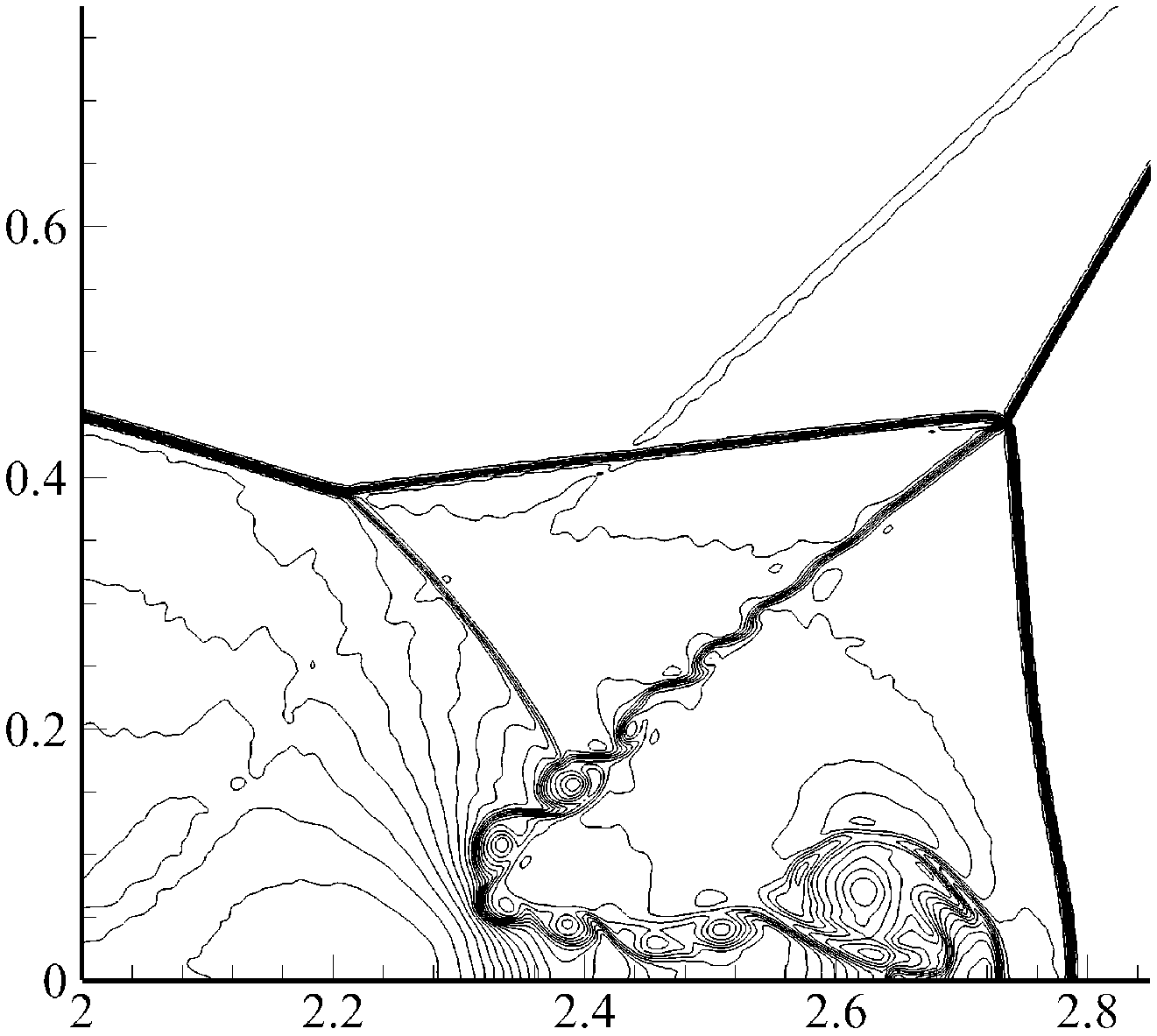}}
  \subfigure[ENO-MR17]{
  \label{FIG:DMR_Density_ENO_MR17_enlarge}
  \includegraphics[width=5.5 cm]{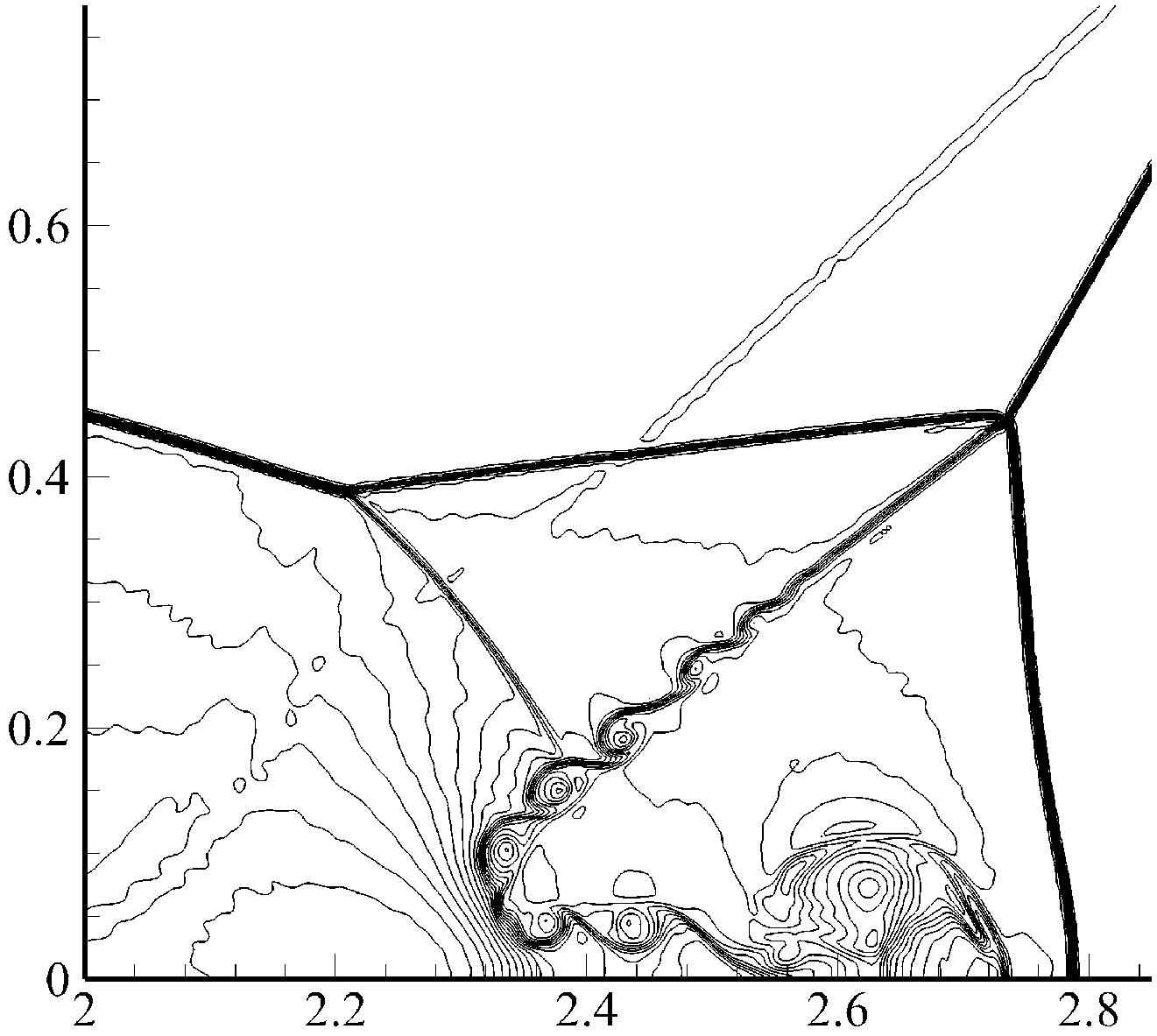}}
  
  \caption{Enlarged view of density contours of the double Mach reflection problem at $t=0.2$ calculated 
  by WENO-AO and ENO-MR schemes with $1201\times301$ mesh points.
  The density contours contain 40 equidistant contours from 2 to 21.}
 \label{FIG:DMR_enlarge}
 \end{figure}

This case simulates the reflection of 
an oblique strong shock impacting on a plate which
was originally proposed by Woodward and Colella \cite{Woodward1984JCP}
and became a benchmark to test the shock-capturing capability and high fidelity of high-order schemes.
The computational domain is $[0, 4]\times[0, 1]$, 
and the initial condition is given as
\begin{equation*}
 \left(\rho, u, v, p\right)=
   \begin{cases}
     \left(1.4, 0, 0, 1\right) & \text{if } x>\frac{1}{6}+\frac{y}{\sqrt{3}}, \\
     \left(8, 8.25sin(60^\circ), -8.25cos(60^\circ), 116.5\right), & \text{otherwise},
   \end{cases}
\end{equation*}
which describes a Mach 10 oblique shock 
inclined at an angle of $60^\circ$ to the horizontal direction.
The post-shocked states are imposed on the left boundary;
Nonreflective boundary conditions are implemented on the right boundary;
The exact motion of the oblique shock is imposed on the top boundary;
Nonreflective and reflective boundary conditions are respectively imposed 
on the bottom boundary for $x\le\frac{1}{6}$ and $x>\frac{1}{6}$.
Fig. \ref{FIG:DMR} and Fig. \ref{FIG:DMR_enlarge} respectively show 
the entire view and the enlarged view of the density contours at $t=0.2$ 
calculated by WENO-AO and ENO-MR schemes with $1201\times301$ mesh points.
All schemes capture the incident shock and reflected shock very well,
but ENO-MR schemes capture more details of the small vortices in the reflection zone.

%%%%%%%%%%%%%%%%%%%%%%%%%%%%%%%%%%%%%%%%%%%%%%%%%%%%%%%%%%%%%%%%%%%%%%%%%%%%%%%%%%%%%%%%%%%%%%%%%%%%
\subsubsection{Rayleigh-Taylor instability}
\begin{figure}
  \centering
  \subfigure[]{
  \label{FIG:RTI128_Density_WENO_AO_5_3}
  \includegraphics[height=5 cm]{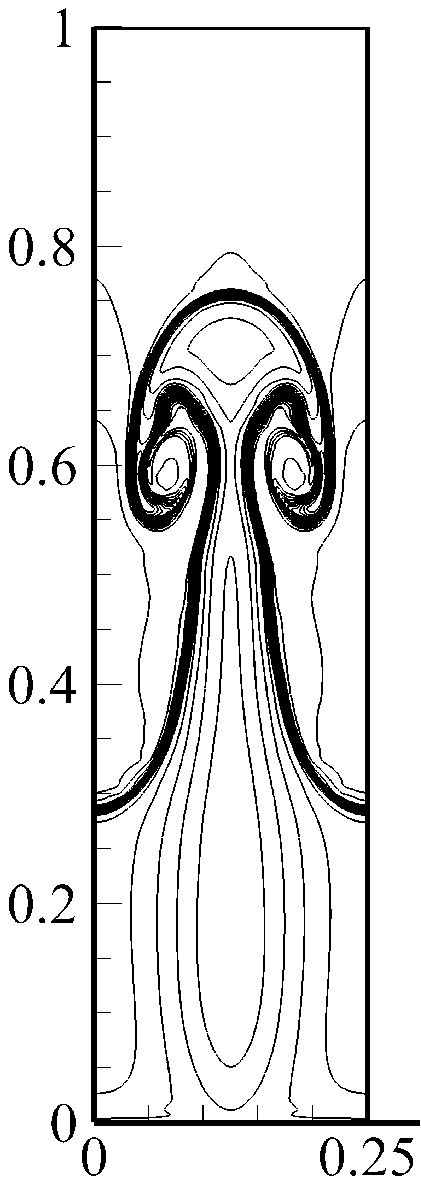}}
  \subfigure[]{
  \label{FIG:RTI128_Density_WENO_AO_9_5_3}
  \includegraphics[height=5 cm]{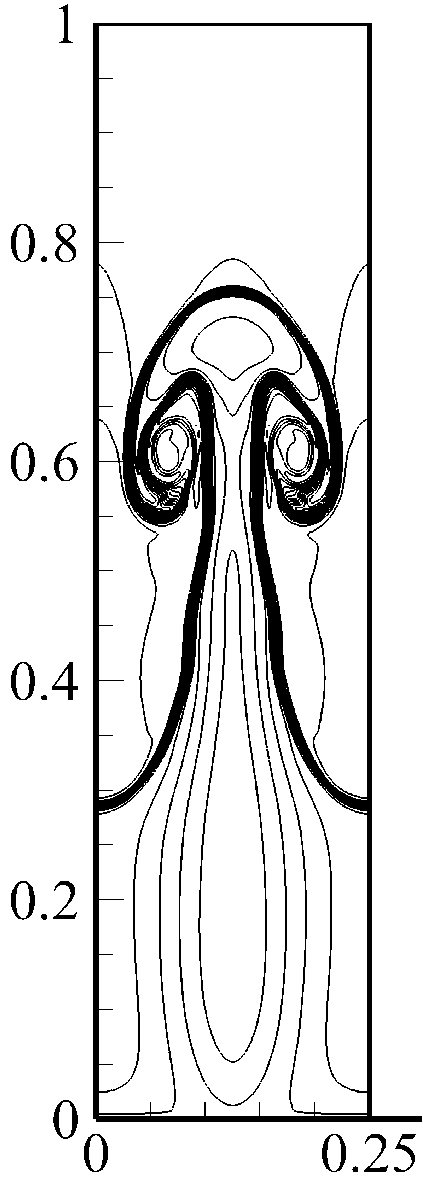}}
  \subfigure[]{
  \label{FIG:RTI128_Density_ENO_MR5}
  \includegraphics[height=5 cm]{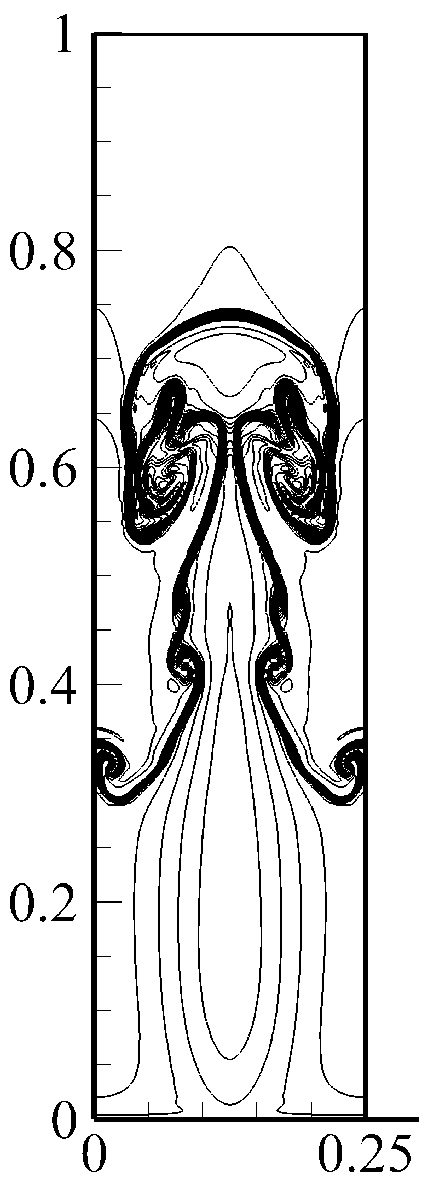}}
  \subfigure[]{
  \label{FIG:RTI128_Density_ENO_MR9}
  \includegraphics[height=5 cm]{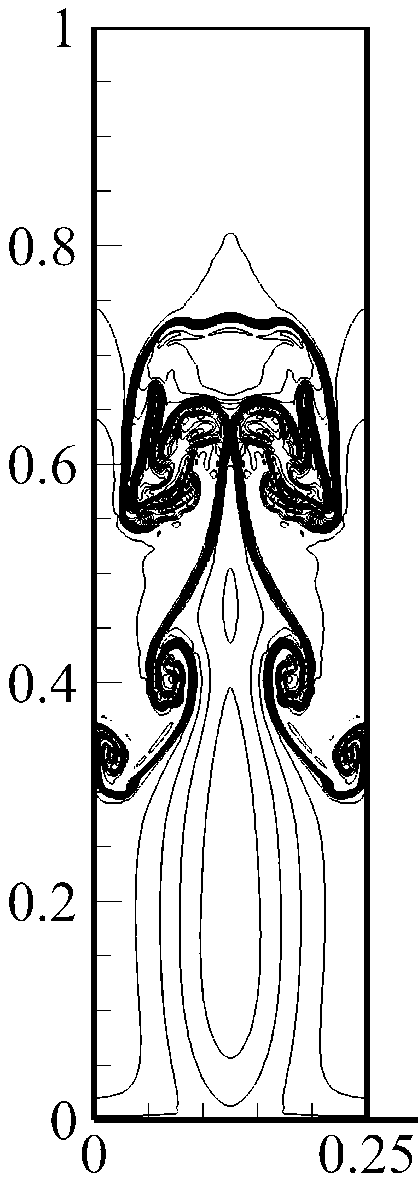}}
  \subfigure[]{
  \label{FIG:RTI128_Density_ENO_MR13}
  \includegraphics[height=5 cm]{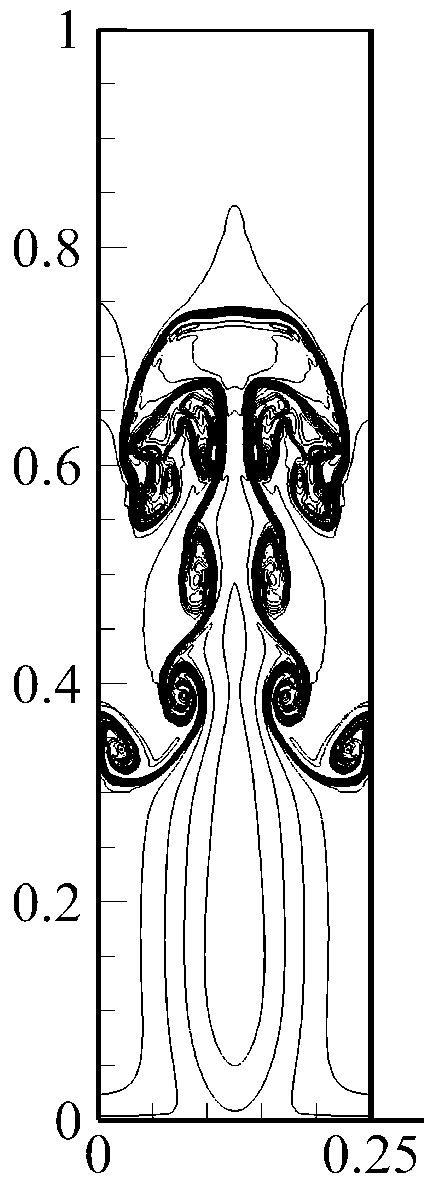}}
  \subfigure[]{
  \label{FIG:RTI128_Density_ENO_MR17}
  \includegraphics[height=5 cm]{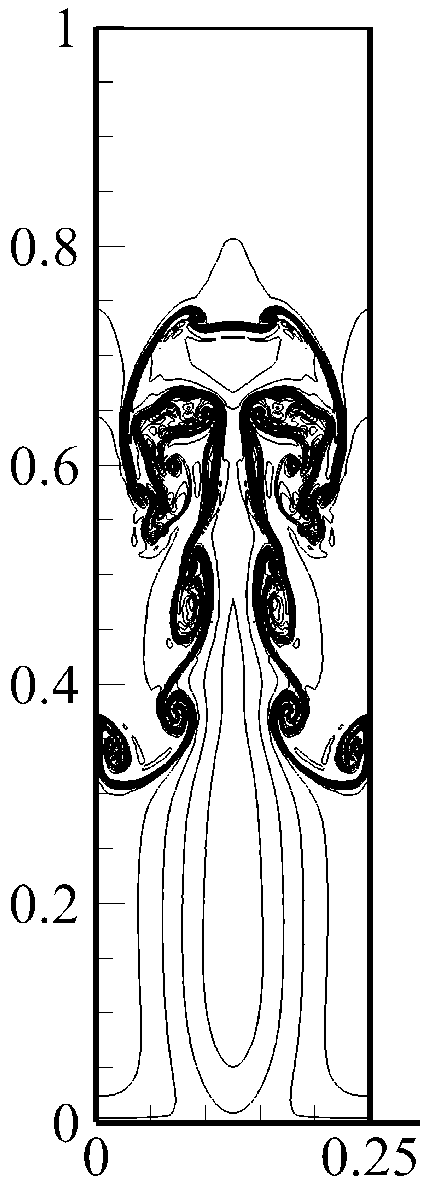}}

  \subfigure[]{
  \label{FIG:RTI200_Density_WENO_AO_5_3}
  \includegraphics[height=5 cm]{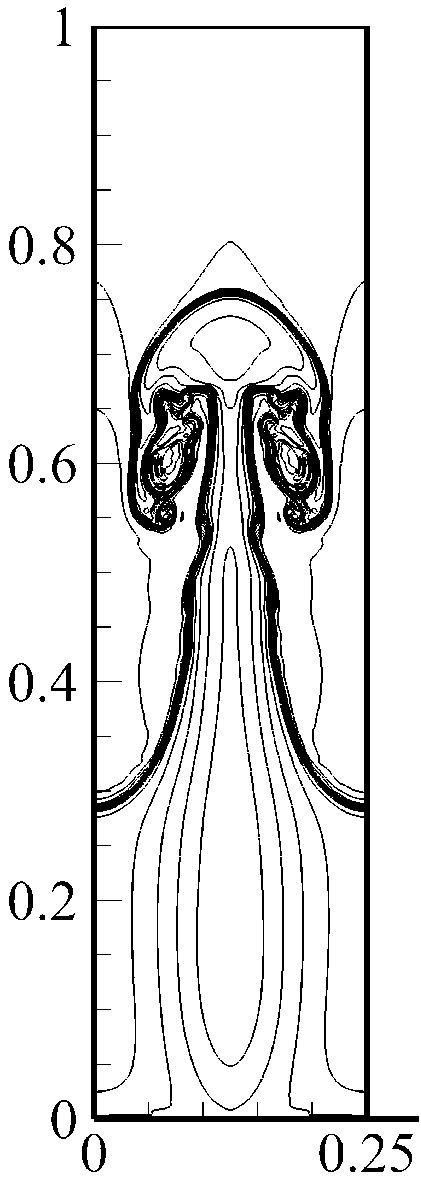}}
  \subfigure[]{
  \label{FIG:RTI200_Density_WENO_AO_9_5_3}
  \includegraphics[height=5 cm]{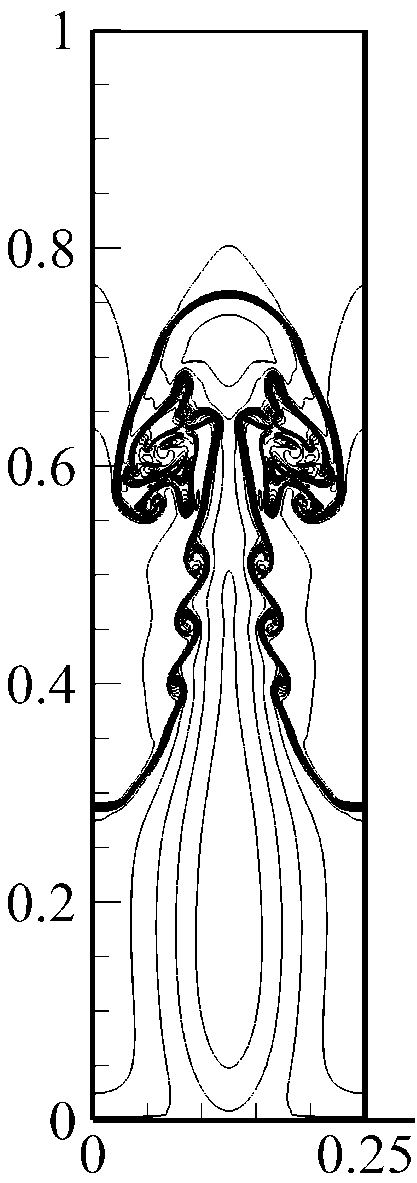}}
  \subfigure[]{
  \label{FIG:RTI200_Density_ENO_MR5}
  \includegraphics[height=5 cm]{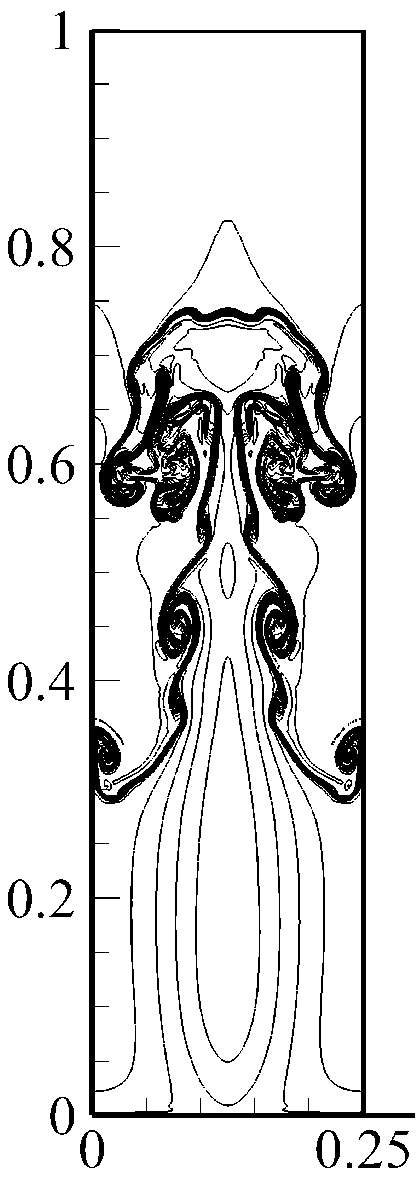}}
  \subfigure[]{
  \label{FIG:RTI200_Density_ENO_MR9}
  \includegraphics[height=5 cm]{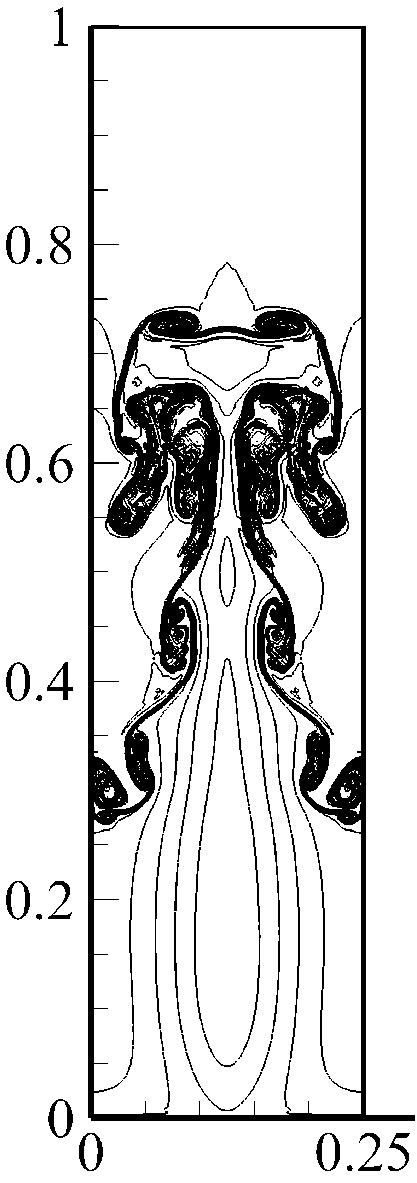}}
  \subfigure[]{
  \label{FIG:RTI200_Density_ENO_MR13}
  \includegraphics[height=5 cm]{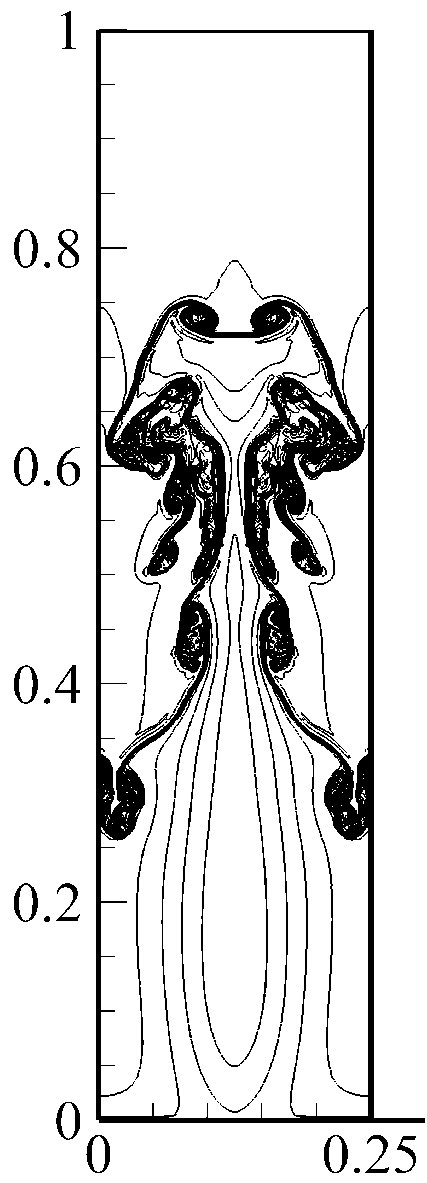}}
  \subfigure[]{
  \label{FIG:RTI200_Density_ENO_MR17}
  \includegraphics[height=5 cm]{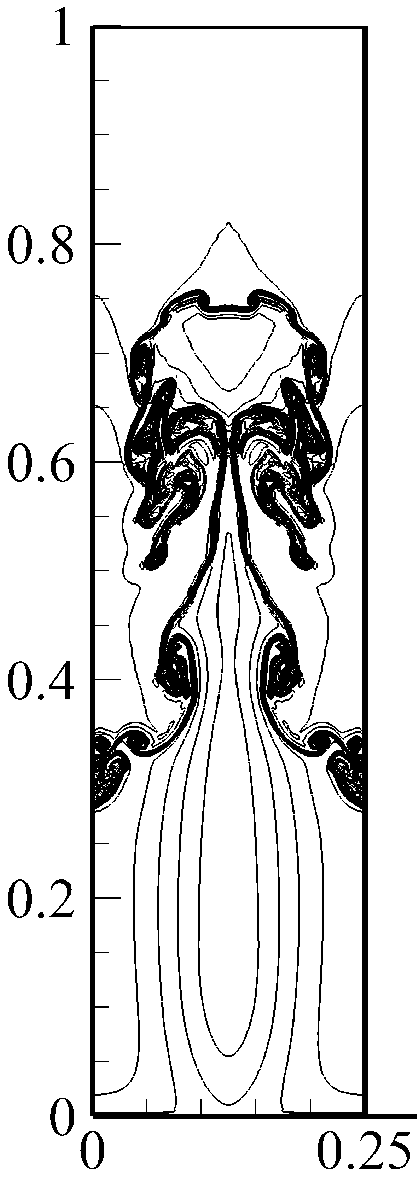}}

  \caption{Density contours of the Rayleigh-Taylor instability problem at $t=1.95$ calculated 
  by WENO-AO and ENO-MR schemes with $129\times 513$ (top row) and $201\times 801$ mesh points (bottom row).
  The results of WENO-AO(5,3), WENO-AO(9,5,3), ENO-MR5, ENO-MR9, ENO-MR13, ENO-MR17 schemes
  are arranged from left to right.
  The density contours contain 20 equidistant contours from 0.9 to 2.2.}
\label{FIG:RTI}
\end{figure}

Rayleigh-Taylor instability problem is another benchmark 
that was widely used to test the high-fidelity properties of high-order numerical schemes.
Following the setup of Shi \emph{et al.} \cite{Shi2003JCP},
the source term $\mathbf{S}=[0,0,\rho,\rho v]^T$ is added to 
the right hand side of the 2D Euler euqations, Eq. (\ref{Eq:2D_Euler_equations}).
The computational domain is $[0,0.25]\times[0,1]$, the ratio of specific heats $\gamma=5/3$,
and the initial condition is given by
\begin{equation*}
  \left(\rho, u, v, p\right)=
    \begin{cases}
      \left(2, 0, -0.025\sqrt{\frac{\gamma p}{\rho}}cos(8\pi x), 2y+1\right) & \text{if } y<0.5, \\
      \left(1, 0, -0.025\sqrt{\frac{\gamma p}{\rho}}cos(8\pi x), y+1.5\right), & \text{otherwise}.
    \end{cases}
\end{equation*}
The development of Rayleigh-Taylor instabilities will induce complex fingering structures 
which are sensitive to small numerical errors.
Low-dissipation schemes easily lose symmetry
which can be fixed by carefully implementing the computer programs \cite{Fleischmann2019symmetry}. 
In this study, we enforce the symmetry by a simple strategy from \cite{Wang2020AWENO}.
Fig. \ref{FIG:RTI} shows the density contours at $t=1.95$ calculated 
by WENO-AO and ENO-MR schemes with $129\times 513$ and $201\times 801$ mesh points respectively.
We observe that ENO-MR schemes capture much more details than WENO-AO schemes under the same grid resolution.
High-order ENO-MR schemes also perform better than the ENO-MR5 scheme in this case.

%%%%%%%%%%%%%%%%%%%%%%%%%%%%%%%%%%%%%%%%%%%%%%%%%%%%%%%%
\subsubsection{Computational costs for two-dimensional tests}

\begin{table}
 \footnotesize
 \caption{Normalized computational costs of WENO-AO and ENO-MR schemes for different cases.}
 \label{Table:Computational_Cost}
 \centering
 \begin{tabular}{|c|c|c|c|c|c|c|}
   \hline
   % after \\: \hline or \cline{col1-col2} \cline{col3-col4} ...
     test       & WENO-AO(5,3) & WENO-AO(9,5,3)  &ENO-MR5 &ENO-MR9 &ENO-MR13 &ENO-MR17    \\
     \hline
      RP1       &1             &1.86     &0.62            &1.09   &1.62   &2.21\\
      RP2       &1             &1.90     &0.62            &1.05   &1.62   &2.17\\
      DMR       &1             &1.86     &0.58            &1.02   &1.59   &2.30\\
      RTI128    &1             &1.82     &0.55            &0.84   &1.20   &1.68\\
      RTI200    &1             &1.87     &0.55            &0.81   &1.35   &1.72\\
     \hline
 \end{tabular}
\end{table}

To roughly test the efficiency of the proposed ENO-MR schemes,
we record the normalized computational costs of WENO-AO and ENO-MR schemes
for the first Riemann problem (RP1), the second Riemann problem (RP2), the double Mach reflection (DMR),
and the Rayleigh-Taylor instability (RTI) in Table \ref{Table:Computational_Cost}.
The computational costs are empirically determined by the parallel runtime of C codes 
using 16 CPU cores.
We observe that ENO-MR schemes are more efficient than WENO-AO schemes of the same order.
The computational costs of ENO-MR9 are similar to those of WENO-AO(5,3),
although ENO-MR9 has much more candidate stencils.
Even for the very high-order ENO-MR schemes,
the computational costs are still affordable. 
The high efficiency of ENO-MR schemes benefits from two aspects:
First, smoothness indicators of ENO-MR schemes are simple;
Second, we directly use high-order stencils to reconstruct fluxes in smooth regions,
and do not need to calculate smoothness indicators of remaining stencils.
The second aspect is also why the relative computational costs of ENO-MR schemes have 
obvious fluctuations, because the proportion of smooth regions is different for different problems.  

%%%%%%%%%%%%%%%%%%%%%%%%%%%%%%%%%%%%%%%%%%%%%%%%%%%%%%%%%%%%%%%%%%%%%%%%%%%%%%%%%%%%%%%%
\section{Conclusions}
We construct a class of ENO-MR schemes with increasingly high-order
based on linearly stable candidate stencils of unequal sizes
by using simple smoothness indicators and an efficient selection strategy.
The candidate stencils range from first-order up to the designed very high-order,
and the proposed ENO-MR schemes can adaptively select the optimal stencil
according to the simple smoothness indicators.
Analysis and plenty of numerical examples show that
the constructed ENO-MR schemes achieve optimal order in smooth regions 
while maintaining oscillation-free near strong discontinuities.
Additionally, ENO-MR schemes are more accurate and more efficient than WENO-AO schemes,
particularly for the instabilities of contact discontinuities.
In conclusion, ENO-MR schemes provide an efficient way to 
construct high-order schemes with multi-resolution.
%%%% Acknowledgments %%%%%%%%
\section*{Acknowledgments}
The author would like to acknowledge the financial support of 
the National Natural Science Foundation of China (Contract Nos. 11901602 and 62231016).

\bibliographystyle{unsrtnat}
\bibliography{references}  %%% Uncomment this line and comment out the ``thebibliography'' section below to use the external .bib file (using bibtex) .

%%% Uncomment this section and comment out the \bibliography{references} line above to use inline references.
% \begin{thebibliography}{1}

% 	\bibitem{kour2014real}
% 	George Kour and Raid Saabne.
% 	\newblock Real-time segmentation of on-line handwritten arabic script.
% 	\newblock In {\em Frontiers in Handwriting Recognition (ICFHR), 2014 14th
% 			International Conference on}, pages 417--422. IEEE, 2014.

% 	\bibitem{kour2014fast}
% 	George Kour and Raid Saabne.
% 	\newblock Fast classification of handwritten on-line arabic characters.
% 	\newblock In {\em Soft Computing and Pattern Recognition (SoCPaR), 2014 6th
% 			International Conference of}, pages 312--318. IEEE, 2014.

% 	\bibitem{hadash2018estimate}
% 	Guy Hadash, Einat Kermany, Boaz Carmeli, Ofer Lavi, George Kour, and Alon
% 	Jacovi.
% 	\newblock Estimate and replace: A novel approach to integrating deep neural
% 	networks with existing applications.
% 	\newblock {\em arXiv preprint arXiv:1804.09028}, 2018.

% \end{thebibliography}

\end{document}